\documentclass[11pt,a4paper]{article}

\usepackage[T1]{fontenc}
\usepackage{lmodern}
\usepackage{microtype}
\usepackage{fullpage}

\usepackage{mathtools}            
\usepackage{amssymb}
\usepackage{amsthm}
\usepackage{empheq}

\usepackage[table,dvipsnames]{xcolor}
\usepackage{graphicx}
\usepackage{float}
\usepackage{subcaption}           
\usepackage{pifont}
\usepackage{tikz}
\usetikzlibrary{arrows.meta,calc}

\definecolor{myred}{rgb}{0.8,0,0}
\definecolor{myblue}{rgb}{0,0,0.8}
\definecolor{DeepGreen}{RGB}{92,172,129}
\definecolor{DeepRed}{RGB}{167,67,67}
\definecolor{originalcurve}{RGB}{48,48,48}
\definecolor{originalblue}{RGB}{20,102,192}
\definecolor{originalred}{RGB}{255,51,51}

\usepackage{booktabs}

\providecommand{\ExperimentTableSetup}{%
    \scriptsize
    \setlength{\tabcolsep}{4pt}%
    \renewcommand{\arraystretch}{1.08}%
}

\usepackage{algorithm}
\usepackage{algpseudocode}        

\usepackage{enumerate}
\usepackage{cite}
\usepackage[title]{appendix}

\usepackage{hyperref}
\hypersetup{
    colorlinks=true,
    linkcolor=blue,
    citecolor=blue,
    urlcolor=black
}
\usepackage{cleveref}             

\newcommand{\email}[1]{\href{mailto:#1}{#1}}

\allowdisplaybreaks
\numberwithin{equation}{section}

\renewenvironment{abstract}
    {\par\vspace{0.5\baselineskip}\noindent\textbf{Abstract.}\ \ignorespaces}
    {\par\vspace{0.5\baselineskip}}

\newenvironment{keywords}
    {\par\vspace{0.5\baselineskip}\noindent\textbf{Keywords.}\ \ignorespaces}
    {\par\vspace{0.5\baselineskip}}

\newenvironment{MSCcodes}
    {\par\vspace{0.5\baselineskip}\noindent\textbf{MSC2020.}\ \ignorespaces}
    {\par\vspace{0.5\baselineskip}}

\newtheorem{definition}{Definition}[section]
\newtheorem{theorem}{Theorem}[section]
\newtheorem{lemma}{Lemma}[section]
\newtheorem{proposition}{Proposition}[section]
\newtheorem{corollary}{Corollary}[section]
\newtheorem{assumption}{Assumption}[section]
\newtheorem{remark}{Remark}[section]
\newtheorem{example}{Example}[section]

\newcommand{\mS}{\mathbb{S}}
\newcommand{\Spp}{\mathbb{S}_{++}}
\newcommand{\R}{\mathbb{R}}
\newcommand{\m}{\mathcal{M}}
\newcommand{\HH}{\mathbb{H}^n}
\newcommand{\dd}{\mathrm{d}}

\newcommand{\e}{\mathrm{Exp}}
\newcommand{\Hess}{\operatorname{Hess}}
\newcommand{\grad}{\operatorname{grad}}
\newcommand{\Unif}{\operatorname{Unif}}

\DeclareMathOperator{\T}{\mathrm{T}}

\title{A FISTA-type Riemannian Accelerated Proximal Gradient Method}

\author{%
    Shuailing Feng\thanks{%
        School of Mathematical Sciences, Xiamen University, Xiamen, China;
        School of Mathematics and Computer Science, Yan'an University,
        Yan'an, China (\email{shl.feng@yau.edu.cn}).%
    }
    \and
    Yuhang Jiang\thanks{%
        School of Mathematical Sciences, Xiamen University, Xiamen, China
        (\email{jiangyuhang1@stu.xmu.edu.cn}).%
    }
    \and
    Wen Huang\thanks{%
        Corresponding author.
        School of Mathematical Sciences, Xiamen University, Xiamen, China
        (\email{wen.huang@xmu.edu.cn}).%
    }
    \and
    Shihui Ying\thanks{%
        Shanghai Institute of Applied Mathematics and Mechanics,
        Shanghai University, Shanghai, China;
        School of Mechanics and Engineering Science,
        Shanghai University, Shanghai, China (\email{shying@shu.edu.cn}).%
    }
    \and
    Bisheng Xia\thanks{%
        School of Mathematics and Computer Science, Yan'an University,
        Yan'an, China (\email{bishengxia@163.com}).%
    }
}

\begin{document}
\maketitle

\begin{abstract}
We study accelerated proximal gradient methods for minimizing the sum of a smooth and a nonsmooth function on Riemannian manifolds. 
Accelerated convergence rates under geodesic convexity have  not been established for some existing methods, and global geodesic convexity may fail or be difficult to verify.
We develop a Riemannian accelerated proximal gradient method that combines tangent-space proximal steps with momentum updates based on parallel transport and the exponential map.
Under geodesic $\mu$-strong convexity ($\mu\geqslant0$) of the smooth term and pullback $\rho$-convexity ($\rho\in\mathbb{R}$) of the nonsmooth term, we establish accelerated convergence rates for the objective gap: an $\mathcal{O}(k^{-2})$ rate when $\mu-\rho\geqslant0$ and a linear rate when $\mu-\rho>0$, where $k$ is the iteration count.
The rates depend on sectional curvature bounds and the diameter of a prescribed region  assumed to contain the iterates.
For nonconvex objectives or unknown problem parameters, we further consider an adaptive restart variant with  proximal gradient safeguards. If the initial sublevel set is compact, we prove that every accumulation point of the safeguard reference sequence is stationary and that an $\varepsilon$-stationary point is obtained within $\mathcal{O}(\varepsilon^{-2})$ iterations for any tolerance $\varepsilon>0$.
Both methods also apply to smooth objectives. In this setting, the basic method yields a single Riemannian accelerated gradient scheme for geodesically convex and geodesically $\mu$-strongly convex problems. Its iteration and parameter recurrence extend continuously to $\mu=0$, and a single analysis establishes accelerated rates in both cases.
Numerical experiments on nonsmooth composite and smooth problems demonstrate the effectiveness of the proposed methods.
\end{abstract}

	


\section{Introduction}
		\label{sec:intro}
In this paper, we develop accelerated proximal gradient methods for composite optimization on Riemannian manifolds. We consider
\begin{equation}
	\label{eq:f+h}
\min_{x\in\mathcal{M}} F(x)=f(x)+h(x),
\end{equation}
where $\mathcal{M}$ is a finite-dimensional Riemannian manifold, $f:\mathcal{M}\to\mathbb{R}$ is continuously differentiable, and $h:\mathcal{M}\to\mathbb{R}$ is finite-valued, locally Lipschitz continuous, and possibly nonsmooth. Such problems combine manifold constraints with regularization and arise, for example, in sparse principal component analysis~\cite{genicot2015weakly,zou2006sparse}, and the computation of regularized Riemannian means. The composite structure allows one to linearize the smooth term while retaining the nonsmooth term in a subproblem, and then update the iterate through a suitable mapping to the manifold. This provides a natural framework for applying proximal gradient methods to nonsmooth optimization under manifold constraints.

In Euclidean spaces, exploiting the composite structure has led to well-established accelerated convergence guarantees. Nesterov's accelerated methods and FISTA, developed by Beck and Teboulle, improve the $\mathcal{O}(1/k)$ objective residual bound of proximal gradient methods to $\mathcal{O}(1/k^2)$ for convex composite problems~\cite{nesterov2013gradient,beck2009fast,beck2017first,nesterov1983method,tseng2008accelerated}, where $k$ is the iteration count. Under suitable strong convexity assumptions, accelerated linear convergence rates with a square-root dependence on the condition number are also available~\cite{beck2017first,nesterov2018lectures,aujol2024fista,daspremont2021acceleration,schmidt2011convergence}. For nonconvex composite problems, global convergence of accelerated proximal gradient methods has been established under suitable assumptions~\cite{nesterov2013gradient, ghadimi2016accelerated,li2015accelerated}. The method studied in this paper is related to these approaches, as it reduces to a FISTA-type method in Euclidean space. On manifolds, however, proximal gradient mappings, extrapolation, and the associated convergence analysis must explicitly account for the geometry.

Existing Riemannian proximal gradient methods can be broadly classified into three classes according to the construction of their proximal subproblems. The first two linearize the smooth term in the tangent space but differ in their treatment of the nonsmooth term: one evaluates $h(x+\eta)$ in the ambient space~\cite{chen2020proximal}, whereas the second employs the retraction pullback $h(R_x(\eta))$~\cite{huang2022riemannian}. Here, $\eta$ is a tangent vector and $R_x$ a retraction, both based at $x\in\m$. The point $x+\eta$ need not lie on the manifold, whereas the pullback formulation requires only values of $h$ on the manifold. The third class first takes a Riemannian gradient step for the smooth term and then performs an intrinsic proximal step by regularizing the nonsmooth term with the squared geodesic distance from the resulting point~\cite{feng2022proximal,bergmann2026convex,bergmann2026nonconvex}. A natural question is how to incorporate momentum acceleration into these proximal constructions and establish corresponding accelerated convergence guarantees.

We focus on the second class of proximal models, based on retraction pullbacks. Assuming the relevant geometric inequalities hold, Riemannian proximal gradient methods based on this construction achieve an $\mathcal{O}(1/k)$ convergence rate when the composite objective satisfies pullback convexity~\cite{huang2022riemannian} and converge linearly under pullback strong convexity~\cite{choi2025linear}. Momentum variants using FISTA acceleration parameters, together with corresponding safeguard strategies, have also been proposed~\cite{huang2022riemannian,zhao2024proximal}. Although these methods exhibit acceleration in numerical experiments, corresponding accelerated convergence guarantees have not been established. Developing such guarantees requires accounting for the effects of manifold geometry on the proximal step and momentum updates in both algorithm design and convergence analysis. Specifically, three main challenges remain in establishing an accelerated convergence theory for these methods.

First, establishing accelerated convergence rates for proximal gradient methods based on  retraction pullbacks requires combining geodesic convexity, pullback convexity, and geometric comparisons across different base points. Even when the nonsmooth term is geodesically convex, its pullback under the exponential map need not be convex in the tangent space since the exponential map does not generally map line segments to geodesics. The analysis must therefore combine two separate assumptions: geodesic convexity of the smooth term and pullback convexity of the nonsmooth term. Moreover, parallel transport does not eliminate curvature effects in comparisons across base points, so displacement identities from linear spaces cannot be applied directly. Deriving a suitable recurrence for the accelerated iteration therefore requires combining the one-step proximal estimate with an appropriate choice of momentum parameters and geometric comparison bounds. Complexity lower bounds for geodesically convex functions further show that geometric dependence in convergence rates is generally unavoidable~\cite{criscitiello2022negative,criscitiello2023curvature,hamilton2021nogo,zhang2016first}. The analysis must therefore make explicit how the accelerated rates depend on curvature and on the region containing the iterates.

Second, many optimization problems on manifolds are not globally geodesically convex. Even when suitable convexity conditions hold locally, verifying them and estimating the associated parameters can be difficult. Moreover, momentum iterates may leave the region on which the acceleration analysis relies. A descent mechanism that does not rely on convexity is therefore needed to safeguard the accelerated iteration. Proximal gradient safeguards have been studied for this purpose in~\cite{huang2022riemannian,huang2022extension,zhao2024proximal}. For the method considered in this work,  it remains to establish global convergence and iteration complexity of a safeguarded variant without assuming geodesic convexity of the objective.

Finally, even in the smooth setting, Riemannian Nesterov-type accelerated gradient methods are often developed separately for geodesically convex and geodesically strongly convex objectives~\cite{liu2017accelerated,zhang2018estimate,ahn2020riemannian,jin2022understanding,kim2022accelerated} or adapted from one setting to the other through reductions or restarts~\cite{martinez2022global,martinez2023accelerated}.
In Euclidean space, some accelerated methods use the same update rules and parameter construction for convex and strongly convex objectives, with the corresponding accelerated rates established through a single analysis~\cite{kim2023unifying,nesterov2018lectures,daspremont2021acceleration}.
These results motivate the following question for Riemannian accelerated proximal gradient methods applied to Problem~\eqref{eq:f+h} when the nonsmooth term vanishes ($h\equiv0$). 
Can the method for geodesically convex objectives be recovered by setting the strong convexity parameter to zero in the same update rules and parameter recurrence used for geodesically strongly convex objectives? Can a single convergence estimate, valid in both cases, then yield an accelerated sublinear rate in the geodesically convex case and an accelerated linear rate in the geodesically strongly convex case?

To address these three issues, we propose a FISTA-type Riemannian accelerated proximal gradient method and develop an adaptive restart variant with proximal gradient safeguards. The main contributions are as follows.

\begin{itemize}
\item \textbf{Accelerated convergence for convex composite optimization on manifolds.} 
We develop a FISTA-type Riemannian accelerated proximal gradient method for Problem~\eqref{eq:f+h} based on exponential map pullbacks. For geodesically $\mu$-strongly convex $f$ ($\mu\geqslant0$) and pullback $\rho$-convex $h$ (see \Cref{def:rho-convex}), we establish two convergence guarantees for the objective gap: an $\mathcal{O}(k^{-2})$ rate when $\mu\geqslant\rho$ and an accelerated linear rate when $\mu>\rho$, where $k$ denotes the outer iteration index. These guarantees hold under the stated smoothness and parameter conditions on a bounded geodesically convex region containing a minimizer and all iterates. The rates depend explicitly on the sectional curvature bounds and the diameter of the region. Both rates are proved using a curvature-corrected potential function constructed for the proposed method.

\item \textbf{Safeguarded adaptive restart for convex and nonconvex Riemannian optimization.}
We further propose an adaptive restart variant
of the Riemannian accelerated proximal gradient method, which combines proximal gradient safeguards with momentum restarts and adaptive updates of the proximal coefficient and safeguard interval. The method requires no prior knowledge of the smoothness or convexity parameters. Under the stated assumptions, including geodesic smoothness of $f$ and compactness of the initial sublevel set, we prove that every accumulation point of the safeguard reference sequence is stationary. Moreover, a safeguard direction with a norm of at most $\varepsilon$ is obtained within $\mathcal{O}(\varepsilon^{-2})$ safeguard calls and outer iterations. These guarantees require neither convexity nor the parameter conditions used to establish accelerated convergence.

\item \textbf{A unified algorithm and analysis for smooth Riemannian acceleration.} 
For Problem~\eqref{eq:f+h} with $h\equiv0$, setting $\rho=0$ reduces the proposed accelerated proximal gradient method to a Riemannian accelerated gradient method. The algorithm for geodesically convex objectives is obtained directly by setting the strong convexity parameter $\mu$ of $f$ to zero in the iteration and parameter recurrence.  Using the same potential function, a single analysis establishes accelerated sublinear convergence for $\mu\geqslant 0$ and accelerated linear convergence for $\mu>0$. For the adaptive restart variant, the proximal gradient safeguards reduce to gradient safeguards, and the corresponding stationarity and complexity guarantees remain valid.
\end{itemize}

We evaluate the proposed methods on two nonsmooth composite problems and two smooth problems. Each group includes tests on manifolds with nonpositive and nonnegative sectional curvature. The numerical results demonstrate the efficiency of the proposed methods.

The remainder of the paper is organized as follows. \Cref{sec:notation} presents the geometric preliminaries and problem setting. \Cref{sec:RAPG-convex} introduces the Riemannian accelerated proximal gradient method and establishes its accelerated convergence rates under convexity assumptions. \Cref{sec:restart} develops the adaptive restart variant with proximal gradient safeguards and provides stationarity and iteration complexity guarantees for convex and nonconvex problems. 
\Cref{sec:experiment} reports numerical experiments, and \Cref{sec:conclusions} concludes the paper. The appendices 
provide supporting technical proofs.

\section{ Preliminaries}
\label{sec:notation}

This section presents the geometric background, relevant function classes, and  Riemannian Hessian bounds for the squared distance function. For background on Riemannian geometry and optimization, see~\cite{absil2008optimization,boumal2023introduction,lee2018introduction}.

\subsection{Riemannian geometry background}

Throughout this paper, $\m$ denotes a finite-dimensional, simply connected smooth manifold equipped with a Riemannian metric. The tangent space at $x\in\m$ is denoted by $\T_x\m$, and the tangent bundle is $\T\m=\bigsqcup_{x\in\m}\T_x\m$. A vector field on $\m$ is a map $V:\m\to\T\m$ satisfying $V(x)\in\T_x\m$ for every $x\in\m$. For tangent vectors $u,v\in\T_x\m$, the Riemannian inner product is denoted by $\langle u,v\rangle_x$, and the induced norm is given by $\|v\|_x=\sqrt{\langle v,v\rangle_x}$. The subscript $x$ is omitted when clear from context.
For $x,y\in\m$, the Riemannian distance $\dd(x,y)$ is the infimum of the lengths of piecewise smooth curves joining $x$ and $y$. The diameter of a nonempty subset $\Omega\subseteq\m$ is defined by $\operatorname{diam}(\Omega)=\sup_{x,y\in\Omega}\dd(x,y)$.
For $r>0$, the open ball in $\m$ centered at $x$ is denoted by $ B(x,r)= \{y \in \mathcal{M} : \dd(y, x)<r\}$. Similarly, for $v\in\T_x\m$, the open ball in $\T_x\m$ centered at $v$ is denoted by $B(v,r)=\left\{u \in \mathrm{T}_{x} \mathcal{M} : \left\|u-v\right\|_{x}< r\right\}$. Their closures are denoted by $\overline{B}(x,r)$ and $\overline{B}(v,r)$, respectively.

A smooth curve $\gamma:I\to\m$ is a geodesic if $\nabla_{{\gamma'}(t)}{\gamma'}(t)=0$ for all $t\in I$, where $\nabla$ denotes the Levi--Civita connection. Every sufficiently short geodesic segment minimizes the length between its endpoints.
For $x\in\m$ and $v\in\T_x\m$, let $\gamma:I_{x,v}\to\m$ be the unique maximal geodesic satisfying $\gamma(0)=x$ and $\gamma'(0)=v$, where $I_{x,v}$ is its maximal interval of existence containing $0$. The exponential map $\mathrm{Exp}:\mathcal D_{\mathrm{Exp}}\to\m$ is defined by $\mathrm{Exp}(x,v)=\mathrm{Exp}_x(v)=\gamma(1)$, where $\mathcal D_{\mathrm{Exp}}=\{(x,v)\in\T\m:1\in I_{x,v}\}$. This domain is an open neighborhood of the zero section of $\T\m$ and need not coincide with $\T\m$. If $\m$ is geodesically complete, then $\mathcal D_{\mathrm{Exp}}=\T\m$.	


A subset $\Omega\subseteq\m$ is called uniquely geodesically convex if any two points $x,y\in\Omega$ are joined by a unique minimizing geodesic in $\m$, and this geodesic lies entirely in $\Omega$. 
Let $\gamma_{x,y}:[0,1]\to\Omega$ denote this geodesic, with $\gamma_{x,y}(0)=x$ and $\gamma_{x,y}(1)=y$. Consequently, the inverse exponential map is well defined on $\Omega$, with $\mathrm{Exp}_x^{-1}(y)={\gamma'}_{x,y}(0)$. It satisfies $\|\mathrm{Exp}_x^{-1}(y)\|_x=\|\mathrm{Exp}_y^{-1}(x)\|_y=\dd(x,y)$.
For a smooth curve $\gamma:[0,1]\to\m$ with $\gamma(0)=x$ and $\gamma(1)=y$, parallel transport along $\gamma$ with respect to the Levi--Civita connection is denoted by $\mathrm{PT}_{1\gets0}^{\gamma}:\T_x\m\to\T_y\m$. It is a linear isometry and therefore preserves inner products and norms. When $x,y\in\Omega$, we write $\Gamma_x^y=\mathrm{PT}_{1\gets0}^{\gamma_{x,y}}$ for parallel transport along their unique minimizing geodesic.
		

Let $\mathcal D_R\subseteq\T\m$ be an open set containing $(x,0_x)$ for every $x\in\m$, where $0_x$ denotes the zero vector in $\T_x\m$.
A smooth map $R:\mathcal D_R\to\m$, written as $R(x,v)=R_x(v)$, is called a retraction if $R_x(0_x)=x$ and $\mathrm{D}R_x(0_x)[v]=v$ for all $(x,v)\in \mathcal D_R$. For each $x\in\m$, the domain of $R_x$ is $\mathcal D_x=\{v\in\T_x\m:(x,v)\in\mathcal D_R\}$. We call $R$ globally defined if $\mathcal D_R=\T\m$. The exponential map on its natural domain is a second-order retraction~\cite{absil2008optimization,boumal2023introduction}.

\subsection{Function classes}
For a smooth function $f:\m\to\R$, we denote its Riemannian gradient at $x\in\m$ by $\operatorname{grad}f(x)\in\T_x\m$. A function $h:\m\to\R$ is locally Lipschitz if, for every $x\in\m$, there exist $r>0$ and $\ell>0$ such that $|h(z)-h(y)|\leqslant\ell\dd(z,y)$ holds for all $y,z\in B(x,r)$.

For a locally Lipschitz function $h$, let $\hat h_x=h\circ R_x:\mathcal D_x\to\R$ denote its pullback by $R_x$. The Riemannian subdifferential of $h$ at $x$, understood in the Clarke sense, is given by $\partial h(x)=\partial\hat h_x(0_x)$, where the right-hand side is the Clarke subdifferential computed in $\T_x\m$ equipped with the inner product $\langle\cdot,\cdot\rangle_x$; see~\cite{hosseini2011generalized,hosseini2018line}. Any vector $v\in\partial h(x)$ is called a Riemannian subgradient of $h$ at $x$. A point $x\in\m$ is called stationary for $h$ if $0_x\in\partial h(x)$. Every local minimizer of $h$ is stationary~\cite{grohs2016subgradient}. If $f$ is smooth and $h$ is locally Lipschitz, the sum rule gives $\partial(f+h)(x)=\operatorname{grad}f(x)+\partial h(x)$; see~\cite[Proposition~3.1]{hosseini2011generalized} and~\cite[Corollary~1 in Section~2.3]{Clarke1990optimization}.

We now introduce retraction-convex sets, together with retraction convexity and smoothness for functions.

\begin{definition}[Retraction-convex set]
Let $\Omega\subseteq\m$. For $x,y\in\Omega$, define
    \[
    	\mathcal A_R^\Omega(x,y)
    	=
    	\left\{
    	\eta\in\mathcal D_{x}:
    	R_x(\eta)=y
    	\ \text{and}\
    	R_x(t\eta)\in\Omega
    	\ \text{for all }t\in[0,1]
    	\right\}.
    \]
The set $\Omega$ is called retraction-convex with respect to $R$ if $\mathcal A_R^\Omega(x,y)\neq\varnothing$ for every ordered pair $(x,y)\in\Omega\times\Omega$.
\end{definition}

\begin{definition}[Retraction convexity and smoothness]
\label{def:retraction-convex-smooth}
Let $\Omega\subseteq\m$ be retraction-convex with respect to $R$, and let $f:\m \to\R$.

For $\mu\geqslant0$, the function $f$ is called $\mu$-strongly retraction-convex with respect to $R$, or simply $\mu$-strongly $R$-convex, on $\Omega$ if, for every $x,y\in\Omega$ and every $\eta\in\mathcal A_R^\Omega(x,y)$,
\begin{equation}
    f\bigl(R_x(t\eta)\bigr)
    \leqslant
    (1-t)f(x)+tf(y)
    -\frac{\mu}{2}t(1-t)\|\eta\|_x^2,
    \qquad t\in[0,1].
    \label{eq:R-strong-convexity}
\end{equation}
When $\mu=0$, we simply say that $f$ is retraction-convex with respect to $R$, or $R$-convex, on $\Omega$.

For $L > 0$, the function $f$ is called $L$-retraction-smooth on $\Omega$ with respect to $R$ if $f$ is differentiable and satisfies
\begin{equation}
    f(y)
    \leqslant
    f(x)
    +\left\langle\grad f(x),\eta\right\rangle_x
    +\frac{L}{2}\|\eta\|_x^2
    \label{eq:R-smoothness}
\end{equation}
for every $x,y\in\Omega$ and every $\eta\in\mathcal A_R^\Omega(x,y)$.
\end{definition}

	
	

\begin{remark}
\label{rem:R-smoothness-convention}
Condition~\eqref{eq:R-smoothness} is a one-sided quadratic upper bound; see~\cite[Definition~1]{huang2022riemannian}. The two-sided remainder bound
\begin{equation}
    \left|
    f(y)-f(x)-\left\langle\grad f(x),\eta\right\rangle_x
    \right|
    \leqslant
    \frac{L}{2}\|\eta\|_x^2
    \label{eq:R-two-sided-smoothness}
\end{equation}
for every $x,y\in\Omega$ and every $\eta\in\mathcal A_R^\Omega(x,y)$ implies~\eqref{eq:R-smoothness}; see~\cite[Assumption~2.6]{boumal2019global}. Throughout this paper, $L$-retraction-smoothness refers to the one-sided condition~\eqref{eq:R-smoothness}, which is sufficient for our analysis.
\end{remark}

\begin{remark}
\label{rem:R-first-order-Exp}
Suppose that $f$ is differentiable and $\mu$-strongly $R$-convex on $\Omega$. Then, for every $x,y\in\Omega$ and every $\eta\in\mathcal A_R^\Omega(x,y)$,
\begin{equation}
    f(y)
    \geqslant
    f(x)
    +\left\langle\grad f(x),\eta\right\rangle_x
    +\frac{\mu}{2}\|\eta\|_x^2.
    \label{eq:R-first-order-strong-convexity}
\end{equation}
%
When $R=\mathrm{Exp}$ and $\Omega$ is uniquely geodesically convex, the vector $\eta=\mathrm{Exp}_x^{-1}(y)$ associated with the unique minimizing geodesic from $x$ to $y$ belongs to $\mathcal A_{\mathrm{Exp}}^\Omega(x,y)$ and satisfies $\|\eta\|_x=\dd(x,y)$. 
Substituting this vector into~\eqref{eq:R-first-order-strong-convexity} and~\eqref{eq:R-smoothness} yields the standard inequalities for geodesic $\mu$-strong convexity and geodesic $L$-smoothness of $f$, respectively; see, for example,~\cite{zhang2016first}. When $\mu=0$, the former reduces to the first-order inequality for geodesic convexity of $f$.
\end{remark}

We next introduce a notion of convexity for functions through their pullbacks to tangent spaces. For each $x\in\Omega$, set $\mathcal D_{x,\Omega}=\{\eta\in\mathcal D_x:R_x(\eta)\in\Omega\}$.

\begin{definition}[Pullback $\rho$-convexity]
\label{def:rho-convex}
Let $\Omega\subseteq\m$, and let $h$ be locally Lipschitz on an open neighborhood of $\Omega$. For $\rho\in\R$, the function $h$ is said to be pullback $\rho$-convex on $\Omega$ with respect to $R$ if, for every $x\in\Omega$, the function $\widetilde h_x(\eta):=h(R_x(\eta))+\frac{\rho}{2}\|\eta\|_x^2$ is convex on every convex subset $U_x\subseteq\mathcal D_{x,\Omega}$.
\end{definition}

\begin{remark}
\label{rem:pullback-subgradient}
Pullback $\rho$-convexity extends the Euclidean notion of $\rho$-convexity to retraction pullbacks; see~\cite{vial1983strong}.
Under our sign convention, the pullback $h\circ R_x$ is $\rho$-weakly convex for $\rho>0$, convex for $\rho=0$, and $(-\rho)$-strongly convex for $\rho<0$ on each convex subset of $\mathcal D_{x,\Omega}$.  The cases $\rho>0$ and $\rho=0$ are referred to as \emph{weak retraction convexity} in~\cite[Definition~3]{choi2025linear} and \emph{retraction convexity} in~\cite[Definition~2]{huang2022riemannian}, respectively. To avoid confusion with $R$-convexity as defined in \Cref{def:retraction-convex-smooth}, we use the term \emph{pullback $\rho$-convexity}, or simply \emph{pullback convexity} when $\rho=0$.
\end{remark}

\begin{remark}
Suppose that $h$ is pullback $\rho$-convex on $\Omega$. Fix $x\in\Omega$, and let $U_x\subseteq\mathcal D_{x,\Omega}$ be open and convex in $\T_x\m$. For every $\omega\in U_x$, the Clarke subdifferential $\partial\widetilde h_x(\omega)$ is nonempty and satisfies $\partial\widetilde h_x(\omega)=\partial(h\circ R_x)(\omega)+\rho\omega$. Consequently, for every $\omega,\eta\in U_x$ and every $g\in\partial(h\circ R_x)(\omega)$,
    $h(R_x(\eta))
    \geqslant
    h(R_x(\omega))
    +\left\langle g,\eta-\omega\right\rangle_x
    -\frac{\rho}{2}\|\eta-\omega\|_x^2$.
\end{remark}

For twice continuously differentiable functions $h$ on $\m$, sufficient conditions for local convexity of their pullbacks are given in~\cite[Lemma~3]{huang2022riemannian}. Under the stated assumptions, \Cref{lem:weakly-convex} establishes local $\rho$-weak convexity of the pullbacks of possibly nondifferentiable functions for some $\rho>0$.

\subsection{Riemannian Hessian bounds for the squared distance function}
\label{subsec:zeta-and-delta}

Let $\Omega\subseteq\m$ be a bounded uniquely geodesically convex set with diameter $D_\Omega=\operatorname{diam}(\Omega)$. Assume that, for every $p\in\Omega$, the function $q_p=\frac12\dd^2(\cdot,p)$ is smooth on a neighborhood of $\Omega$. Suppose that the sectional curvature on $\Omega$ is bounded below by $\kappa_{\min}$ and above by $\kappa_{\max}$, with $D_\Omega<\pi/\sqrt{\kappa_{\max}}$ when $\kappa_{\max}>0$. Define
\begin{equation}
    \label{eq:eigenvalue-zeta}
    \zeta(D_\Omega)=
    \begin{cases}
        \sqrt{-\kappa_{\min}}\,D_\Omega\,
        \coth\!\left(\sqrt{-\kappa_{\min}}\,D_\Omega\right),
        & \text{if } \kappa_{\min}<0,\\[2mm]
        1,
        & \text{if } \kappa_{\min}\geqslant0,
    \end{cases}
\end{equation}
and
\begin{equation}
    \label{eq:eigenvalue-delta}
    \delta(D_\Omega)=
    \begin{cases}
        1,
        & \text{if } \kappa_{\max}\leqslant0,\\[2mm]
        \sqrt{\kappa_{\max}}\,D_\Omega\,
        \cot\!\left(\sqrt{\kappa_{\max}}\,D_\Omega\right),
        & \text{if } \kappa_{\max}>0.
    \end{cases}
\end{equation}
At $D_\Omega=0$, both functions are defined by continuity as $1$.

Hessian comparison yields
    $\delta(D_\Omega)\|v\|_x^2
    \leqslant
    \left\langle\operatorname{Hess}q_p(x)[v],v\right\rangle_x
    \leqslant
    \zeta(D_\Omega)\|v\|_x^2$
for every $p,x\in\Omega$ and $v\in\T_x\m$; see~\cite[Appendix~B, Lemma~2]{alimisis2020continuous} and~\cite[proof of Theorem~3.15]{lezcano2020curvature}. Thus, $\delta(D_\Omega)$ and $\zeta(D_\Omega)$ provide uniform lower and upper bounds on the eigenvalues of $\operatorname{Hess}q_p(x)$.
For fixed curvature bounds, $\delta(D_\Omega)\leqslant1\leqslant\zeta(D_\Omega)$. Moreover, $\zeta(D_\Omega)$ is increasing in $D_\Omega$ when $\kappa_{\min}<0$, whereas $\delta(D_\Omega)$ is decreasing on $(0,\pi/\sqrt{\kappa_{\max}})$ when $\kappa_{\max}>0$. In Euclidean space, taking $\kappa_{\min}=\kappa_{\max}=0$ gives $\zeta(D_\Omega)=\delta(D_\Omega)=1$. 

We also use the notation $\zeta_\Omega=\zeta(D_\Omega)$ and
$\delta_\Omega=\delta(D_\Omega)$.
When $\Omega$ is clear from the context, we omit the subscript and
write $D$, $\zeta$, and $\delta$.

\section{Riemannian accelerated proximal gradient method for convex optimization}
	\label{sec:RAPG-convex}
	This section develops a Riemannian accelerated proximal gradient method and establishes its accelerated convergence rates through potential function analysis. The proof relies on inequalities for the weighted objective residual and norm comparisons across tangent spaces.
	
\subsection{Algorithm development}
\label{subsec:Algo}
\subsubsection{Riemannian proximal mapping}
	\label{subsubsec:proximal map}
	
In the Riemannian proximal gradient framework, we adopt the proximal mapping introduced in~\cite{huang2022riemannian}, defined as follows:
\begin{equation}\label{eq:eta_x}
\begin{cases}
\ell_x(\eta)=\left\langle\operatorname{grad}f(x),\eta\right\rangle+\tfrac{\theta L}{2}\|\eta\|^2+h(R_x(\eta)), \quad \eta\in\mathcal{D}_x,\\
\eta_x\in\mathcal{D}_x \text{ is a stationary point of }\ell_x\text{ satisfying }\ell_x(\eta_x)\leqslant\ell_x(0_x),\\
y=R_x(\eta_x), \qquad \text{(proximal point)}
\end{cases}
\end{equation}
where $f:\m\to\R$ is differentiable, $h:\m\to\R$ is continuous, and $\theta L$ is a prescribed positive constant.

\begin{remark}
Several notions of Riemannian proximal mappings have been proposed, including those in~\cite{chen2020proximal,huang2022riemannian} and the intrinsic proximal mapping introduced in~\cite{feng2022proximal,bergmann2026nonconvex,bergmann2026convex}. We adopt the mapping of~\cite{huang2022riemannian} because its associated proximal gradient inequality plays a central role in our convergence analysis. In particular, this inequality allows us to establish an $\mathcal{O}(1/k)$ convergence rate for the Riemannian proximal gradient method for geodesically convex problems on a broader class of Riemannian manifolds. Although the intrinsic proximal mapping also satisfies a proximal gradient inequality, Bergmann et al.~\cite{bergmann2026convex} observed that this inequality alone does not yield analogous convergence rate guarantees, even on Hadamard manifolds.
\end{remark}
	
To the best of our knowledge, the well-definedness of the Riemannian proximal mapping defined by~\eqref{eq:eta_x} has not been explicitly addressed in the literature. The following proposition establishes the existence of a proximal point, thereby ensuring that this mapping is well defined.

\begin{proposition}[Existence of a proximal point]
\label{prop:existence-proximal-point}
Suppose that $h:\m\to\R$ is locally Lipschitz and bounded below. Fix $x\in\m$ and define $\mathcal L_x=\{\eta\in\mathcal D_x:\ell_x(\eta)\leqslant\ell_x(0_x)\}$. Assume that one of the following conditions holds:
\begin{enumerate}[(i)]
    \item $R_x$ is globally defined, that is, $\mathcal D_x=\T_x\m$.
    \item $R_x$ is defined on an open subset $\mathcal D_x\subsetneq\T_x\m$, and the sublevel set $\mathcal L_x$  is compact.
\end{enumerate}
Then $\ell_x$ admits a global minimizer $\eta_x\in\mathcal D_x$ satisfying $0\in\partial\ell_x(\eta_x)$ and $\ell_x(\eta_x)\leqslant\ell_x(0_x)$. Consequently, $R_x(\eta_x)$ is a proximal point in the sense of~\eqref{eq:eta_x}.
\end{proposition}

\begin{proof}
Since $h$ is locally Lipschitz and $R_x$ is smooth, $\ell_x$ is locally Lipschitz, and hence continuous, on $\mathcal D_x$. We show that the sublevel set $\mathcal L_x$ is compact under either condition.

In case~(i), choose $m_h\in\R$ such that $h(y)\geqslant m_h$ for all $y\in\m$. Then, for every $\eta\in\T_x\m$,
\begin{align*}
    \ell_x(\eta)
    &=\langle\grad f(x),\eta\rangle +\frac{\theta L}{2}\|\eta\|_x^2+h(R_x(\eta))\notag \\
    &\geqslant-\|\grad f(x)\| \|\eta\| + \frac{\theta L}{2}\|\eta\|^2+m_h.
\end{align*}
Since $\theta L>0$, this estimate implies that $\ell_x(\eta)\to+\infty$ as $\|\eta\| \to+\infty$. Thus $\ell_x$ is coercive on $\T_x\m$. Together with continuity, this shows that $\mathcal L_x$ is closed and bounded, and hence compact. In case~(ii), compactness holds by assumption.

Since $0_x\in\mathcal L_x$, this set is nonempty. By continuity and compactness, $\ell_x$ attains its minimum on $\mathcal L_x$ at some $\eta_x\in\mathcal L_x$. In particular,
    $\ell_x(\eta_x)\leqslant\ell_x(0_x)$.
For every $\eta\in\mathcal D_x\setminus\mathcal L_x$, the definition of $\mathcal L_x$ gives
   $ \ell_x(\eta)>\ell_x(0_x)\geqslant\ell_x(\eta_x)$.
Hence $\eta_x$ is a global minimizer of $\ell_x$ over $\mathcal D_x$.

Finally, since $\mathcal D_x$ is open and $\ell_x$ is locally Lipschitz, the Clarke Fermat rule yields $0\in\partial\ell_x(\eta_x)$. Together with the preceding value inequality, this proves that $R_x(\eta_x)$ is a proximal point in the sense of~\eqref{eq:eta_x}.
\end{proof}
	
\subsubsection{Algorithm description}
\label{subsubsec:alg}
We first state the assumptions on the objective function over a region $\Omega\subseteq\mathcal{M}$ that underlie the proposed algorithm and its convergence analysis. Although the preceding notions are formulated for a general retraction, we use the exponential map in the algorithm and analysis below unless otherwise stated.

\begin{assumption}
\label{assum:omega}
Let $\Omega\subseteq\m$ be an open, uniquely geodesically convex set with compact closure $\overline{\Omega}$. Assume that $\operatorname{diam}(\Omega)\leqslant D<\infty$ and that the sectional curvature of $\m$ on $\Omega$ is bounded below by $\kappa_{\min}$ and above by $\kappa_{\max}$.
\end{assumption}

Geometric assumptions of this type are commonly imposed in the analysis of accelerated first-order methods for geodesically convex optimization on Riemannian manifolds; see, e.g.,~\cite{ahn2020riemannian,alimisis2021momentum,kim2022accelerated}.

\begin{assumption}\label{assum:f-l-smooth,h-lipschitz}
		The function $ f $ is geodesically $ L $-smooth on $\Omega$, and $h$ is locally Lipschitz on~$\Omega$.  
	\end{assumption}
	
\begin{assumption}\label{assum:f,h-convex}
		The function $f$ is geodesically $\mu$-strongly convex on $\Omega$
		for some $\mu\geqslant0$, and $h$ is pullback $\rho$-convex on
		$\Omega$ with respect to the exponential map.
	\end{assumption}
	
Under Assumptions~\ref{assum:omega}--\ref{assum:f,h-convex}, we present a Riemannian accelerated proximal gradient method for Problem~\eqref{eq:f+h}; see~\Cref{algrm:rapg-sc}. The same scheme covers both the strongly convex regime $\mu>\rho$ and the convex regime $\mu=\rho$, with the parameters for the strongly convex case reducing continuously to those for the convex case as $\mu-\rho\downarrow0$ (equivalently, $q\downarrow0$).

\begin{algorithm}[!htb]
\caption{A Riemannian accelerated proximal gradient method}
\label{algrm:rapg-sc}
\begin{algorithmic}[1]
\Require Initial point $x_0\in\m$; parameters $(L,\mu,\rho,\zeta,\xi,\theta,A_0)$ satisfying $L>\mu\geqslant\rho$, $\mu\geqslant0$, $\xi\geqslant\zeta\geqslant1$, $\theta\geqslant1$, $\theta>\frac{\rho+(\mu-\rho)\xi}{L}$, and $A_0\geqslant\frac{\xi(\xi-1)}{1-q\xi}$, where $q=\frac{\mu-\rho}{\theta L-\rho}$.
\State Initialize $s_0=0_{x_0}\in\mathrm T_{x_0}\m$.
\For{$k=0,1,2,\ldots$}
\State Compute $A_{k+1}$, $\beta_k$, $\gamma_k$, and $\tau_k$ as follows:\label{alg:st3}
\begin{align}
A_{k+1} &= \frac{\xi+2\xi A_k+\sqrt{\xi^2+4\xi^2A_k+4q\xi A_k^2}}{2(\xi-q)}, \label{eq:A}\\
\beta_k &= \frac{\xi+qA_k}{\xi+qA_{k+1}}, \label{eq:beta_k}\\
\gamma_k &= \frac{A_{k+1}-A_k}{\xi+qA_{k+1}}, \label{eq:gamma_k}\\
\tau_k &= \frac{\beta_kA_{k+1}}{\gamma_kA_k+\beta_kA_{k+1}}. \label{eq:tau_k}
\end{align}
\State Set $y_k=\e_{x_k}(\tau_ks_k)$.\label{alg:st6}
\State Compute a stationary point $\eta_{y_k}\in\mathrm T_{y_k}\m$ of the model $\ell_{y_k}(\eta)=\left\langle\operatorname{grad}f(y_k),\eta\right\rangle+\frac{\theta L}{2}\|\eta\|^2+h\bigl(\e_{y_k}(\eta)\bigr)$ satisfying $\ell_{y_k}(\eta_{y_k})\leqslant\ell_{y_k}(0_{y_k})$.\label{alg:st5}
\State Set $x_{k+1}=\e_{y_k}(\eta_{y_k})$.\label{alg:st9}
\State Set $v_{y_k}=\beta_k(1-\tau_k)\Gamma_{x_k}^{y_k}(s_k)+\gamma_k\eta_{y_k}$.\label{alg:st10}
\State Set $s_{k+1}=\Gamma_{y_k}^{x_{k+1}}\bigl(v_{y_k}-\eta_{y_k}\bigr)$.\label{alg:s_k+1}
\EndFor
\end{algorithmic}
\end{algorithm}

\begin{figure}[!htb]
\centering
\resizebox{0.65\linewidth}{!}{%
\begin{tikzpicture}[x=0.008cm,y=-0.008cm,line cap=round,line join=round,
  font=\fontsize{20}{24}\selectfont,
  geodesic/.style={draw=originalcurve,line width=2.05pt},
  auxiliarycurve/.style={geodesic,dash pattern=on 3.0pt off 2.4pt},
  dashedvector/.style={draw=originalcurve,line width=1.35pt,dash pattern=on 1.9pt off 1.7pt,-{Latex[length=3.6mm,width=2.3mm]}},
  redvector/.style={draw=originalred,line width=1.55pt,-{Latex[length=3.6mm,width=2.3mm]}},
  every node/.style={inner sep=0pt}]
\path[use as bounding box] (0,0) rectangle (2307,1227);
\coordinate (X) at (141.51,771.02);
\coordinate (Y) at (655.42,378.26);
\coordinate (Z) at (977.45,235.76);
\coordinate (Xnext) at (993.50,834.31);
\coordinate (Ynext) at (1443.57,987.13);
\coordinate (Znext) at (1571.65,1164.23);

\draw[geodesic] (X) .. controls (312.81333,606.21556) and (484.11667,475.29556) .. (Y);
\draw[auxiliarycurve] (Y) .. controls (762.76333,317.45488) and (870.10667,269.95488) .. (Z);
\draw[geodesic] (Y) .. controls (820,498) and (940,620) .. (Xnext);
\draw[geodesic] (Xnext) .. controls (1143.52333,834.31) and (1293.54667,830.06361) .. (Ynext);
\draw[auxiliarycurve] (Ynext) .. controls (1486.26333,1031.82763) and (1528.95667,1089.58912) .. (Znext);

\draw[dashedvector] (Y)--(1003,142);
\draw[dashedvector] (Y)--(1196,772);
\draw[redvector] (Y)--(1736,693);
\draw[dashedvector] (Xnext)--(1712,839);

\foreach \p in {X,Y,Z,Xnext,Ynext,Znext}
  \filldraw[fill=originalblue,draw=black,line width=0.4pt] (\p) circle[radius=1.9pt];

\node[anchor=east] at (122,771) {$x_k$};
\node[anchor=south east] at (637,355) {$y_k$};
\node[anchor=west] at (1003,245) {$z_k=\e_{x_k}(s_k)$};
\node[anchor=east] at (969,835) {$x_{k+1}$};
\node[anchor=east] at (1409,1004) {$y_{k+1}$};
\node[anchor=east] at (1539,1165) {$z_{k+1}$};

\node[anchor=south west] at (998,143) {$\operatorname{Exp}_{y_k}^{-1}(z_k)$};
\node[anchor=base west] at (1174,678) {$\eta_{y_k}$};
\node[anchor=south] at (1571,609) {$v_{y_k}$};
\node[anchor=north west] at (1640,851) {$s_{k+1}=\Gamma_{y_k}^{x_{k+1}}\!\left(v_{y_k}-\eta_{y_k}\right)$};
\end{tikzpicture}
}
\caption{Geometric illustration of the Riemannian accelerated proximal gradient updates. The dashed curve segments represent auxiliary geodesic extensions to $z_k$ and $z_{k+1}$ when defined. }
\label{fig:rapg2}
\end{figure}

The constants $L$, $\mu$, and $\rho$ characterize the smoothness and convexity properties of $f$ and $h$ in~\eqref{eq:f+h}, while $\zeta$ depends on the geometry of the manifold and the region $\Omega$. The algorithmic parameters $\xi$, $\theta$, and $A_0$ are chosen to satisfy the admissibility conditions derived from the convergence analysis. In particular, the condition $\theta>\frac{\rho+(\mu-\rho)\xi}{L}$ ensures that $1-q\xi>0$, so the lower bound $A_0\geqslant\frac{\xi(\xi-1)}{1-q\xi}$ is well defined. The update formulas for $\beta_k$, $\gamma_k$, and $\tau_k$ likewise follow from the Lyapunov analysis presented in the subsequent sections.

Steps~\ref{alg:st6}--\ref{alg:s_k+1} of \Cref{algrm:rapg-sc} generate the iterate sequences $\{x_k\}$ and $\{y_k\}$ and the auxiliary sequence $\{s_k\}$, where $s_k\in\mathrm T_{x_k}\m$. Whenever $\e_{x_k}(s_k)$ is defined, we associate with $s_k$ the auxiliary point $z_k=\e_{x_k}(s_k)$. For such $k$, the relation $y_k=\e_{x_k}(\tau_k s_k)$ implies that $x_k$, $y_k$, and $z_k$ lie on the geodesic $t\mapsto\e_{x_k}(t s_k)$. With the inverse image of $z_k$ at $y_k$ selected along this geodesic and parallel transport taken along the same geodesic, we obtain $\e_{y_k}^{-1}(z_k)=(1-\tau_k)\Gamma_{x_k}^{y_k}(s_k)$. Consequently, $v_{y_k}$ admits the geometric representation $v_{y_k}=\beta_k\e_{y_k}^{-1}(z_k)+\gamma_k\eta_{y_k}$; see Figure~\ref{fig:rapg2}.

For notational convenience, we henceforth use $z_k=\e_{x_k}(s_k)$ and $s_k=\e_{x_k}^{-1}(z_k)$ as formal notation, without assuming that $\e_{x_k}(s_k)$ is defined.

\begin{remark}
\label{remark:theta}
When $\kappa_{\min}\geqslant0$, the curvature constant defined in Section~\ref{subsec:zeta-and-delta} satisfies $\zeta=1$. Choosing $\xi=1$ reduces the lower bound on $A_0$ to $A_0\geqslant0$. Moreover, $\frac{\rho+(\mu-\rho)\xi}{L}=\frac{\mu}{L}<1$, so the choice $(\xi,\theta,A_0)=(1,1,0)$ is admissible.
When $\kappa_{\min}<0$, the parameters satisfy $\xi\geqslant\zeta\geqslant 1$. Since $\theta L-\rho>0$, the condition $\theta>\frac{\rho+(\mu-\rho)\xi}{L}$ is equivalent to $q\xi<1$. This ensures that the lower bound $\frac{\xi(\xi-1)}{1-q\xi}$ for $A_0$ is well defined and nonnegative. In particular, $\theta=1$ remains admissible whenever $\rho+(\mu-\rho)\xi<L$.
\end{remark}

\begin{remark}[Relation to smooth Riemannian acceleration]
The iterative structure of \Cref{algrm:rapg-sc} is inspired by the smooth Riemannian accelerated gradient methods in~\cite{kim2022accelerated}, but its acceleration parameters follow a different recurrence. Specifically, a single recurrence for $A_k$ determines the coefficients in both the strongly convex regime $\mu>\rho$ and the convex regime $\mu=\rho$. The parameter formulas are well defined and continuous at $q=0$, so the convex case follows by direct substitution. In contrast, setting the strong convexity parameter to zero in the strongly convex method of~\cite{kim2022accelerated} does not recover its convex counterpart. In the smooth setting $h\equiv0$ and $\rho=0$, \Cref{algrm:rapg-sc} reduces to a Riemannian accelerated gradient method whose parameter recurrence likewise covers both geodesically convex and geodesically strongly convex objectives.
\end{remark}

\begin{remark}[Euclidean specialization]
\label{remark:Euclidean-specialization}
Suppose that $\m=\mathbb{R}^n$ and that the smoothness and convexity assumptions on $f$ and $h$ hold globally. With $z_k=x_k+s_k$, Steps~\ref{alg:st6}--\ref{alg:s_k+1} of \Cref{algrm:rapg-sc} reduce to
\begin{equation*}
\begin{cases}
y_k=x_k+\tau_k(z_k-x_k),\\
x_{k+1}=\operatorname{prox}_{h/(\theta L)}\left(y_k-\dfrac{1}{\theta L}\nabla f(y_k)\right),\\
z_{k+1}=\beta_k z_k+(1-\beta_k)y_k+\gamma_k(x_{k+1}-y_k),
\end{cases}
\end{equation*}
where
\(
\operatorname{prox}_{h/(\theta L)}(u)=\operatorname*{arg\,min}_{x\in\mathbb{R}^n}\left\{h(x)+\frac{\theta L}{2}\|x-u\|^2\right\}
\)
is single-valued because $\theta L>\rho$. This formulation also accommodates $\rho$-weakly convex $h$ for $\rho>0$ and $(-\rho)$-strongly convex $h$ for $\rho<0$, subject to the algorithm's parameter conditions. For $\mu\geqslant\rho$, the objective $f+h$ is convex, and it is $(\mu-\rho)$-strongly convex when $\mu>\rho$.

For convex $h$, choosing $\rho=0$, $\xi=\theta=1$, and $A_0=0$ recovers the method in~\cite[Algorithm~19]{daspremont2021acceleration} initialized with $L_0=L$. Global $L$-smoothness ensures that its backtracking condition is satisfied with $L_{k+1}=L$ at every iteration, so $q_{k+1}=\mu/L=q$. The recurrence for $A_k$ and the coefficient $\tau_k$ coincide with those in the cited algorithm, whose coefficient $\delta_k$ satisfies
\[
\gamma_k=\frac{A_{k+1}-A_k}{1+qA_{k+1}}=\delta_k,
\qquad
\beta_k=\frac{1+qA_k}{1+qA_{k+1}}=1-q\delta_k.
\]
Since $s_0=0$ gives $z_0=x_0$, the two methods have identical initializations and iterate updates.
\end{remark}

\begin{remark}[Why a direct Riemannian extension of FISTA does not telescope]
Consider the Riemannian extension of FISTA in~\cite[Algorithm~3]{huang2019riemannianpreprint}. Starting from $y_0=x_0$ and $t_0=1$, with $\alpha_k=(t_k-1)/t_{k+1}$ and $t_{k+1}=\frac{1+\sqrt{1+4t_k^2}}{2}$, its updates are
\[
x_{k+1}=R_{y_k}(\eta_{y_k}),\qquad 
y_{k+1}=R_{y_k}\bigl((1+\alpha_k)\eta_{y_k}-\alpha_kR_{y_k}^{-1}(x_k)\bigr).
\]
Here, $\eta_{y_k}$ is a local minimizer of $\ell_{y_k}$ with proximal coefficient $L$, satisfying $\ell_{y_k}(\eta_{y_k})\leqslant\ell_{y_k}(0_{y_k})$. Assume that all required inverse retractions are well defined. Under pullback convexity of $f$ and $h$ and the smoothness bound required for the proximal gradient inequality, we obtain, for $k\geqslant0$,
\begin{equation}
\label{eq:100}
\frac{2}{L}\left[t_k^2\bigl(F(x_{k+1})-F(x_*)\bigr)-t_{k-1}^2\bigl(F(x_k)-F(x_*)\bigr)\right]\leqslant\|\widehat{W}_k\|_{y_k}^2-\|\widetilde{W}_{k+1}\|_{y_k}^2,
\end{equation}
where $x_*$ is a minimizer, $t_{-1}=0$, 
$\widehat{W}_k=(t_k-1)R_{y_k}^{-1}(x_k)+R_{y_k}^{-1}(x_*)$, and $\widetilde{W}_{k+1}=\widehat{W}_k-t_k\eta_{y_k}$.

In Euclidean space, the scheme reduces to FISTA, and extrapolation gives
\[
\widehat{W}_{k+1}=(t_k-1)x_k+x_*-t_kx_{k+1}=\widetilde{W}_{k+1}.
\]
Thus, the right-hand side of~\eqref{eq:100} telescopes, yielding the $\mathcal{O}(1/k^2)$ rate; see~\cite[Lemma~4.1]{beck2009fast}.

On a manifold, the extrapolation does not generally ensure $\|\widehat{W}_{k+1}\|_{y_{k+1}}=\|\widetilde{W}_{k+1}\|_{y_k}$. Summation therefore leaves intermediate differences $\|\widehat{W}_k\|_{y_k}^2-\|\widetilde{W}_k\|_{y_{k-1}}^2$ that are not controlled by~\eqref{eq:100}, so this inequality alone does not establish an accelerated rate. In \Cref{algrm:rapg-sc}, the auxiliary tangent-vector recurrence and the metric distortion estimates for the exponential map provide the comparisons between successive potential terms needed for our convergence analysis.
\end{remark}

\subsubsection{Properties of the parameter sequences}

\begin{lemma}\label{lem:parameters}
	The parameters $\beta_{k}$, $\gamma_{k}$, $\tau_{k}$, and $A_{k}$ in \Cref{algrm:rapg-sc} satisfy the following properties: 
	\begin{enumerate}[(i)]
		   \item \label{Ak}   
			$A_{k+1} \geqslant A_{k} +1 > A_{k} \geqslant A_{0}\geqslant 0$, $\forall\,k\geqslant 0$, and 
			$A_0 = 0 \iff \left( A_0 = \frac{\xi(\xi - 1)}{1 - q \xi}  \text{ and } \xi = 1  \right)$.
			\item \label{A1} 
			$A_{1} \geqslant \frac{\xi^{2}}{1 - q\xi}$, and $A_{1}=\frac{\xi^{2}}{1 - q\xi} \iff A_0 = \frac{\xi(\xi - 1)}{1 - q \xi}$.
			\item \label{beta} 
			$\beta_{k} \in (0,1]$, $\forall\, k\geqslant 0$, and $\beta_{k}=1\iff \mu=\rho $. 
			\item \label{gamma} 
			$\gamma_{k} >1$,  $\forall\, k\geqslant 1$, and $\gamma_{0}\geqslant 1$, with $\gamma_{0}= 1\iff  A_{0}=\frac{\xi(\xi - 1)}{1 - q \xi}$.
			\item \label{tau}
			$\tau_{k} \in (0,1)$, $\forall\, k\geqslant 1$, and $\tau_{0}\in(0,1]$, with $\tau_{0}=1 \iff A_0 = 0$. 
			\item \label{eq:A_k+1_square}
			$A_{k+1}$ is the larger root of the following equation:
			\begin{equation}\label{eq:A_k+1}
				\frac{A_{k+1}}{\xi}=\frac{\left(A_{k+1}-A_{k}\right)^{2}}{\xi + q A_{k+1} }, \quad \forall\, k\geqslant 0.
			\end{equation}
			Consequently,
			\begin{equation}
				A_{k}  =  A_{k+1} - \sqrt{A_{k+1}\left (1+ \frac{q A_{k+1}}{\xi}\right )}, \quad \forall\, k\geqslant 0. \label{eq:Ak}
			\end{equation}
			Moreover,
			\begin{equation}
				1-\beta_{k}=q\gamma_{k},  \quad \forall\, k\geqslant 0; \qquad  \frac{\tau_{k}}{1-\tau_{k}}=\frac{\beta_{k}A_{k+1}}{\gamma_{k}A_{k}}, \quad \forall\, k\geqslant 1. \label{eq:beta-gamma-tau}
			\end{equation}
	\end{enumerate}
\end{lemma}

\begin{proof}
The parameter conditions imply $\xi\geqslant1$, $0\leqslant q<1/\xi$, and $A_0\geqslant a_*:=\xi(\xi-1)/(1-q\xi)\geqslant0$. Define
$\psi(t)=\frac{\xi\left(1+2t+\sqrt{1+4t+4qt^2/\xi}\right)}{2(\xi-q)}$ 
for  $t\geqslant0$.
Then $A_{k+1}=\psi(A_k)$ by~\eqref{eq:A}. The map $\psi$ is strictly increasing and satisfies $\psi(t)-t\geqslant\xi/(\xi-q)\geqslant1$, proving~\eqref{Ak}. Direct substitution gives $\psi(a_*)=\xi^2/(1-q\xi)$, so strict monotonicity also proves~\eqref{A1}, including its equality condition.

The recurrence~\eqref{eq:A} is the larger-root formula for~\eqref{eq:A_k+1}. Since $A_{k+1}>A_k$, taking the positive square root yields~\eqref{eq:Ak}. Furthermore, $\gamma_k>0$ and
\[
\gamma_k^2\stackrel{\eqref{eq:gamma_k}}{=} \frac{(A_{k+1}-A_k)^2}{(\xi+qA_{k+1})^2}\stackrel{\eqref{eq:A_k+1}}{=} \frac{A_{k+1}}{\xi(\xi+qA_{k+1})}.
\]
Thus, $\gamma_k\geqslant1$ if and only if $A_{k+1}\geqslant\xi^2/(1-q\xi)$. Items~\eqref{Ak} and~\eqref{A1} therefore give~\eqref{gamma} and its equality condition.

Finally, the definitions give
\[
\beta_k\stackrel{\eqref{eq:beta_k}}{=}\frac{\xi+qA_k}{\xi+qA_{k+1}},\qquad \tau_k\stackrel{\eqref{eq:tau_k}}{=}\frac{\beta_k A_{k+1}}{\gamma_k A_k+\beta_k A_{k+1}}.
\]
The first formula, together with $A_{k+1}>A_k$ and $q=0\iff\mu=\rho$, yields~\eqref{beta}. The second gives $0<\tau_k\leqslant1$, with $\tau_k=1$ if and only if $A_k=0$. Since $\{A_k\}$ is strictly increasing, this occurs precisely when $k=0$ and $A_0=0$. In that case, $A_0\geqslant a_*$ forces $\xi=1$, completing the proof of~\eqref{tau}.
Both identities in~\eqref{eq:beta-gamma-tau} follow directly from these formulas and together yield~\eqref{eq:A_k+1_square}.
\end{proof}

\subsection{Inequalities for the weighted objective residual}
\label{subsec:ineq-function}

For the subsequent analysis, define
\begin{align}
\widetilde{h}_x(\eta)
&=h\bigl(\e_x(\eta)\bigr)
+\frac{\rho}{2}\|\eta\|_x^2,
\label{eq:h-tilde}\\
P_k
&=\xi(\theta L-\rho)+(\mu-\rho)A_k.
\label{eq:Pk}
\end{align}

The following lemma extends the fundamental proximal gradient inequality (see~\cite[Lemma~2.3]{beck2009fast} and~\cite[Theorem~10.16]{beck2017first}) to $\mu$-strongly $R$-convex $f$ and pullback $\rho$-convex $h$. When $\mu=\rho=0$, it recovers the Riemannian inequality in~\cite[Lemma~4]{huang2022riemannian}.
It is stated for a general retraction $R$ and will be
applied with $R=\e$.
	
	\begin{lemma}[Fundamental proximal gradient inequality]
			\label{lem:proximal inequality}
			Let $\Omega\subseteq\m$ be retraction-convex in the existence
			sense with respect to $R$. Suppose that $f$ is
			$L$-retraction-smooth and $\mu$-strongly $R$-convex on $\Omega$
			(or $R$-convex when $\mu=0$), and that $h$ is pullback
			$\rho$-convex on $\Omega$ with respect to $R$. 
			Fix $\theta\geqslant1$ and $y\in\Omega$, and let $\eta_y\in\mathcal D_y$ be a stationary point of $\ell_y$
			satisfying $\ell_y(\eta_y)\leqslant\ell_y(0)$.
			Set $z=R_y(\eta_y)$ and suppose that $\eta_y\in\mathcal A_R^\Omega(y,z)$.
			Then, for any $w_y\in\mathcal D_y$ with
			$x=R_y(w_y)$ and $w_y\in\mathcal A_R^\Omega(y,x)$, it holds that
			\begin{equation*}
				F(x)-F(z)
				\geqslant
				\frac{\theta L-\rho}{2}
				\|w_y-\eta_y\|^2
				-
				\frac{\theta L-\mu}{2}
				\|w_y\|^2.
			\end{equation*}
	\end{lemma}

	\begin{proof}
		By the definition of $\widetilde h_y$, we have
		\begin{equation*}
			\begin{aligned}
				h(x)-h(z)   = h\left(R_{y}\left(w_y\right) \right)
				- h\left(R_{y}\left(\eta_{y}\right)\right)   
				 = \widetilde{h}_{y}(w_{y}) -\widetilde{h}_{y}(\eta_{y})  - \frac{\rho }{2}\left (\left \|w_{y}\right \|^{2} - \|\eta_{y}\|^{2}\right ). 
			\end{aligned}
		\end{equation*}
		Moreover, since $ \eta_{y} $ is a stationary point of $ \ell_{y}(\eta) $, we obtain
			$0\in \operatorname{grad} f\left(y\right)+ (\theta L-\rho ) \eta_{y}+ \partial\widetilde{h}_{y}(\eta_{y})$.
By the convexity of $\widetilde h_y$, the subgradient inequality gives
		\begin{equation}\label{eq:18}
			h(x)-h(z)
			\geqslant \left\langle \operatorname{grad} f\left(y\right)+ (\theta L-\rho )\eta_{y}, \eta_{y} 
			-w_{y} \right\rangle  
			- \frac{\rho }{2}\left (\left \|w_{y}\right \|^{2} - \left \|\eta_{y}\right \|^{2}\right ).
		\end{equation}

		On the other hand, the $\mu$-strong $R$-convexity and $L$-retraction-smoothness of $f$, together with $\theta\geqslant1$, give
		\[
		\begin{aligned}
		f(x)-f(z) 
		\geqslant{}& 
		\left\langle \grad f(y),
		w_y-\eta_y\right\rangle 
		+\frac{\mu}{2}\|w_y\|_y^2
		-\frac{\theta L}{2}\|\eta_y\|_y^2.
		\end{aligned}
		\]
		Adding these inequalities yields
		\begin{equation*}
			\begin{aligned}
				& \qquad F(x) - F(z)  = f(x)-f(z) + h(x)-h(z) \\
				& \ \ \geqslant  \; \frac{\mu - \rho}{2} \left\| w_{y}\right \|^{2} 
				- \frac{\theta L -\rho}{2}\left\| \eta_{y}\right \|^{2} 
				+ (\theta L-\rho ) \left \langle \eta_{y}, \eta_{y} 
				-w_{y} \right\rangle  \\
				& \ \ =  \; \frac{\theta L -\rho}{2} \left\| w_{y} - \eta_{y}\right \|^{2} 
				- \frac{\theta L - \mu}{2} \left\| w_{y}\right \|^{2},
			\end{aligned}
		\end{equation*}
		which completes the proof.
		\end{proof}

We next derive inequalities relating the weighted objective residuals at $x_k$ and $x_{k+1}$ under the following assumption.

\begin{assumption}
\label{assum:iterates}
Problem~\eqref{eq:f+h} admits a minimizer $x_*\in\Omega$, and the iterates generated by \Cref{algrm:rapg-sc} satisfy $x_k,y_k\in\Omega$ for all $k\geqslant0$.
\end{assumption}

\begin{remark}
\label{remark:boundedness}
For accelerated gradient methods in Euclidean convex optimization, arguments based on estimate sequences or Lyapunov functions yield bounded iterates under the usual global regularity assumptions, provided that a minimizer exists; see, e.g.,~\cite[Theorem~2.2.4]{nesterov2018lectures} and~\cite{bot2025iterates,jang2026point,salzo2025iterates}. One may then choose a sufficiently large bounded convex region containing a minimizer and both iterate sequences.
On general Riemannian manifolds, the convexity, smoothness, and geometric conditions required for acceleration may hold only on a restricted region. Accelerated convergence analyses therefore often work on a bounded, uniquely geodesically convex set and require the iterates to remain within it; see~\cite{zhang2018estimate,ahn2020riemannian,alimisis2021momentum,kim2022accelerated}. Specialized algorithms can enforce such containment for smooth optimization in hyperbolic and spherical settings~\cite{martinez2022global} and in certain classes of Hadamard manifolds~\cite{martinez2023accelerated}. For the proposed method on general manifolds, containment in $\Omega$ remains an explicit assumption.
\end{remark}

\begin{lemma}\label{lem:ineq-function}
Suppose that Assumptions~\ref{assum:omega}--\ref{assum:iterates} hold. Then the quantities generated by \Cref{algrm:rapg-sc} satisfy, for every $k\geqslant0$,
\[
\begin{aligned}
&A_k\bigl(F(x_k)-F(x_*)\bigr)-A_{k+1}\bigl(F(x_{k+1})-F(x_*)\bigr)\\
&\quad\geqslant\frac{P_{k+1}}{2}\left\|v_{y_k}-\e_{y_k}^{-1}(x_*)\right\|_{y_k}^2-\frac{P_k}{2}\left\|\e_{y_k}^{-1}(z_k)-\e_{y_k}^{-1}(x_*)\right\|_{y_k}^2\\
&\quad\quad+\sigma_k\|v_{y_k}\|_{y_k}^2+t_k\left\|\e_{y_k}^{-1}(z_k)\right\|_{y_k}^2,
\end{aligned}
\]
where $P_k$ is defined in~\eqref{eq:Pk} and
\[
\sigma_k=\frac{P_{k+1}}{2}(\xi-1),\qquad
t_k=
\begin{cases}
0, & k=0,\\[1mm]
\displaystyle\frac{\xi P_k^2}{2P_{k+1}}\left(\frac{qA_{k+1}}{A_k}-1\right)+\frac{P_k}{2}, & k\geqslant1.
\end{cases}
\]
Here, $\e_{y_k}^{-1}(z_k):=(1-\tau_k)\Gamma_{x_k}^{y_k}(s_k)$ is formal notation; the existence of $z_k=\e_{x_k}(s_k)$ is not assumed.
\end{lemma}

\begin{proof}
Fix $k\geqslant0$ and set
\[
    u_k=\e_{y_k}^{-1}(x_k),\qquad
    w_k=\e_{y_k}^{-1}(z_k),\qquad
    r_k=\e_{y_k}^{-1}(x_*).
\]
The expression for $w_k$ uses the formal convention for $z_k$.
%
Applying \Cref{lem:proximal inequality} at $y=y_k$, first with
$x=x_k$ and then with $x=x_*$, gives
\[
    \begin{aligned}
        F(x_k)-F(x_{k+1})
        & \geqslant \frac{\theta L-\rho}{2}\|u_k-\eta_{y_k}\|^2
          -\frac{\theta L-\mu}{2}\|u_k\|^2, \\
        F(x_*)-F(x_{k+1})
        & \geqslant \frac{\theta L-\rho}{2}\|r_k-\eta_{y_k}\|^2
          -\frac{\theta L-\mu}{2}\|r_k\|^2.
    \end{aligned}
\]
Multiplying these inequalities by the nonnegative weights $A_k$ and
$A_{k+1}-A_k$, respectively, and adding yields
\begin{equation}\label{eq:3-28}
    A_k\bigl(F(x_k)-F(x_*)\bigr)
    -A_{k+1}\bigl(F(x_{k+1})-F(x_*)\bigr)
    \geqslant N_k,
\end{equation}
where
\[
    \begin{aligned}
        N_k ={}&
        \frac{(\mu-\rho)A_k}{2}\|u_k\|^2
        +\frac{(\mu-\rho)(A_{k+1}-A_k)}{2}\|r_k\|^2 
         +\frac{(\theta L-\rho)A_{k+1}}{2}\|\eta_{y_k}\|^2 \\
        & -(\theta L-\rho)A_k\langle\eta_{y_k},u_k\rangle
          -(\theta L-\rho)(A_{k+1}-A_k)
           \langle\eta_{y_k},r_k\rangle.
    \end{aligned}
\]

To express $N_k$ in terms of $v_{y_k}$, $w_k$, and $r_k$, we use
Step~\ref{alg:st10} of \Cref{algrm:rapg-sc}, which gives
\begin{equation}\label{eq:eta}
    \eta_{y_k}=\frac{v_{y_k}-\beta_k w_k}{\gamma_k}, \quad k\geqslant 0.
\end{equation}
Step~\ref{alg:st6} of \Cref{algrm:rapg-sc}
and~\eqref{eq:beta-gamma-tau} also give
\begin{equation}\label{eq:e_(xyz)}
    u_k=-\frac{\tau_k}{1-\tau_k}w_k
        =-\frac{\beta_k A_{k+1}}{\gamma_k A_k}w_k, \quad k\geqslant 1.
\end{equation}
The parameter definitions and the recurrence for $\{A_k\}$ yield, for every $k\geqslant 0$,
\begin{equation}\label{eq:combine}
    \begin{alignedat}{2}
        & \beta_k
          \stackrel{\eqref{eq:beta_k}\eqref{eq:Pk}}{=}
          \frac{P_k}{P_{k+1}},
        & \qquad &
          \frac{(\theta L-\rho)(A_{k+1}-A_k)}{\gamma_k}
          \stackrel{\eqref{eq:gamma_k}\eqref{eq:Pk}}{=}
          P_{k+1},
        \\[0.5ex]
        & P_{k+1}-P_k
          \stackrel{\eqref{eq:Pk}}{=}
          (\mu-\rho)(A_{k+1}-A_k),
        & \qquad &
          \frac{(\theta L-\rho)A_{k+1}}{\gamma_k^2}
          \stackrel{\eqref{eq:gamma_k}\eqref{eq:A_k+1}\eqref{eq:Pk}}{=}
          \xi P_{k+1}.
    \end{alignedat}
\end{equation}

Consider first $k\geqslant1$. Substitute~\eqref{eq:eta}
and~\eqref{eq:e_(xyz)} into $N_k$. Using
$q=\frac{\mu-\rho}{\theta L-\rho}$ and the identity
\[
    \|v_{y_k}-\beta_k w_k\|^2
    +2\beta_k\langle v_{y_k}-\beta_k w_k,w_k\rangle
    =\|v_{y_k}\|^2-\beta_k^2\|w_k\|^2,
\]
we obtain
\[
    \begin{aligned}
        N_k \stackrel{\hphantom{\eqref{eq:combine}}}{=}{}&
        \frac{(\theta L-\rho)A_{k+1}}{2\gamma_k^2}
        \left[\|v_{y_k}\|^2
        +\beta_k^2\left(\frac{qA_{k+1}}{A_k}-1\right)
         \|w_k\|^2\right] \\
        & +\frac{(\mu-\rho)(A_{k+1}-A_k)}{2}\|r_k\|^2 
         -\frac{(\theta L-\rho)(A_{k+1}-A_k)}{\gamma_k}
          \langle v_{y_k}-\beta_k w_k,r_k\rangle \\
        \stackrel{\eqref{eq:combine}}{=}{}&
        \frac{\xi P_{k+1}}{2}\|v_{y_k}\|^2
        +\frac{\xi P_k^2}{2P_{k+1}}
         \left(\frac{qA_{k+1}}{A_k}-1\right)\|w_k\|^2 \\
        & +\frac{P_{k+1}-P_k}{2}\|r_k\|^2
          -P_{k+1}\langle v_{y_k},r_k\rangle
          +P_k\langle w_k,r_k\rangle.
    \end{aligned}
\]
Completing the squares in the terms involving $r_k$ gives
\[
    \begin{aligned}
        N_k ={}&
        \frac{P_{k+1}}{2}\|v_{y_k}-r_k\|^2
        -\frac{P_k}{2}\|w_k-r_k\|^2 \\
        & +\frac{P_{k+1}}{2}(\xi-1)\|v_{y_k}\|^2 
         +\left[
          \frac{\xi P_k^2}{2P_{k+1}}
          \left(\frac{qA_{k+1}}{A_k}-1\right)
          +\frac{P_k}{2}\right]\|w_k\|^2.
    \end{aligned}
\]
The coefficients of $\|v_{y_k}\|^2$ and $\|w_k\|^2$ in the last
two terms are $\sigma_k$ and $t_k$, respectively.

For $k=0$, the initialization $s_0=0$ and $y_0=x_0$ gives
$u_0=w_0=0$. Thus~\eqref{eq:eta} reduces to
$\eta_{y_0}=v_{y_0}/\gamma_0$. Substituting into $N_0$ and
using~\eqref{eq:combine}, we obtain
\[
    \begin{aligned}
        N_0 ={}&
        \frac{\xi P_1}{2}\|v_{y_0}\|^2
        +\frac{P_1-P_0}{2}\|r_0\|^2
        -P_1\langle v_{y_0},r_0\rangle \\
        ={}&
        \frac{P_1}{2}\|v_{y_0}-r_0\|^2
        -\frac{P_0}{2}\|r_0\|^2
        +\sigma_0\|v_{y_0}\|^2.
    \end{aligned}
\]
Since $w_0=0$, the same formula for $N_k$ applies at $k=0$ with
$t_0=0$. Substituting this formula into~\eqref{eq:3-28} proves
the claim for all $k\geqslant0$.
\end{proof}

\subsection{Inequalities for norms across different tangent spaces}
\label{subsec:ineq-norm}

To complete the convergence rate analysis, we use the metric distortion estimates of Kim and Yang~\cite[Lemmas~5.2 and~5.3]{kim2022accelerated} to compare the norms of tangent vectors based at different points. For ease of reference, we state the two estimates together in the following lemma. \Cref{fig:metric distortion} illustrates the corresponding geometric configuration.

\begin{figure}[!htb]
\centering
\resizebox{\linewidth}{!}{%
\begin{tikzpicture}[
  x=0.012cm,y=-0.012cm,
  font=\fontsize{12}{14.4}\selectfont,
  md geometry/.style={draw=black!80,line width=0.95pt,line cap=round},
  md connector/.style={draw=black!80,line width=0.9pt,line cap=round,dash pattern=on 0.3pt off 1pt},
  md vector/.style={line width=0.9pt,line cap=round,-{Latex[length=2.1mm,width=0.95mm]}},
  md dashed/.style={dash pattern=on 2.5pt off 1.5pt},
  md label/.style={inner sep=0pt,outer sep=0pt},
  md formula/.style={md label},
  md point/.style={circle,draw=black!65,fill=mdblue,line width=0.25pt,inner sep=0.85pt}
]
\definecolor{mdblue}{RGB}{17,112,190}
\definecolor{mdgreen}{RGB}{87,170,44}
\definecolor{mdred}{RGB}{255,56,61}
\definecolor{mdcyan}{RGB}{118,207,255}

\begin{scope}
  \coordinate (AX) at (99,375);
  \coordinate (AP) at (323,117);
  \coordinate (AY) at (369,294);
  \coordinate (AZ) at (567,368);
  \coordinate (AB) at (458,118);
  \coordinate (AQ) at (343.12,200.24); 

  \draw[md geometry] (AX) .. controls (142,261) and (229,164) .. (AP);
  \draw[md geometry] (AP) .. controls (438,172) and (523,267) .. (AZ);
  \draw[md geometry] (AX) .. controls (184,314.1) and (275,284.4) .. (AY);
  \draw[md geometry] (AY) .. controls (444,301.7) and (511,326) .. (AZ);
  \draw[md connector] (AP) .. controls (358,160) and (380,231) .. (AY);

  \draw[md vector,draw=black!85] (AX) -- (AQ);
  \draw[md vector,md dashed,draw=mdgreen] (AX) -- (AB);
  \draw[md vector,draw=mdred] (AY) -- (603,318);

  \foreach \point in {AX,AP,AY,AZ} {\node[md point] at (\point) {};}
  \node[md label] at (118,389) {$x$};
  \node[md label] at (305,107) {$p$};
  \node[md label] at (361,314) {$y$};
  \node[md label] at (583,367) {$z$};

  \node[md formula,anchor=west] at (432,164) {$v_x=b=\operatorname{Exp}_x^{-1}(z)$};
  \node[md formula,anchor=west] at (195,255) {$\operatorname{Exp}_x^{-1}(y)=r v_x$};
  \node[md formula,anchor=west] at (480,280) {$v_y=\operatorname{Exp}_y^{-1}(z)$};
  \node[md label] at (380,408) {$(\mathrm{a})\ a=0$};
\end{scope}

\begin{scope}[shift={(620,0)}]
  \coordinate (BX) at (99,375);
  \coordinate (BP) at (323,117);
  \coordinate (BY) at (369,294);
  \coordinate (BZ) at (567,368);
  \coordinate (BB) at (462,118);
  \coordinate (BQ) at (345.84,200.24); 
  \coordinate (BV) at (408,104);

  \draw[md geometry] (BX) .. controls (142,261) and (229,164) .. (BP);
  \draw[md geometry] (BP) .. controls (438,172) and (523,267) .. (BZ);
  \draw[md geometry] (BX) .. controls (184,314.1) and (275,284.4) .. (BY);
  \draw[md geometry] (BY) .. controls (444,301.7) and (511,326) .. (BZ);
  \draw[md connector] (BP) .. controls (358,160) and (380,231) .. (BY);

  \draw[md vector,draw=black!85] (BX) -- (BQ);
  \draw[md vector,draw=black!85] (BX) -- (BV);
  \draw[md vector,md dashed,draw=mdgreen] (BX) -- (BB);
  \draw[md vector,md dashed,draw=mdred] (BQ) -- (BV);
  \draw[md vector,md dashed,draw=mdcyan] (BB) -- (BV);
  \draw[md vector,draw=mdred] (BY) -- (565,268);

  \foreach \point in {BX,BP,BY,BZ} {\node[md point] at (\point) {};}
  \node[md label] at (118,389) {$x$};
  \node[md label] at (305,107) {$p$};
  \node[md label] at (361,314) {$y$};
  \node[md label] at (583,367) {$z$};

  \node[md formula,anchor=west] at (362,90) {$v_x=a+b=a+\operatorname{Exp}_x^{-1}(z)$};
  \node[md label] at (432,124) {$a$};
  \node[md formula,anchor=west] at (435,157) {$\operatorname{Exp}_x^{-1}(z)=b$};
  \node[md formula,anchor=west] at (248,265) {$\operatorname{Exp}_x^{-1}(y)=r b$};
  \node[md formula,anchor=west] at (495,295) {$v_y=\Gamma_x^y\bigl(v_x-\operatorname{Exp}_x^{-1}(y)\bigr)$};
  \node[md label] at (385,408) {$(\mathrm{b})\ a\neq0$};
\end{scope}
\end{tikzpicture}
}
\caption{Geometric configurations underlying \Cref{lem:metric-distortion}. The point $z$ is included solely as a schematic representation of $\e_x(b)$; the lemma does not require $\e_x(b)$ to be defined.}
\label{fig:metric distortion}
\end{figure}

\begin{lemma}
\label{lem:metric-distortion}
Suppose that Assumption~\ref{assum:omega} holds, and let $\zeta$ and $\delta$ be the constants defined in Section~\ref{subsec:zeta-and-delta}. Let $x,y\in\Omega$, $a,b\in\T_x\m$, and $r\in[0,1]$ satisfy $\e_x^{-1}(y)=rb$. Define
$v_x=a+b$,
and 
$v_y=\Gamma_x^y\bigl(v_x-\e_x^{-1}(y)\bigr)$.
Then, for every $p\in\Omega$ and every $\xi\geqslant\zeta$, the following hold:

\begin{enumerate}[(i)]
	\item \label{metric-item1} For $r\in[0,1]$ and $a=0$, we have 
	\begin{align*}
		\left\|v_{y}-\e_{y}^{-1}(p)\right\|^{2}+(\xi-1)\left\|v_{y}\right\|^{2} \leqslant & \left\|v_{x}-\e_{x}^{-1}(p)\right\|^{2}+(\xi-1)\left\|v_{x}\right\|^{2}.  
	\end{align*} 
	\item \label{metric-item2} For $r\in[0,1)$, we have 
	\begin{align*}
		\left\|v_{y}-\e_{y}^{-1}(p)\right\|^{2}+(\xi-1)\left\|v_{y}\right\|^{2} \leqslant & \left\|v_{x}-\e_{x}^{-1}(p)\right\|^{2}+(\xi-1)\left\|v_{x}\right\|^{2}  \notag \\
		& +\frac{\xi-\delta}{2}\left(\frac{1}{1-r}-1\right)\|a\|^{2}. 
	\end{align*}
\end{enumerate} 
\end{lemma}

Applying \Cref{lem:metric-distortion} to the updates in \Cref{algrm:rapg-sc} yields the following inequality relating weighted sums of squared norms in different tangent spaces.

\begin{lemma}\label{lem:ineq-norm}
Suppose that Assumptions~\ref{assum:omega} and~\ref{assum:iterates} hold. Then the quantities generated by \Cref{algrm:rapg-sc} satisfy, for every $k\geqslant0$,
\[
    \begin{aligned}
        & \frac{P_{k+1}}{2}
          \Bigl(
              \bigl\|v_{y_k}-\e_{y_k}^{-1}(x_*)\bigr\|^2
              +(\xi-1)\|v_{y_k}\|^2
          \Bigr) \\
        & \quad -\frac{P_k}{2}
          \Bigl(
              \bigl\|\e_{y_k}^{-1}(z_k)-\e_{y_k}^{-1}(x_*)\bigr\|^2
              +(\xi-1)\bigl\|\e_{y_k}^{-1}(z_k)\bigr\|^2
          \Bigr) \\
        & \geqslant \frac{P_{k+1}}{2}
          \Bigl(
              \bigl\|\e_{x_{k+1}}^{-1}(z_{k+1})-\e_{x_{k+1}}^{-1}(x_*)\bigr\|^2
              +(\xi-1)\bigl\|\e_{x_{k+1}}^{-1}(z_{k+1})\bigr\|^2
          \Bigr) \\
        & \quad -\frac{P_k}{2}
          \Bigl(
              \bigl\|\e_{x_k}^{-1}(z_k)-\e_{x_k}^{-1}(x_*)\bigr\|^2
              +(\xi-1)\bigl\|\e_{x_k}^{-1}(z_k)\bigr\|^2
          \Bigr) 
         +r_k\bigl\|\e_{y_k}^{-1}(z_k)\bigr\|^2,
    \end{aligned}
\]
where
\[
    r_k=
    \begin{cases}
        0, & k=0, \\[2mm]
        \displaystyle
        -\frac{(\xi-\delta)P_{k+1}}{4}\frac{\beta_k^2}{\gamma_k-1}
        +\frac{(\xi-\zeta)P_k}{2}
        \left(\frac{1}{(1-\tau_k)^2}-1\right),
        & k\geqslant1.
    \end{cases}
\]
As before, expressions involving $z_k$ are interpreted through $\e_{x_k}^{-1}(z_k)=s_k$ and $\e_{y_k}^{-1}(z_k)=(1-\tau_k)\Gamma_{x_k}^{y_k}(s_k)$, without assuming that $\e_{x_k}(s_k)$ is defined.
\end{lemma}

\begin{proof}
Let $w_k=\e_{y_k}^{-1}(z_k)$. The formal convention and the algorithmic updates give
\[
    w_k=(1-\tau_k)\Gamma_{x_k}^{y_k}(s_k),
    \qquad
    s_{k+1}=\Gamma_{y_k}^{x_{k+1}}\bigl(v_{y_k}-\eta_{y_k}\bigr).
\]
We first consider $k\geqslant1$. By \Cref{lem:parameters}, $\gamma_k>1$ and $0<\tau_k<1$.

Apply \Cref{lem:metric-distortion}\eqref{metric-item2} with
\[
    \begin{aligned}
        x  = y_k,\quad y=x_{k+1},\quad p=x_*, \quad
        a  = \beta_k w_k,\quad 
		b =\gamma_k\eta_{y_k},\quad r=\frac{1}{\gamma_k}.
    \end{aligned}
\]
The tangent vectors in \Cref{lem:metric-distortion}\eqref{metric-item2} satisfy
    $v_x=v_{y_k}$
    and
    $v_y=\Gamma_{y_k}^{x_{k+1}}\bigl(v_{y_k}-\eta_{y_k}\bigr)=s_{k+1}$.
Consequently,
\begin{equation}\label{eq:6-2s}
    \begin{aligned}
        & \bigl\|s_{k+1}-\e_{x_{k+1}}^{-1}(x_*)\bigr\|^2
          +(\xi-1)\|s_{k+1}\|^2 \\
         \leqslant {}&
          \bigl\|v_{y_k}-\e_{y_k}^{-1}(x_*)\bigr\|^2
          +(\xi-1)\|v_{y_k}\|^2
          +\frac{(\xi-\delta)\beta_k^2}{2(\gamma_k-1)}\|w_k\|^2.
    \end{aligned}
\end{equation}

Next, apply \Cref{lem:metric-distortion}\eqref{metric-item1} with $\zeta$ in place of its parameter $\xi$ and
\[
    \begin{aligned}
        x  = x_k,\quad y=y_k,\quad p=x_*, \quad 
        a  = 0,\quad b=s_k,\quad r=\tau_k.
    \end{aligned}
\]
In this application, $v_x=s_k$ and $v_y=(1-\tau_k)\Gamma_{x_k}^{y_k}(s_k)=w_k$. Hence,
    \begin{align}\label{eq:3-2}
         \bigl\|s_k-\e_{x_k}^{-1}(x_*)\bigr\|^2
          +(\zeta-1)\|s_k\|^2 
         \geqslant 
          \bigl\|w_k-\e_{y_k}^{-1}(x_*)\bigr\|^2
          +(\zeta-1)\|w_k\|^2.
    \end{align}
Since parallel transport preserves norms,
    $\|s_k\|^2=\frac{\|w_k\|^2}{(1-\tau_k)^2}$.
Adding $(\xi-\zeta)\|s_k\|^2$ to both sides of~\eqref{eq:3-2} and using the above identity yields
\begin{equation}\label{eq:7-1s}
    \begin{aligned}
        & \bigl\|s_k-\e_{x_k}^{-1}(x_*)\bigr\|^2
          +(\xi-1)\|s_k\|^2 \\
         \geqslant {}&
          \bigl\|w_k-\e_{y_k}^{-1}(x_*)\bigr\|^2
          +(\xi-1)\|w_k\|^2 
          +(\xi-\zeta)
          \left(\frac{1}{(1-\tau_k)^2}-1\right)\|w_k\|^2.
    \end{aligned}
\end{equation}
Combining~\eqref{eq:6-2s} and~\eqref{eq:7-1s} with the positive weights $P_{k+1}/2$ and $P_k/2$ gives the claimed inequality. Indeed, the resulting coefficient of $\|w_k\|^2$ is
\[
    -\frac{P_{k+1}}{2}\frac{(\xi-\delta)\beta_k^2}{2(\gamma_k-1)}
    +\frac{P_k}{2}(\xi-\zeta)
    \left(\frac{1}{(1-\tau_k)^2}-1\right)
    =r_k.
\]

It remains to consider $k=0$. The initialization gives $s_0=w_0=0$ and $y_0=x_0$, so $v_{y_0}=\gamma_0\eta_{y_0}$. Apply \Cref{lem:metric-distortion}\eqref{metric-item1} with
\[
    \begin{aligned}
        x  = y_0,\quad y=x_1,\quad p=x_*, \quad
        a  = 0,\quad b=v_{y_0},\quad r=\frac{1}{\gamma_0}\in(0,1].
    \end{aligned}
\]
The transported vector is again $s_1$, and therefore
\begin{equation}\label{eq:6-2s0}
    \begin{aligned}
        & \bigl\|s_1-\e_{x_1}^{-1}(x_*)\bigr\|^2
          +(\xi-1)\|s_1\|^2 
         \leqslant 
          \bigl\|v_{y_0}-\e_{y_0}^{-1}(x_*)\bigr\|^2
          +(\xi-1)\|v_{y_0}\|^2.
    \end{aligned}
\end{equation}
Moreover, $s_0=w_0=0$ and $x_0=y_0$ imply
\begin{equation}\label{eq:7-1s0}
    \begin{aligned}
        & \bigl\|s_0-\e_{x_0}^{-1}(x_*)\bigr\|^2
          +(\xi-1)\|s_0\|^2 
         = 
          \bigl\|w_0-\e_{y_0}^{-1}(x_*)\bigr\|^2
          +(\xi-1)\|w_0\|^2.
    \end{aligned}
\end{equation}
Combining~\eqref{eq:6-2s0} and~\eqref{eq:7-1s0} with weights $P_1/2$ and $P_0/2$ proves the result with $r_0=0$, including the boundary cases $\gamma_0=1$ and $\tau_0=1$.
\end{proof}

\subsection{Convergence analysis}
\label{subsec:main-results}

For \Cref{algrm:rapg-sc}, we can now introduce a potential function defined as
\begin{equation}
\label{eq:potential-sc-1}
\begin{aligned}
\Phi_k
={}&A_k\bigl(F(x_k)-F(x_*)\bigr)\\
&+\frac{P_k}{2}
\Bigl(
\bigl\|\e_{x_k}^{-1}(z_k)-\e_{x_k}^{-1}(x_*)\bigr\|^2
+(\xi-1)\bigl\|\e_{x_k}^{-1}(z_k)\bigr\|^2
\Bigr),
\end{aligned}
\end{equation}
where $P_k=(\theta L-\rho)(\xi+qA_k)$, as in~\eqref{eq:Pk}, and $\e_{x_k}^{-1}(z_k)=s_k$ is understood as formal notation, as before. By \Cref{lem:parameters} and the parameter conditions in \Cref{algrm:rapg-sc}, we have $A_k\geqslant0$ and $P_k>0$ for all $k\geqslant0$.

\begin{remark}
The recursions for $A_k$, $\beta_k$, $\gamma_k$, and $\tau_k$ in \Cref{algrm:rapg-sc} are derived by requiring $\Phi_{k+1}\leqslant\Phi_k$. In the proof of \Cref{thm:potential-decrease} below, we impose two conditions: (i) the coefficients of inner product terms of indeterminate sign vanish, and (ii) the coefficients of the additional squared norm terms have a common sign that allows these terms to be discarded in establishing $\Phi_{k+1}\leqslant\Phi_k$.
Since $\xi\geqslant1$ and $P_k>0$, we have $\Phi_k\geqslant A_k\bigl(F(x_k)-F(x_*)\bigr)$. Combining this with $\Phi_k\leqslant\Phi_0$ yields
\(F(x_k)-F(x_*)\leqslant\frac{\Phi_0}{A_k}\)
whenever $A_k>0$,
and hence a convergence rate of $\mathcal{O}(1/A_k)$.
\end{remark}

\subsubsection{Decrease of the potential function}
In \Cref{thm:potential-decrease}, we establish the decrease of the potential function $\Phi_k$ under Condition~\eqref{eq:D11}. Parameter choices ensuring this condition at every iteration are given in \Cref{prop:cond-sc}. For convenience, define
\begin{equation}\label{eq:Gk}
    G_k=1+\frac{qA_k}{\xi},\qquad k\geqslant0.
\end{equation}
Since $A_k\geqslant0$, $q\geqslant0$, and $\xi\geqslant1$, we have $G_k\geqslant1$ for all $k\geqslant0$.

\begin{theorem}[Decrease of the potential function]\label{thm:potential-decrease}
Suppose that Assumptions~\ref{assum:omega}--\ref{assum:iterates} hold. Then the potential function $\Phi_k$ defined in~\eqref{eq:potential-sc-1} satisfies $\Phi_{k+1}\leqslant\Phi_k$ for all $k\geqslant0$, provided that the following condition holds for every $k\geqslant1$:
\begin{equation}\label{eq:D11}
    H_{k+1}
    =4(\xi-\zeta)\bigl(\sqrt{A_{k+1}}-\xi\sqrt{G_{k+1}}\bigr)
    -(\xi-\delta)\bigl(\sqrt{A_{k+1}}-\sqrt{G_{k+1}}\bigr)
    \geqslant0.
\end{equation}
\end{theorem}

\begin{proof}
Combining \Cref{lem:ineq-function,lem:ineq-norm} and using $\sigma_k=\frac{P_{k+1}}{2}(\xi-1)$, we obtain
\[
    \Phi_k-\Phi_{k+1}
    \geqslant b_k\bigl\|\e_{y_k}^{-1}(z_k)\bigr\|^2,
    \qquad k\geqslant0,
\]
where
\begin{equation}\label{eq:b_2}
    b_k=r_k+\frac{P_k}{2}(\xi-1)+t_k.
\end{equation}
At $k=0$, we have $r_0=t_0=0$ and $\e_{y_0}^{-1}(z_0)=0$ by initialization. Thus, $\Phi_1\leqslant\Phi_0$.

Now fix $k\geqslant1$. By \Cref{lem:parameters}, we have $A_k>0$, $\beta_k>0$, $\gamma_k>1$, and $\tau_k\in(0,1)$. Moreover, \eqref{eq:Pk} gives $P_{k+1}\geqslant P_k>0$. Hence, the definition of $t_{k}   = \frac{\xi P_{k}^{2}}{2P_{k+1}}\left( \frac{qA_{k+1}}{A_{k} } - 1 \right)  + \frac{P_{k} }{2}$ yields
\[
    \begin{aligned}
        & \frac{P_k}{2}(\xi-1)+t_k 
        = \frac{\xi P_k^2}{2P_{k+1}}
              \left(
                  \frac{qA_{k+1}}{A_k}
                  +\frac{P_{k+1}-P_k}{P_k}
              \right)
              \geqslant 0.
    \end{aligned}
\]
It therefore remains to show that $r_k\geqslant0$.

For $k\geqslant1$, the definition of $r_k$ gives
\begin{equation}\label{eq:r_k}
    r_k
    =-\frac{(\xi-\delta)P_{k+1}}{4}
      \frac{\beta_k^2}{\gamma_k-1}
    +\frac{(\xi-\zeta)P_k}{2}
      \left(\frac{1}{(1-\tau_k)^2}-1\right).
\end{equation}
By~\eqref{eq:Ak} and the definitions of $\beta_k$, $\gamma_k$, and $G_k$, we have
\begin{subequations}
    \begin{empheq}[left=\empheqlbrace]{align}
        A_k
        &\overset{\eqref{eq:Ak}\eqref{eq:Gk}}{=} A_{k+1}-\sqrt{A_{k+1}G_{k+1}},
        \label{eq:10-1}\\
        \beta_k
        &\overset{\eqref{eq:beta_k}\eqref{eq:Pk}}{=} \frac{P_k}{P_{k+1}},
        \label{eq:A13}\\
        \gamma_k
        &\overset{\eqref{eq:gamma_k}\eqref{eq:Gk}}{=} \frac{A_{k+1}-A_k}{\xi G_{k+1}}
        \overset{\eqref{eq:10-1}}{=} \frac{\sqrt{A_{k+1}}}{\xi\sqrt{G_{k+1}}}.
        \label{eq:8-11}
    \end{empheq}
\end{subequations}
Also, \eqref{eq:tau_k} implies
\begin{equation}\label{eq:A15}
    \frac{1}{(1-\tau_k)^2}-1
    =2\frac{\beta_kA_{k+1}}{\gamma_kA_k}
     +\left(\frac{\beta_kA_{k+1}}{\gamma_kA_k}\right)^2
    \geqslant2\frac{\beta_kA_{k+1}}{\gamma_kA_k}.
\end{equation}
Substituting~\eqref{eq:A13} and~\eqref{eq:A15} into~\eqref{eq:r_k} and using $\xi\geqslant\zeta$, we have
\[
    \begin{aligned}
        r_k
        & \geqslant
          -\frac{(\xi-\delta)\beta_kP_k}{4(\gamma_k-1)}
          +\frac{(\xi-\zeta)\beta_kP_kA_{k+1}}{\gamma_kA_k}\\
        & = \frac{\beta_kP_k}{4\gamma_k(\gamma_k-1)A_k}
            \Bigl[
                4(\xi-\zeta)A_{k+1}(\gamma_k-1)
                -(\xi-\delta)\gamma_kA_k
            \Bigr]\\
        & = \frac{\beta_kP_k\sqrt{A_{k+1}}}
                 {4(\gamma_k-1)A_k}
            H_{k+1}
          \geqslant 0.
    \end{aligned}
\]
%
The last equality follows from~\eqref{eq:10-1}, \eqref{eq:8-11}, and~\eqref{eq:D11}, and the last inequality from $H_{k+1}\geqslant 0$ and the positivity of the prefactor. Hence $b_k\geqslant 0$ by~\eqref{eq:b_2}, and therefore $\Phi_{k+1}\leqslant\Phi_k$.
\end{proof}

The following proposition gives explicit parameter conditions for~\eqref{eq:D11} to hold. Conditions~\eqref{cond:a} and~\eqref{cond:b} apply in Euclidean space and, more generally, on flat manifolds, where the sectional curvature is zero everywhere. Condition~\eqref{cond:c} allows for nonzero sectional curvature.

\begin{proposition}[Sufficient parameter conditions]
\label{prop:cond-sc}
Suppose that the parameters satisfy the requirements of \Cref{algrm:rapg-sc}. Then~\eqref{eq:D11} holds for every $k\geqslant 1$ in each of the following cases.
\begin{enumerate}[(i)]
    \item \label{cond:a} $\zeta=\delta=\xi=1$.

    \item \label{cond:b} $\zeta=\delta=1$ and $\xi>1$. Let $\theta\geqslant1$ satisfy
        $\theta>\frac{\rho}{L}
        +\frac{(4\xi-1)^2(\mu-\rho)}{9L\xi}$.
    With $N_1=9-(4\xi-1)^2q/\xi$, assume that
    \begin{equation}\label{eq:A0-case2}
        A_0\geqslant\frac{4(4\xi-1)(\xi-1)}{N_1}.
    \end{equation}

    \item \label{cond:c} $\zeta>\delta$ and $\xi\geqslant\zeta+\frac{\zeta-\delta}{\lambda-1}$ for some $\lambda\in(1,4)$. Let $\theta\geqslant1$ satisfy
        $\theta>\frac{\rho}{L}
        +\frac{(4\xi-\lambda)^2(\mu-\rho)}{(4-\lambda)^2L\xi}$.
    With $N_2=(\tfrac4\lambda-1)^2-(\tfrac{4\xi}{\lambda}-1)^2q/\xi$, assume that
    \begin{equation}\label{eq:A0-case3}
        A_0\geqslant\frac{4(4\xi-\lambda)(\xi-1)}{\lambda^2N_2}.
    \end{equation}
\end{enumerate}
\end{proposition}

\begin{proof}
Fix $k\geqslant 1$. In Case~\eqref{cond:a}, both terms in $H_{k+1}$ vanish. For the remaining cases, set
\[
    (a,b)=
    \begin{cases}
        (3,\,4\xi-1), & \text{in Case~\eqref{cond:b}},\\[1mm]
        \left(\tfrac4\lambda-1,\,\tfrac{4\xi}{\lambda}-1\right),
        & \text{in Case~\eqref{cond:c}}.
    \end{cases}
\]
In both cases, $b>a>0$, and it suffices to establish
\begin{equation}\label{eq:cond-sc-sqrt}
    a\sqrt{A_{k+1}}\geqslant b\sqrt{G_{k+1}}.
\end{equation}
Indeed, in Case~\eqref{cond:b}, we have $H_{k+1}=(\xi-1)(a\sqrt{A_{k+1}}-b\sqrt{G_{k+1}})$. In Case~\eqref{cond:c}, the bound on $\xi$ gives $\lambda(\xi-\zeta)\geqslant\xi-\delta>0$. Since $b-a\xi=\xi-1\geqslant0$, \eqref{eq:cond-sc-sqrt} also implies $\sqrt{A_{k+1}}-\xi\sqrt{G_{k+1}}\geqslant0$. Consequently,
\[
    \begin{aligned}
        \frac{H_{k+1}}{\xi-\delta}
        & \geqslant \frac4\lambda
            \left(\sqrt{A_{k+1}}-\xi\sqrt{G_{k+1}}\right)
            -\left(\sqrt{A_{k+1}}-\sqrt{G_{k+1}}\right) 
         = a\sqrt{A_{k+1}}-b\sqrt{G_{k+1}}.
    \end{aligned}
\]

Using $G_{k+1}=1+qA_{k+1}/\xi$, we can write~\eqref{eq:cond-sc-sqrt} equivalently as $NA_{k+1}\geqslant b^2$, where $N=a^2-b^2q/\xi$. In both cases, the strict bound on $\theta$ takes the form $\theta L-\rho>b^2(\mu-\rho)/(a^2\xi)$, which ensures $N>0$.

Suppose that $A_{k+1}<b^2/N$. Then $a\sqrt{A_{k+1}}<b\sqrt{G_{k+1}}$. Since $A_{k+1}>0$, the recurrence~\eqref{eq:Ak} gives
\[
    \begin{aligned}
        A_k
         = A_{k+1}\left(1-\sqrt{\frac{G_{k+1}}{A_{k+1}}}\right)
         < A_{k+1}\left(1-\frac{a}{b}\right)
          < \frac{b(b-a)}{N}.
    \end{aligned}
\]
Because $A_k\geqslant A_0$, this is impossible if $A_0\geqslant b(b-a)/N$. Substituting the corresponding values of $a$, $b$, and $N$ gives~\eqref{eq:A0-case2} and~\eqref{eq:A0-case3}.
\end{proof}

\begin{remark}[Choice of $\theta$ and $A_0$]
\label{rem:q-branches}
When $q=0$, we have $\mu=\rho$. Since $L>\rho$, the smallest admissible choice $\theta=1$ satisfies the bounds on $\theta$ in \Cref{prop:cond-sc}. In this case, $N_1=9$ and $N_2=(\tfrac4\lambda-1)^2$, so~\eqref{eq:A0-case3} reduces to
\begin{equation}\label{eq:A0-q0}
    A_0\geqslant
    \frac{4(4\xi-\lambda)(\xi-1)}{(4-\lambda)^2}.
\end{equation}
When $q>0$, the bound in~\eqref{eq:A0-case3} can be large if $N_2$ is close to zero. Fix $\nu\in(0,1)$ and impose the margin $N_2\geqslant\nu(\tfrac4\lambda-1)^2$. Substituting $q=(\mu-\rho)/(\theta L-\rho)$ into this margin gives the bound on $\theta$ below. Using the margin in~\eqref{eq:A0-case3} gives the corresponding sufficient bound on $A_0$.
\begin{equation}\label{eq:theta-nu}
    \begin{aligned}
        \theta
         \geqslant \max\left\{
            1,\,
            \frac{\rho}{L}
            +\frac{(4\xi-\lambda)^2(\mu-\rho)}
            {(1-\nu)(4-\lambda)^2L\xi}
        \right\},\qquad 
        A_0
         \geqslant
        \frac{4(4\xi-\lambda)(\xi-1)}
        {\nu(4-\lambda)^2}.
    \end{aligned}
\end{equation}
For fixed $\xi$ and $\lambda$, increasing $\nu$ lowers the sufficient bound on $A_0$ and can raise the bound on $\theta$. We recommend $\nu=\tfrac{1}{5}$ as the default value.
\end{remark}

\begin{remark}[Choice of $\xi$ and $\lambda$]
\label{rem:xi-lambda}
Suppose that $\zeta>\delta$. For each $\lambda\in(1,4)$, choose $\xi(\lambda)=\zeta+\frac{\zeta-\delta}{\lambda-1}$, the smallest value satisfying the lower bound on $\xi$ in Case~\eqref{cond:c} of \Cref{prop:cond-sc}.
After substituting $\xi=\xi(\lambda)$ into~\eqref{eq:A0-q0}, the derivative of its right-hand side with respect to $\lambda$ has the sign of
   $ (\zeta-1)(\lambda-1)^2+(\zeta-\delta)(2\lambda-5)$.
This polynomial is strictly increasing on $(1,4)$ and has opposite signs at the endpoints. Its unique zero therefore gives the minimizing choice
\[
    \lambda_*=1+\frac{3}{1+\sqrt{1+\frac{3(\zeta-1)}{\zeta-\delta}}},
    \qquad \xi_*=\xi(\lambda_*).
\]
For fixed $\nu\in(0,1)$, the lower bound on $A_0$ in~\eqref{eq:theta-nu} is $\frac{1}{\nu}$ times that in~\eqref{eq:A0-q0}$ $ and hence has the same minimizer. Since $0\leqslant\zeta-1\leqslant\zeta-\delta$, we have $\lambda_*\in[2,\frac{5}{2}]$, with
\[
    (\lambda_*,\xi_*)=
    \begin{cases}
        (2,\,2\zeta-1), & \text{if }\delta=1,\ \zeta>1,\\[1mm]
        \left(\frac{5}{2},\,\frac{5-2\delta}{3}\right), & \text{if }\zeta=1,\ \delta<1.
    \end{cases}
\]
As $\lambda\downarrow1$, the two $A_0$ bounds diverge because $\xi(\lambda)\to+\infty$. As $\lambda\uparrow4$, they diverge because their denominators tend to zero while their numerators have positive limits, since $\xi(\lambda)\to\zeta+\frac{\zeta-\delta}{3}>1$. When $\mu>\rho$, the lower bound on $\theta$ in~\eqref{eq:theta-nu} also diverges at both endpoints.
We therefore recommend moderate values, such as $\lambda\in[2,\frac{5}{2}]$.
\end{remark}

\subsubsection{Convergence rates}

We first bound the growth of $\{A_k\}$ to express the objective residual estimate in terms of the iteration index.

\begin{lemma}\label{lem:growth rate of Ak}
Let $\{A_k\}$ be the sequence generated by \Cref{algrm:rapg-sc}. Then, for every $k\geqslant1$,
\[
    A_k\geqslant\left(\sqrt{A_0}+\frac{k}{2}\right)^2
    \qquad\text{and}\qquad
    A_k\geqslant\frac{A_0+\frac12}{\left(1-\sqrt{\frac{q}{\xi}}\right)^k}.
\]
\end{lemma}

\begin{proof}
By~\eqref{eq:Ak} and~\eqref{eq:Gk}, we have $A_{k+1}-A_k=\sqrt{A_{k+1}G_{k+1}}$ for every $k\geqslant0$. Since $G_{k+1}\geqslant1$ and $A_{k+1}>A_k\geqslant0$, it follows that
\[
    \sqrt{A_{k+1}}-\sqrt{A_k}
    =\frac{A_{k+1}-A_k}{\sqrt{A_{k+1}}+\sqrt{A_k}}
    \geqslant\frac{\sqrt{A_{k+1}}}{\sqrt{A_{k+1}}+\sqrt{A_k}}
    \geqslant\frac{1}{2}.
\]
Summing these inequalities gives $\sqrt{A_k}\geqslant\sqrt{A_0}+\frac{k}{2}$, which yields the first bound.

For the second bound, rationalizing the difference $A_{k+1}-A_k-\sqrt{\frac{q}{\xi}}A_{k+1}$ gives
\[
    \begin{aligned}
       \left(1-\sqrt{\frac{q}{\xi}}\right)A_{k+1}-A_k 
	   ={} &
	   A_{k+1}-A_k-\sqrt{\frac{q}{\xi}}A_{k+1} 
        = \frac{(A_{k+1}-A_k)^2 - \frac{q}{\xi} A_{k+1}^2 }{A_{k+1}-A_k+\sqrt{\frac{q}{\xi}}A_{k+1}} \\
        \stackrel{\eqref{eq:A_k+1}}{=} & \frac{A_{k+1}}{A_{k+1}-A_k+\sqrt{\frac{q}{\xi}}A_{k+1}} 
        \geqslant \frac{1}{1+\sqrt{\frac{q}{\xi}}}
        \geqslant\frac{1}{2},
    \end{aligned}
\]
where the inequalities follow from $A_{k+1}>A_k\geqslant 0$ and $0\leqslant q<\xi$. At $k=0$, this gives $A_1\geqslant\frac{A_0+\frac{1}{2}}{1-\sqrt{\frac{q}{\xi}}}$. Iterating $A_{k+1} > \frac{A_k}{1-\sqrt{\frac{q}{\xi}}}$ from $A_1$ then yields the result.
\end{proof}

The decrease of the potential function converts these growth estimates into bounds on the objective residual.

\begin{theorem}[Convergence rate]\label{thm:rate-rapg}
Suppose that Assumptions~\ref{assum:omega}--\ref{assum:iterates} hold and that the parameters satisfy one of the conditions in \Cref{prop:cond-sc}. Then, for every $k\geqslant1$, the iterates generated by \Cref{algrm:rapg-sc} satisfy
\[
    F(x_k)-F(x_*)
    \leqslant C\min\left\{
        \frac{2}{(2\sqrt{A_0}+k)^2},\,
        \frac{1}{1+2A_0}\left(1-\sqrt{\frac{q}{\xi}}\right)^k
    \right\},
\]
where
\(
    C=2A_0\bigl(F(x_0)-F(x_*)\bigr)
    +(\theta L-\rho)(\xi+qA_0)\bigl\|\e_{x_0}^{-1}(x_*)\bigr\|^2.
\)
\end{theorem}

\begin{proof}
By \Cref{prop:cond-sc}, condition~\eqref{eq:D11} holds for every $k\geqslant1$. Thus, \Cref{thm:potential-decrease} gives $\Phi_k\leqslant\Phi_0$. Since the quadratic terms in~\eqref{eq:potential-sc-1} are nonnegative and $A_k>0$ for $k\geqslant1$, we have
\begin{equation}\label{eq:59}
    F(x_k)-F(x_*)\leqslant\frac{\Phi_k}{A_k}\leqslant\frac{\Phi_0}{A_k},
    \qquad k\geqslant1.
\end{equation}
Using $s_0=0$ and $P_0=(\theta L-\rho)(\xi+qA_0)$ in the definition of $\Phi_0$ gives $2\Phi_0=C$. Substituting the two bounds from \Cref{lem:growth rate of Ak} into~\eqref{eq:59} yields the stated estimate.
\end{proof}

\begin{remark}\label{remark:conv}
The first term in \Cref{thm:rate-rapg} gives an $\mathcal{O}(k^{-2})$ bound on the objective residual for all $\mu\geqslant\rho$. When $\mu>\rho$, the second term gives a geometric bound with factor $1-\sqrt{\frac{q}{\xi}}$. The parameter $\xi$ can be chosen using only the geometric constants $\zeta$ and $\delta$, as in \Cref{rem:xi-lambda}. The dependence of the geometric factor on the objective enters through $q=\frac{\mu-\rho}{\theta L-\rho}$. See \Cref{remark:theta} for settings in which $\theta$ can be chosen independently of the objective function.
\end{remark}

\section{An adaptive restart variant for convex and nonconvex optimization}
\label{sec:restart}

This section develops an adaptive restart variant with proximal gradient safeguards and analyzes its convergence and iteration complexity.

\subsection{Adaptive restart and safeguard mechanisms}

Many manifold optimization problems are not globally geodesically convex, and verifying local convexity and estimating the associated parameters can be difficult. Moreover, the requirement in \Cref{algrm:rapg-sc} that the iterates remain in a prescribed region may be difficult to enforce. These issues motivate the adaptive restart variant presented in \Cref{algrm:safeguard-rapg2}, which combines restart ideas~\cite{odonoghue2015adaptive} with a proximal gradient safeguard adapted from~\cite[Algorithm~5]{huang2022riemannian}.

The accelerated updates follow Steps~\ref{alg:st6}--\ref{alg:s_k+1} of \Cref{algrm:rapg-sc}. The acceleration parameters $(A_{k+1},\beta_k,\gamma_k,\tau_k)$ may be chosen by any rule that makes these updates well defined.

The safeguard in \Cref{alg:Safeguard} is invoked after each block of $N_i$ accelerated iterations. Starting from the reference point $\tilde{x}_i$, it computes a proximal gradient direction $\eta_{\tilde{x}_i}$ and backtracks from $\alpha_i=1$ to enforce
\begin{equation}
\label{eq:ar-safeguard-decrease}
    F\bigl(\e_{\tilde{x}_i}(\alpha_i\eta_{\tilde{x}_i})\bigr) \leqslant F(\tilde{x}_i)-\sigma\alpha_i\|\eta_{\tilde{x}_i}\|^{2}.
\end{equation}
If this condition still fails after $N_{\mathrm{ls}}$ backtracking reductions, $L$ is increased, the direction is recomputed, and backtracking restarts at the same reference point. Thus, the exact Lipschitz constant need not be known in advance. The resulting safeguard point is $z_i=\e_{\tilde{x}_i}(\alpha_i\eta_{\tilde{x}_i})$.

If $F(z_i)<F(x_k)$, the safeguard replaces $x_k$ by $z_i$ and restarts the accelerated iteration by resetting $s_k=0_{x_k}$ and $A_k=A_0$; it also increases $L$ when $N_i<N_{\max}$. Otherwise, $x_k$ is retained, with the same reset performed only if backtracking has increased $L$. In either case, the next reference point is $\tilde{x}_{i+1}=x_k$, so that
   $ F(\tilde{x}_{i+1}) \leqslant F(z_i) \leqslant F(\tilde{x}_i)$.

To limit the additional proximal gradient computations, the safeguard interval is shortened by one after an accepted replacement and lengthened by one otherwise, subject to the bounds $N_{\min}$ and $N_{\max}$.

\begin{algorithm}[!htb]
\caption{Riemannian accelerated proximal gradient method with adaptive restart}
\label{algrm:safeguard-rapg2}
\begin{algorithmic}[1]
    \Require Initial point $x_0\in\mathcal M$; integers $1\leqslant N_{\min}\leqslant N_0\leqslant N_{\max}$ specifying the safeguard intervals; initial Lipschitz estimate $L_{\mathrm{init}}>0$; growth factor $\tau_L>1$; line-search parameters $\sigma,\iota\in(0,1)$; maximum number of backtracking reductions per direction $N_{\mathrm{ls}}\in\mathbb N$; parameters $\theta\geqslant1$ and $A_0\geqslant0$.
    \State Initialize $\tilde{x}_0=x_0$, $s_0=0_{x_0}$, $L=L_{\mathrm{init}}$, $i=0$, and $j=N_0$. \Comment{$i$: safeguard count; $j$: next safeguard iteration}
    \For{$k=0,1,2,\ldots$}
        \If{$k=j$} \label{reARPG:st3}
            \State $(\tilde{x}_{i+1},x_k,s_k,A_k,N_{i+1},L)\gets\text{\Cref{alg:Safeguard}}(\tilde{x}_i,x_k,s_k,A_k,N_i,L)$.
            \State Set $j=j+N_{i+1}$ and $i=i+1$. \label{reARPG:st4}
        \EndIf \label{reARPG:st6}
        \State Choose $(A_{k+1},\beta_k,\gamma_k,\tau_k)$ using the current $L$, e.g., by the recurrence in Step~\ref{alg:st3} of \Cref{algrm:rapg-sc} with inputs for which it is well defined. \label{reARPG:st8}
        \State Compute $(y_k,x_{k+1},s_{k+1})$ according to Steps~\ref{alg:st6}--\ref{alg:s_k+1} of \Cref{algrm:rapg-sc}. \label{reARPG:st9}
    \EndFor
\end{algorithmic}
\end{algorithm}

\begin{algorithm}[!htb]
\caption{Adaptive safeguard procedure}
\label{alg:Safeguard}
\begin{algorithmic}[1]
    \Require $(\tilde{x}_i,x_k,s_k,A_k,N_i,L)$.
    \Ensure $(\tilde{x}_{i+1},x_k,s_k,A_k,N_{i+1},L)$.
    \State Set $\tilde{L}=L$.
    \Repeat
        \State Using the current $L$, compute a stationary point $\eta_{\tilde{x}_i}\in\mathrm{T}_{\tilde{x}_i}\mathcal M$ of $\ell_{\tilde{x}_i}$ satisfying $\ell_{\tilde{x}_i}(\eta_{\tilde{x}_i})\leqslant\ell_{\tilde{x}_i}(0_{\tilde{x}_i})$. \label{alg:Safeguard:st1}
        \State Set $\alpha_i=1$ and $i_{\mathrm{ls}}=0$.
        \While{$i_{\mathrm{ls}}<N_{\mathrm{ls}}$ \textbf{and} \eqref{eq:ar-safeguard-decrease} fails} \label{alg:Safeguard:st3}
            \State Set $\alpha_i=\iota\alpha_i$ and $i_{\mathrm{ls}}=i_{\mathrm{ls}}+1$.
        \EndWhile \label{alg:Safeguard:st5}
        \State Set $z_i=\e_{\tilde{x}_i}(\alpha_i\eta_{\tilde{x}_i})$.
        \If{\eqref{eq:ar-safeguard-decrease} fails} \label{alg:Safeguard:st6}
            \State Set $L=\tau_L L$. \label{alg:Safeguard:st8}
        \EndIf
    \Until{\eqref{eq:ar-safeguard-decrease} holds}
    \If{$F(z_i)<F(x_k)$} \label{alg:Safeguard:st9}
	\Comment{Accept the safeguard point}
        \If{$N_i<N_{\max}$}
            \State Set $L=\tau_L L$. \label{alg:Safeguard:st12}
        \EndIf
        \State Set $x_k=z_i$ and reset $s_k=0_{x_k}$ and $A_k=A_0$. \label{alg:Safeguard:st13}
		 \Comment{Replace the iterate and restart}
        \State Set $N_{i+1}=\max\{N_i-1,N_{\min}\}$. \label{alg:Safeguard:st14}
    \Else
        \If{$L>\tilde{L}$}
            \State Reset $s_k=0_{x_k}$ and $A_k=A_0$. \Comment{Restart at the retained iterate}
        \EndIf
        \State Set $N_{i+1}=\min\{N_i+1,N_{\max}\}$. \label{alg:Safeguard:st17}
    \EndIf
    \State Set $\tilde{x}_{i+1}=x_k$. \label{alg:Safeguard:st19}
\end{algorithmic}
\end{algorithm}

\subsection{Convergence and complexity analysis}
\label{subsec:global-convergence}

Let $R$ be a smooth retraction on $\m$. For $\eta_x\in\mathrm{T}_x\m$ and $y=R_x(\eta_x)$, let $\mathcal T_{\eta_x}=\mathrm{D}R_x(\eta_x)$ denote the differentiated retraction from $\mathrm{T}_x\m$ to $\mathrm{T}_y\m$. Its metric adjoint $\mathcal T_{\eta_x}^\sharp$ maps $\mathrm{T}_y\m$ to $\mathrm{T}_x\m$ and satisfies
\(
    \langle\mathcal T_{\eta_x}^\sharp v_y,u_x\rangle_x=\langle v_y,\mathcal T_{\eta_x}u_x\rangle_y
\)
for any \(u_x\in\mathrm{T}_x\m\) and \( v_y\in\mathrm{T}_y\m\).
When $\mathcal T_{\eta_x}^\sharp$ is invertible, let $\mathcal T_{\eta_x}^{-\sharp}=(\mathcal T_{\eta_x}^\sharp)^{-1}$. In the safeguard analysis below, we take $R=\e$.

We analyze \Cref{algrm:safeguard-rapg2} under Assumption~\ref{assum:f-l-smooth,h-lipschitz} with $\Omega=\m$.   
For the $i$th call to \Cref{alg:Safeguard}, let $L_i$ denote the value of $L$ used to compute the direction $\eta_{\tilde{x}_i}$ from which $z_i$ is formed. This is the value at the end of the line search, before the comparison of $F(z_i)$ with $F(x_k)$. 
We first establish  descent estimates for the safeguard step with an adaptive smoothness estimate.

\begin{lemma}\label{le2}
Suppose that Assumption~\ref{assum:f-l-smooth,h-lipschitz} holds with $\Omega=\m$. For \Cref{algrm:safeguard-rapg2}, each safeguard search terminates after finitely many increases of $L$. Moreover,
\begin{enumerate}[(i)]
    \item \label{le2:1}
    At every safeguard call,
    \(
        F(\tilde{x}_{i+1})\leqslant F(z_i)\leqslant F(\tilde{x}_i)-\sigma\iota^{N_{\mathrm{ls}}}\|\eta_{\tilde{x}_i}\|^2.
    \)
    \item \label{le2:2}
    If $\eta_{\tilde{x}_i}=0_{\tilde{x}_i}$, then $\tilde{x}_i$ is stationary for Problem~\eqref{eq:f+h}.
    \item \label{le2:3}
    If $F$ is bounded below, then
    \(
        \lim_{i\to\infty}\|\eta_{\tilde{x}_i}\|=0.
    \)
\end{enumerate}
\end{lemma}

\begin{proof}
    Let $L_f$ be a geodesic smoothness constant of $f$, and let $\eta_{\tilde{x}_i}$ be the safeguard direction computed with the trial estimate $L$. Geodesic smoothness and the decrease condition $\ell_{\tilde{x}_i}(\eta_{\tilde{x}_i})\leqslant\ell_{\tilde{x}_i}(0_{\tilde{x}_i})$ yield
\[
    \begin{aligned}
        F\bigl(\e_{\tilde{x}_i}(\eta_{\tilde{x}_i})\bigr)-F(\tilde{x}_i)
        & \leqslant \ell_{\tilde{x}_i}(\eta_{\tilde{x}_i})-\ell_{\tilde{x}_i}(0_{\tilde{x}_i})-\frac{\theta L-L_f}{2}\|\eta_{\tilde{x}_i}\|^2 
        \leqslant -\frac{\theta L-L_f}{2}\|\eta_{\tilde{x}_i}\|^2.
    \end{aligned}
\]
Thus the full step $\alpha_i=1$ satisfies the sufficient decrease condition whenever $\theta L\geqslant L_f+2\sigma$. Since each unsuccessful search increases $L$ by the factor $\tau_L>1$, the safeguard search terminates after finitely many increases.

At termination, $\alpha_i\geqslant\iota^{N_{\mathrm{ls}}}$, and the sufficient decrease condition gives
\[
    F(z_i)\leqslant F(\tilde{x}_i)-\sigma\alpha_i\|\eta_{\tilde{x}_i}\|^2
    \leqslant F(\tilde{x}_i)-\sigma\iota^{N_{\mathrm{ls}}}\|\eta_{\tilde{x}_i}\|^2.
\]
By the safeguard acceptance rule, $\tilde{x}_{i+1}=z_i$ if $F(z_i)<F(x_k)$ and $\tilde{x}_{i+1}=x_k$ otherwise. Hence $F(\tilde{x}_{i+1})\leqslant F(z_i)$, proving~\eqref{le2:1}.

Since $\eta_{\tilde{x}_i}$ is a Clarke stationary point of $\ell_{\tilde{x}_i}$, the smooth sum rule and the Clarke chain rule give
\begin{equation}\label{e22}
    0\in\grad f(\tilde{x}_i)+\theta L_i\eta_{\tilde{x}_i}+\mathcal T_{\eta_{\tilde{x}_i}}^\sharp\partial h\bigl(\e_{\tilde{x}_i}(\eta_{\tilde{x}_i})\bigr).
\end{equation}
If $\eta_{\tilde{x}_i}=0_{\tilde{x}_i}$, then $\e_{\tilde{x}_i}(0_{\tilde{x}_i})=\tilde{x}_i$ and $\mathcal T_{0_{\tilde{x}_i}}^\sharp=\mathrm{id}$. Thus~\eqref{e22} reduces to $0\in\grad f(\tilde{x}_i)+\partial h(\tilde{x}_i)$, proving~\eqref{le2:2}.

Finally, if $F$ is bounded below, then~\eqref{le2:1} implies that $\{F(\tilde{x}_i)\}$ converges. Consequently,
\[
    0\leqslant\sigma\iota^{N_{\mathrm{ls}}}\|\eta_{\tilde{x}_i}\|^2
    \leqslant F(\tilde{x}_i)-F(\tilde{x}_{i+1})\to 0,
\]
which proves~\eqref{le2:2}. 
\end{proof}

The following lemma gives a sequential formulation of the vector field limit in~\cite[Lemma~2]{huang2022riemannian}. 

\begin{lemma}\label{RPG:le10}
Let $R$ be a smooth retraction and $V$ a continuous vector field on $\m$. Suppose that $\lim_{j\to\infty}x_j=x$ and $\eta_{x_j}\in\mathrm{T}_{x_j}\m$ satisfy $\lim_{j\to\infty}\|\eta_{x_j}\|=0$. Set $y_j=R_{x_j}(\eta_{x_j})$. Then $\lim_{j\to\infty}y_j=x$, and $\mathcal T_{\eta_{x_j}}^{-\sharp}$ exists and is uniformly bounded in operator norm for all sufficiently large $j$. Moreover,
\(
    \lim_{j\to\infty}\|V_{y_j}-\mathcal T_{\eta_{x_j}}^{-\sharp}V_{x_j}\|=0.
\)
\end{lemma}

\begin{proof}
By continuity of $R$, $\lim_{j\to\infty}y_j=R_x(0_x)=x$. Consider a smooth local orthonormal frame near $x$. Since $R$ is smooth and $\mathrm{D}R_x(0_x)=\mathrm{id}$, the matrices representing $\mathcal T_{\eta_{x_j}}^\sharp$ converge to the identity. Hence these operators are invertible for all sufficiently large $j$, and the matrices of their inverses also converge to the identity. In particular, the inverses are uniformly bounded in operator norm. By continuity of $V$, the components of $V_{x_j}$ and $V_{y_j}$ in the same frame both converge to those of $V_x$. The asserted limit follows.
\end{proof}

We next establish the stationarity of accumulation points of the safeguard sequence $\{\tilde{x}_i\}$ and a complexity bound for the safeguard directions.  Neither result relies on the specific choices of the acceleration parameters $A_k$, $\beta_k$, $\gamma_k$, and $\tau_k$.

\begin{theorem}
\label{thm:convergence}
Suppose that Assumption~\ref{assum:f-l-smooth,h-lipschitz} holds with $\Omega=\mathcal M$. Then the following statements hold for the safeguard sequence $\{\tilde{x}_i\}$ generated by \Cref{algrm:safeguard-rapg2}.
\begin{enumerate}[(i)]
    \item \label{thm: safe1}
    If the initial sublevel set $\Omega_{x_0}=\{x\in\m: F(x)\leqslant F(x_0)\}$ is compact, and $\sup_{i\geqslant0}L_i<\infty$, then $\{\tilde{x}_i\}$ has at least one accumulation point, and every accumulation point is stationary for Problem~\eqref{eq:f+h}.
    \item \label{thm: safe2}
    Suppose that $F_{\inf}=\inf_{x\in\m}F(x)>-\infty$. For any $\epsilon>0$, a safeguard point $\tilde{x}_i$ satisfying $\|\eta_{\tilde{x}_i}\|\leqslant\epsilon$ is identified within
    \(
        M_\epsilon=\max\left\{1,\left\lceil\frac{F(\tilde{x}_0)-F_{\inf}}{\sigma\iota^{N_{\mathrm{ls}}}\epsilon^2}\right\rceil\right\}
    \)
    calls to \Cref{alg:Safeguard}. The corresponding outer iteration index satisfies $k\leqslant N_{\max}M_\epsilon$.
\end{enumerate}
\end{theorem}

\begin{proof}
\eqref{thm: safe1}. By \Cref{le2}\eqref{le2:1} and $\tilde{x}_0=x_0$, all safeguard points lie in the compact set $\Omega_{x_0}$. Hence $\{\tilde{x}_i\}$ has at least one accumulation point. Since $F$ is continuous, it attains a minimum on $\Omega_{x_0}$, which is also its global minimum. Thus \Cref{le2}\eqref{le2:3} gives
\(
    \lim_{i\to\infty}\|\eta_{\tilde{x}_i}\|=0.
\)

Let $\tilde{x}_*$ be any accumulation point, and choose a subsequence satisfying $\lim_{j\to\infty}\tilde{x}_{i_j}=\tilde{x}_*$. By \Cref{RPG:le10},
\[
    \lim_{j\to\infty}\e_{\tilde{x}_{i_j}}(\eta_{\tilde{x}_{i_j}})=\tilde{x}_*,
\]
and $\mathcal T_{\eta_{\tilde{x}_{i_j}}}^{-\sharp}$ exists and is uniformly bounded in operator norm for all sufficiently large $j$. Applying the inverse adjoint to~\eqref{e22} and using the smooth sum rule, we obtain
\begin{equation*}
    \begin{aligned}
        r_j :=\grad f\bigl(\e_{\tilde{x}_{i_j}}(\eta_{\tilde{x}_{i_j}})\bigr)-\mathcal T_{\eta_{\tilde{x}_{i_j}}}^{-\sharp}\bigl(\grad f(\tilde{x}_{i_j})+\theta L_{i_j}\eta_{\tilde{x}_{i_j}}\bigr)
        \in\partial F\bigl(\e_{\tilde{x}_{i_j}}(\eta_{\tilde{x}_{i_j}})\bigr).
    \end{aligned}
\end{equation*}
Consequently,
\begin{equation}\label{eq:partialF_limit}
    \begin{aligned}
        \|r_j\| \leqslant \bigl\|\grad f\bigl(\e_{\tilde{x}_{i_j}}(\eta_{\tilde{x}_{i_j}})\bigr)-\mathcal T_{\eta_{\tilde{x}_{i_j}}}^{-\sharp}\grad f(\tilde{x}_{i_j})\bigr\|
        +\theta L_{i_j}\|\mathcal T_{\eta_{\tilde{x}_{i_j}}}^{-\sharp}\| \cdot \|\eta_{\tilde{x}_{i_j}}\|.
    \end{aligned}
\end{equation}
The first term on the right of~\eqref{eq:partialF_limit} tends to zero by \Cref{RPG:le10} with $V=\grad f$. The second term tends to zero because $L_{i_j}$ and $\|\mathcal T_{\eta_{\tilde{x}_{i_j}}}^{-\sharp}\|$ are uniformly bounded for all sufficiently large $j$ and $\lim_{j\to\infty}\|\eta_{\tilde{x}_{i_j}}\|=0$. Hence
$\bigl(\e_{\tilde{x}_{i_j}}(\eta_{\tilde{x}_{i_j}}),\ r_j\bigr)\to(\tilde{x}_*,0_{\tilde{x}_*})$ in the tangent bundle.
The closed graph property of the Clarke subdifferential~\cite[Theorem~2.2(c)]{hosseini2018line} then yields $0_{\tilde{x}_*}\in\partial F(\tilde{x}_*)$. Thus every accumulation point of $\{\tilde{x}_i\}$ is stationary.

\eqref{thm: safe2}. Suppose that $\|\eta_{\tilde{x}_i}\|>\epsilon$ for $i=0,\ldots,M_\epsilon-1$. Summing the descent inequalities in \Cref{le2}\eqref{le2:1} yields
\[
    \begin{aligned}
        F(\tilde{x}_0)-F_{\inf}
        & \geqslant F(\tilde{x}_0)-F(\tilde{x}_{M_\epsilon}) 
         \geqslant \sigma\iota^{N_{\mathrm{ls}}}\sum_{i=0}^{M_\epsilon-1}\|\eta_{\tilde{x}_i}\|^2
         >\sigma\iota^{N_{\mathrm{ls}}}M_\epsilon\epsilon^2.
    \end{aligned}
\]
This contradicts the definition of $M_\epsilon$, which gives $\sigma\iota^{N_{\mathrm{ls}}}M_\epsilon\epsilon^2\geqslant F(\tilde{x}_0)-F_{\inf}$.
 Thus the first safeguard call satisfying $\|\eta_{\tilde{x}_i}\|\leqslant\epsilon$ has index $i\leqslant M_\epsilon-1$. Since $N_j\leqslant N_{\max}$ for all $j$, the corresponding outer iteration index satisfies
\(
    k=\sum_{j=0}^{i}N_j\leqslant(i+1)N_{\max}\leqslant N_{\max}M_\epsilon.
\)
\end{proof}

Under the assumptions of \Cref{thm:convergence}\eqref{thm: safe1}, finitely many safeguard replacements imply finitely many restarts. The iterates from the last restart onward admit the rate in \Cref{thm:rate-rapg} under the local assumptions below.

\begin{corollary}\label{coro:GloCovRate}
Suppose that the assumptions of \Cref{thm:convergence}\eqref{thm: safe1} hold, and let $x_*$ be an accumulation point of the safeguard sequence $\{\tilde{x}_i\}$ and a minimizer of Problem~\eqref{eq:f+h}. Assume that the safeguard replaces the current iterate only finitely many times. Let $K$ be the outer iteration index of the final restart, with $K=0$ if no restart occurs. Suppose further that
\begin{enumerate}[(a)]
    \item there exists a neighborhood $\mathcal N_{x_*}$ of $x_*$ such that Assumption~\ref{assum:f,h-convex} and the geometric assumptions of \Cref{thm:rate-rapg} hold with $\Omega=\mathcal N_{x_*}$;
    \item the iterates generated by \Cref{algrm:safeguard-rapg2} satisfy $x_k,y_k\in\mathcal N_{x_*}$ for all $k\geqslant K$;
    \item from iteration $K$ onward, the acceleration parameters follow the recurrence in \Cref{algrm:rapg-sc}, with fixed parameters satisfying one of the conditions in \Cref{prop:cond-sc}.
\end{enumerate}
Then there exist constants $C_1,C_2>0$,  such that
\[
    F(x_k)-F(x_*)\leqslant\min\left\{\frac{C_1}{(k-K)^2}, \ C_2\left(1-\sqrt{\frac{q}{\xi}}\right)^{k-K}\right\},\qquad k>K,
\]
where $q$ is defined as in \Cref{thm:rate-rapg} with the eventual value of $L$.
\end{corollary}

\begin{proof}
By \Cref{thm:convergence}\eqref{thm: safe1}, $\sup_{i\geqslant0}L_i<\infty$. Since $L$ is nondecreasing and each increase multiplies it by $\tau_L>1$, only finitely many increases can occur. 
Together with the assumed finite number of safeguard replacements, this implies that the acceleration is restarted only finitely many times. Hence $K$ is well defined.

At the final restart, the safeguard sets $s_K=0_{x_K}$ and $A_K=A_0$. These conditions also hold at $K=0$ if no restart occurs. Since every increase of $L$ triggers a restart, $L$ remains constant from iteration $K$ onward. Subsequent safeguard calls neither replace the current iterate nor modify the acceleration state. Thus, by assumption~(c), the accelerated updates coincide with those of \Cref{algrm:rapg-sc} initialized at $x_K$.

Assumptions~(a)--(c) ensure that the resulting iterates satisfy the conditions of \Cref{thm:rate-rapg}. Applying \Cref{algrm:rapg-sc} with iteration count $k-K$ yields the stated bound.
\end{proof}

\begin{remark}[Smooth specialization]\label{remark:ar-rapg-smooth}
When $h\equiv0$, the proximal subproblem at $x$ has the unique solution $\eta_x=-\grad f(x)/(\theta L)$, where $L$ is the current smoothness estimate. Thus \Cref{algrm:safeguard-rapg2} reduces to a Riemannian accelerated gradient method with adaptive restart. Under the assumptions of \Cref{thm:convergence}\eqref{thm: safe1}, every accumulation point $\tilde{x}_*$ of the safeguard sequence satisfies $\grad f(\tilde{x}_*)=0$.
At each safeguard call, $\|\grad f(\tilde{x}_i)\|=\theta L_i\|\eta_{\tilde{x}_i}\|$. Hence, if $f$ is bounded below and $\sup_{i\geqslant0}L_i<\infty$, \Cref{thm:convergence}\eqref{thm: safe2} guarantees that a safeguard point satisfying $\|\grad f(\tilde{x}_i)\|\leqslant\epsilon$ is identified within $M_{\epsilon/(\theta\sup_{i\geqslant0}L_i)}$ safeguard calls, yielding an $\mathcal{O}(\epsilon^{-2})$ complexity bound. Under the assumptions of \Cref{coro:GloCovRate}, the rate bound specializes to $\rho=0$ and $q=\mu/(\theta L)$, where $L$ is the eventual value of the smoothness estimate.
\end{remark}

\section{Numerical experiments}
\label{sec:experiment}

We evaluate the proposed methods on four problem classes: two nonsmooth composite problems and two smooth problems. The nonsmooth tests are a regularized Karcher mean problem on hyperbolic space, whose objective is geodesically convex, and nonconvex sparse principal component analysis (SPCA) problem on the oblique manifold. The smooth tests are the Karcher mean problem on the manifold of symmetric positive definite (SPD) matrices, whose objective is globally geodesically convex, and the Rayleigh quotient problem on the sphere, which we test on a region where the objective is geodesically convex.
\Cref{tab:experiment-algorithms} lists the methods compared.

\begin{table}[H]
    \centering
    \caption{Methods compared in the numerical experiments.}
    \label{tab:experiment-algorithms}
    \ExperimentTableSetup
    \begin{tabular}{@{}l p{0.54\linewidth} p{0.20\linewidth}@{}}
        \toprule
        Method & Description & Source \\
        \midrule
        RGD       & Riemannian gradient descent & Standard baseline   \\
        RNAG-C    & Acceleration for geodesically convex objectives
        & \cite[Algorithm~1]{kim2022accelerated} \\
        RNAG-SC   & Acceleration for geodesically strongly convex objectives
        & \cite[Algorithm~2]{kim2022accelerated} \\
        RPG       & Riemannian proximal gradient method
        & \cite[Algorithm~1]{huang2022riemannian} \\
        ManPG     & Manifold proximal gradient method
        & \cite[Algorithm~1]{chen2020proximal} \\
        ManPG-Ada & Adaptive ManPG
        & \cite[Algorithm~2]{chen2020proximal} \\
        RAPG      & Riemannian accelerated proximal gradient method
        & \Cref{algrm:rapg-sc} \\
        AR-RAPG   & RAPG with adaptive restart
        & \Cref{algrm:safeguard-rapg2} \\
        \bottomrule
    \end{tabular}
\end{table}

All methods are implemented in MATLAB R2025a and run on a computer with a 2.3 GHz dual-core Intel Core i5 processor. Each instance uses ten initial points shared by all methods. Iteration counts and solver times are arithmetic means over all ten runs, including runs that reach the iteration limit. A run is counted as successful if the prescribed stopping criterion is satisfied. 

AR-RAPG computes $(A_{k+1},\beta_k,\gamma_k,\tau_k)$ by Step~\ref{alg:st3} of \Cref{algrm:rapg-sc}, with $(\xi,\theta,A_0)=(1,1,0)$. Unless specified otherwise, its safeguard parameters are $(\sigma,\iota,N_{\mathrm{ls}},\tau_L)=(10^{-4},0.5,3,1.1)$ and $(N_0,N_{\min},N_{\max})=(5,2,10)$. 
The objective plots show the relative error
\begin{equation}
    e_k=\frac{|F(x_k)-F_{\mathrm{ref}}|}{|F(x_0)-F_{\mathrm{ref}}|},
    \label{eq:common-relative-error}
\end{equation}
where the reference value $F_{\mathrm{ref}}$ denotes the exact optimal objective value when it is available, and a high-accuracy approximation otherwise.

\subsection{Nonsmooth composite optimization}
\label{subsec:nonsmooth-exp}

The regularized Karcher mean experiments on hyperbolic space compare practical and theoretically admissible acceleration settings under global geodesic strong convexity. The nonconvex SPCA experiments on the oblique manifold assess the AR-RAPG safeguard and compare solver performance.

\subsubsection{Regularized Karcher mean on hyperbolic space}
\label{subsec:hyp}

\paragraph{Problem statement.}
We represent hyperbolic space by the hyperboloid $\mathbb H^n=\{x\in\mathbb R^{n+1}:\langle x,x\rangle_L=-1,\ x_0>0\}$, where $\langle u,v\rangle_L=-u_0v_0+\sum_{j=1}^n u_jv_j$. The Riemannian metric at $x$ is the restriction of the Lorentzian inner product to $\T_x\mathbb H^n=\{u\in\mathbb R^{n+1}:\langle x,u\rangle_L=0\}$. For $n\geqslant2$, this metric has constant sectional curvature $-1$~\cite[Section~7.6]{boumal2023introduction}.
Given $a_1,\ldots,a_N,b\in\mathbb H^n$ and $r,\tau\geqslant0$, we solve
\begin{equation}
    \min_{x\in\mathbb H^n} F(x)
    =\underbrace{\frac{1}{2N}\sum_{i=1}^N\mathrm d^2(x,a_i)}_{f(x)}
    +\underbrace{\tau\operatorname{dist}(x,C)}_{h(x)},
    \qquad C=\overline B(b,r),
    \label{eq:hyp-model}
\end{equation}
where $C$ is the closed geodesic ball of radius $r$ centered at $b$. Each function $\tfrac12\mathrm d^2(\cdot,a_i)$ is geodesically $1$-strongly convex on $\mathbb H^n$. The geodesic convexity of $\mathrm d(\cdot,b)$ implies that $h(x)=\tau\max\{\mathrm d(x,b)-r,0\}$ is geodesically convex.  The exponential pullback $h\circ\e_x$ can be nonconvex, as shown in \Cref{ex:hyp-nonconvex-pullback}.

Let $o=(1,0,\ldots,0)^{\mathsf T}$, $D_\Omega>0$, and $\Omega=\overline B(o,D_\Omega/2)$. When all data points lie in $\Omega$, the Hessian comparison bounds in \Cref{subsec:zeta-and-delta} give the smoothness bound $L=\zeta_\Omega=D_\Omega\coth D_\Omega$ and strong convexity modulus $\mu=\delta_\Omega=1$ for $f$ on $\Omega$. By \Cref{prop:hyp-ball-pullback}, for each $x\in\Omega$, the exponential pullback of $h$ is $\rho(D_\Omega)$-weakly convex on $\{\eta\in\T_x\mathbb H^n:\|\eta\|_x\leqslant D_\Omega\}$, where
\(
    \rho(D_\Omega)=\tau\frac{\sinh D_\Omega\cosh D_\Omega-D_\Omega}{D_\Omega^2}.
\)

\paragraph{Experimental design.}
We compare RPG, RAPG, and AR-RAPG at $D_\Omega\in\{3,3.5,4,4.5\}$ with $\tau=0.01$.
For each diameter, we set the sampling radius to $S=D_\Omega/2-1$ and generate $N=100$ points $a_i=\operatorname{Exp}_o(r_iu_i)$, where $r_i\sim\Unif(0,S)$ and $u_i$ is uniform on the unit sphere of $\T_o\mathbb H^n$.
Since $S<D_\Omega/2$, all data points lie in $\Omega$.
The penalty ball has radius $r=0.8$ and center $b$ with $\mathrm d(o,b)=1.01$.
Each run starts from a point $x_0$ with $\mathrm d(o,x_0)=1.4$ and terminates when the proximal direction satisfies $\|\eta_k\|\leqslant10^{-6}$.
AR-RAPG initializes its smoothness estimate at $L_{\mathrm{init}}=0.5L$.

For all four diameters, $\rho=\rho(D_\Omega)<\mu=1$.
We test the two valid strong convexity moduli $\mu=1$ and $\mu=\rho$ for $f$, which give $q>0$ and $q=0$, respectively, where $q=(\mu-\rho)/(\theta L-\rho)$.
For RAPG, we consider two parameter regimes, termed \emph{practical} and \emph{theory-admissible}.
For each choice of $\mu$, both regimes share the same constants $(\mu,L,\rho)$ and differ only in $(\xi,\theta,A_0)$.
The practical regime, also used by AR-RAPG, takes $(\xi,\theta,A_0)=(1,1,0)$.
For the theory-admissible setting, we set $\nu=1/5$, $\lambda=2$, and $\xi=2\zeta_\Omega-1$ and take $\theta$ and $A_0$ equal to their respective lower bounds specified in \Cref{rem:q-branches,rem:xi-lambda}. 
\Cref{tab:hyp-parameters} lists the resulting parameters.

\begin{table}[!htbp]
 \centering
 \caption{Regional constants and theory-admissible RAPG parameters for \eqref{eq:hyp-model} with $\tau=0.01$ and $\rho=\rho(D_\Omega)$.}
 \label{tab:hyp-parameters}
 \begingroup\footnotesize
 \setlength{\tabcolsep}{4pt}\renewcommand{\arraystretch}{1.06}
 \begin{tabular}{@{}rrrr@{\hspace{12pt}}rrr@{\hspace{12pt}}rrr@{}}
 \toprule
 & & & &\multicolumn{3}{c}{$\mu=1$}&\multicolumn{3}{c}{$\mu=\rho$}\\
 \cmidrule(lr){5-7}\cmidrule(l){8-10}
 $D_\Omega$&$L=\zeta_\Omega$&$\rho$&$\xi$&$\theta$&$q$&$A_0$&$\theta$&$q$&$A_0$\\
 \midrule
    3.0 & 3.015 & 0.109 & 5.030 & 6.066 & 0.049 & 365.09 & 1 & 0 & 73.02 \\
    3.5 & 3.506 & 0.221 & 6.013 & 5.678 & 0.040 & 552.69 & 1 & 0 & 110.54 \\
    4.0 & 4.003 & 0.463 & 7.005 & 4.166 & 0.033 & 781.34 & 1 & 0 & 156.27 \\
    4.5 & 4.501 & 0.998 & 8.002 & 1.000 & 0.001 & 1050.64 & 1 & 0 & 210.13 \\
 \bottomrule
 \end{tabular}
 \endgroup
\end{table}

The three methods share the same solver for their tangent proximal subproblems (e.g., Step~\ref{alg:st5} of \Cref{algrm:rapg-sc}).
At each base point $y$, we restrict the subproblem to the tangent ball $\{\eta\in\T_y\mathbb H^n:\|\eta\|_y\leqslant D_T\}$ with $D_T=1.41$.
We verify that the quadratic regularization parameter $\theta L$, computed with the current value of $L$, exceeds $\rho(D_T)$; by \Cref{prop:hyp-ball-pullback}, this ensures that the subproblem objective is strongly convex on the ball.
Reflection across $V_y=\operatorname{span}\{\operatorname{grad}f(y),\operatorname{Exp}_y^{-1}(b)\}$ preserves both the subproblem objective and the feasible ball, so the unique minimizer lies in $V_y$.
The subproblem therefore reduces to at most two dimensions, which we solve by sequential quadratic programming using MATLAB's \texttt{fmincon}.
In all reported runs, every accepted numerical solution lies strictly inside the tangent ball.

We take $F_{\mathrm{ref}}$ to be the objective value at termination of the intrinsic Riemannian proximal gradient method~\cite{bergmann2026convex}, run with step size $1/L$ until $\mathrm d(x_{k+1},x_k)\leqslant10^{-10}$.
For $t>0$, the intrinsic proximal mapping of $th$ used by this method has the closed form
\[
    \operatorname{prox}_{th}(y)
    =\arg\min_{x\in\mathbb H^n}\left\{h(x)+\frac{\mathrm d^2(x,y)}{2t}\right\}
    =\begin{cases}
        y, & \mathrm d(y,b)\leqslant r,\\
        \operatorname{Exp}_y \left(
        \frac{\min\{t\tau,\ \mathrm d(y,b)-r\}}{\mathrm d(y,b)}
        \operatorname{Exp}_y^{-1}(b)\right), & \mathrm d(y,b)>r.
    \end{cases}
\]

\begin{figure}[!htpb]
    \centering
    \includegraphics[width=\textwidth]{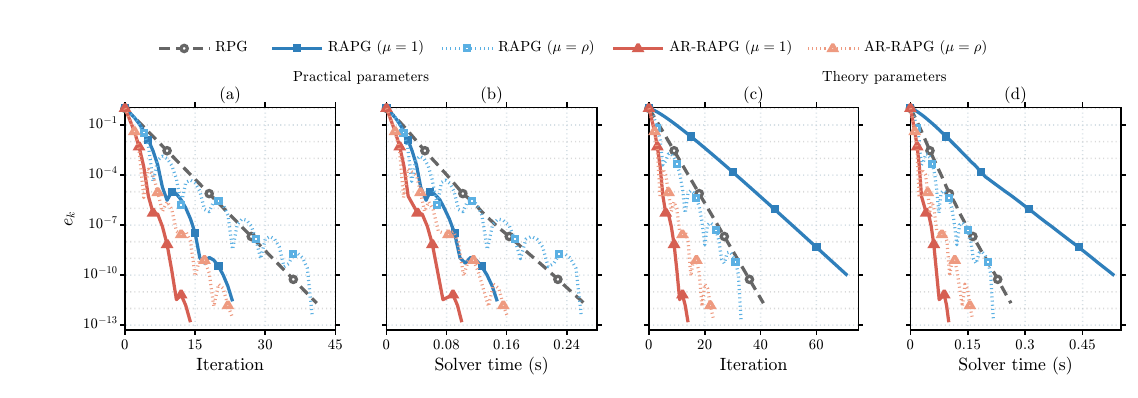}
    \caption{Convergence for the regularized hyperbolic Karcher mean at $(n,D_\Omega,\tau)=(100,4,0.01)$. RAPG uses practical parameters in (a)--(b) and theory-admissible parameters from \Cref{tab:hyp-parameters} in (c)--(d). RPG and AR-RAPG results are shared across the comparisons.}
    \label{fig:hyp-main}
\end{figure}

\paragraph{Comparison results.}
On the main instance in \Cref{fig:hyp-main}, AR-RAPG with $\mu=1$ has the lowest mean iteration count and solve time. Practical RAPG performs better with $\mu=1$ than with $\mu=\rho$, whereas this ordering reverses with theory-admissible parameters. In the latter comparison, $\mu=1$ changes $(\theta,q,A_0)$ simultaneously and reduces the step size through a larger $\theta$. 

In the regional comparison, \Cref{tab:hyp-regions} shows that AR-RAPG with $\mu=1$ attains the lowest mean iteration count and solver time at every diameter, tying in iteration count with $\mu=\rho$ at $D_\Omega=4.5$.
The iteration count of RPG increases with $D_\Omega$, whereas that of RAPG, which uses the conservative theory-admissible parameters in \Cref{tab:hyp-parameters}, varies nonmonotonically.
RAPG with $\mu=1$ requires more iterations than RPG at every tested diameter, which is likely due to these conservative parameters.
For $D_\Omega\leqslant4$, the choice $\theta>1$ shrinks the proximal step $1/(\theta L)$ below the step $1/L$ used by RPG, and the large initial value $A_0$ may further delay the effect of momentum.
At $D_\Omega=4.5$, where $\theta=1$, the proximal step coincides with that of RPG, so the remaining gap more likely attributable to $(q,A_0)$.

The RAPG and AR-RAPG iterates $x_k,y_k$ remained in the chosen region in all reported runs.

\begin{table}[!htbp]
  \centering
  \caption{Hyperbolic regional comparison with $\tau=0.01$. RAPG uses the theory-admissible parameters in \Cref{tab:hyp-parameters}.}
  \label{tab:hyp-regions}
  \footnotesize
  \setlength{\tabcolsep}{1.8pt}
  \renewcommand{\arraystretch}{1.08}
  \begin{tabular}{@{}r *{4}{>{\columncolor[gray]{0.96}}r r@{\hspace{8pt}}}
                      >{\columncolor[gray]{0.96}}r r@{}}
    \toprule
    & \multicolumn{2}{c}{RPG}
    & \multicolumn{4}{c}{RAPG}
    & \multicolumn{4}{c}{AR-RAPG} \\
    \cmidrule(lr){2-3}\cmidrule(lr){4-7}\cmidrule(lr){8-11}
    & & & \multicolumn{2}{c}{$\mu=1$} & \multicolumn{2}{c}{$\mu=\rho$}
        & \multicolumn{2}{c}{$\mu=1$} & \multicolumn{2}{c}{$\mu=\rho$} \\
    $D_\Omega$ & Iterations & Time (s) & Iterations & Time (s)
               & Iterations & Time (s) & Iterations & Time (s) & Iterations & Time (s) \\
    \midrule
    3.0 & 30.0 & 0.203 & 87.0 & 0.654 & 32.0 & 0.208
        & \textbf{10.0} & \textbf{0.065} & 17.0 & 0.120 \\
    3.5 & 35.0 & 0.219 & 91.3 & 0.679 & 36.0 & 0.225
        & \textbf{11.0} & \textbf{0.073} & 21.0 & 0.149 \\
    4.0 & 41.0 & 0.263 & 71.1 & 0.533 & 33.0 & 0.217
        & \textbf{14.0} & \textbf{0.100} & 23.0 & 0.162 \\
    4.5 & 46.0 & 0.315 & 51.0 & 0.341 & 36.0 & 0.237
        & \textbf{28.0} & \textbf{0.201} & \textbf{28.0} & 0.211 \\
    \bottomrule
  \end{tabular}
\end{table}

\subsubsection{Sparse principal component analysis on the oblique manifold}
\label{subsec:spca}

\paragraph{Problem statement.}
Following~\cite{genicot2015weakly}, we consider the weakly correlated SPCA problem
\begin{equation}
    \min_{X\in\operatorname{OB}(p,n)} F(X)
    =\underbrace{\bigl\|X^{\mathsf T}A^{\mathsf T}AX-D_p^2\bigr\|_{\mathrm F}^2}_{f(X)}
    +\underbrace{\lambda\|X\|_1}_{h(X)},
    \label{eq:spca}
\end{equation}
where $\operatorname{OB}(p,n)=\{X\in\mathbb R^{n\times p}:\operatorname{diag}(X^{\mathsf T}X)=\mathbf 1_p\}$. The rows of $A\in\mathbb R^{m\times n}$ are samples, and its columns are features. The diagonal matrix $D_p=\operatorname{diag}(\sigma_1(A),\ldots,\sigma_p(A))$ contains the leading $p$ singular values of $A$ in nonincreasing order. With the Euclidean product metric, $\operatorname{OB}(p,n)$ is a product of spheres.  Its sectional curvature lies in $[0,1]$ by \Cref{lem:oblq}.

We use random Gaussian matrices with independent $\mathcal N(0,1)$ entries and two real data sets, CPTAC breast cancer (BRCA) protein abundances\footnote{\href{https://github.com/SaccoPerfettoLab/PatientProfiler_BRCA}{PatientProfiler BRCA repository}, \texttt{Step1/Brca\_prot.tsv}.} and PIRvision sensor measurements\footnote{\href{https://doi.org/10.24432/C56W5M}{PIRvision dataset, UCI Machine Learning Repository}.}. For BRCA, we retain all 122 samples and 7767 completely observed, nonconstant protein features. For PIRvision, we use the first office dataset, retaining all 7651 observations and 56 numerical features after excluding the date, time, and activity label. We center each feature column and normalize it to unit Euclidean norm before computing $D_p$.

\paragraph{Experimental design.}
We compare RPG, ManPG, ManPG-Ada, RAPG, and AR-RAPG, using the exponential map as the retraction in all five methods. RAPG and AR-RAPG use the empirical parameter settings $(\xi,\theta,A_0)=(1,1,0)$ and $(\mu,\rho)=(1,0.5)$. For ManPG and ManPG-Ada, we retain the parameter settings and line search conditions of the original Stiefel manifold implementations~\cite{chen2020proximal}.

Since $\operatorname{OB}(p,n)$ is a product of spheres, the Riemannian proximal subproblems arising in RPG, RAPG, and AR-RAPG (e.g., Step~\ref{alg:st5} of \Cref{algrm:rapg-sc}) for~\eqref{eq:spca} decompose into independent proximal subproblems on the sphere. These subproblems are solved by the fixed-point iteration developed in~\cite{huang2022riemannian}. For ManPG and ManPG-Ada, the corresponding subproblems are solved by the semismooth Newton method used in~\cite{chen2020proximal}.

For each run, we form the initial point by perturbing each of the leading $p$ right singular vectors of $A$ by $0.05$ radians along a random tangent direction. Both experiments use the stopping criterion
\(
    \mathcal R_k:=\frac{L_k\|\eta_k\|_{\mathrm F}}{\sqrt{np}}\leqslant 2\times10^{-4},
\)
where $\eta_k$ is the proximal direction at iteration $k$ and $L_k$ is the smoothness estimate used to compute it. 
Sparsity is the percentage of entries with $|X_{ij}|\leqslant10^{-5}$.
For each instance, $F_{\mathrm{ref}}$ denotes the smallest objective value recorded across all methods and iterates.

\paragraph{Effect of the safeguard in AR-RAPG.}
To compare RAPG and AR-RAPG with matching initial steps, we fix two Gaussian instances with $(m,p,\lambda)=(20,4,1)$ and $n\in\{500,5000\}$. RAPG uses $L=rL_{\mathrm{base}}$, and AR-RAPG starts with $L_{\mathrm{init}}=rL_{\mathrm{base}}$, where $L_{\mathrm{base}}=4\|D_p^2\|_{\mathrm F}^2$ and $r\in\{1,0.5,0.25\}$.   
We retain the residual tolerance above and use an iteration limit of $10000$.

\begin{figure}[H]
	\centering
	\includegraphics[width=\linewidth]{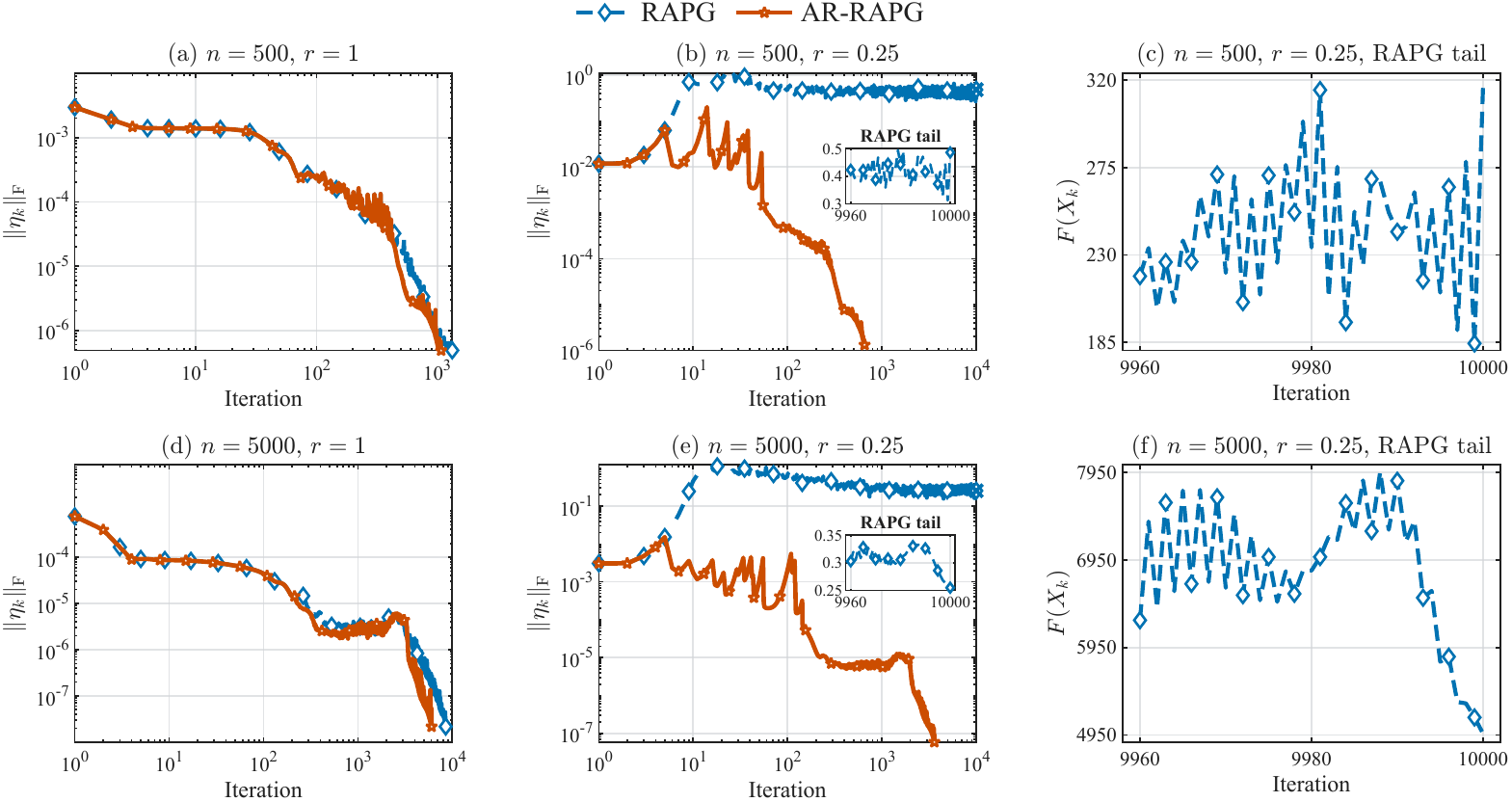}
	\caption{Safeguard comparison for Gaussian SPCA. Crosses mark curves that include runs reaching the iteration limit.}
	\label{fig:spca-safeguard}
\end{figure}

For $r=1$ and $r=0.5$, both methods satisfy the tolerance in all runs, and AR-RAPG requires fewer iterations without increasing $L$, with safeguard replacements and momentum resets occurring along the way.
At $r=0.25$, RAPG reaches the iteration limit in every run, whereas AR-RAPG succeeds in all runs, increasing $L$ four times on average; see \Cref{tab:spca-safeguard}.
The mean trajectories in \Cref{fig:spca-safeguard} show an overall decrease in the proximal direction norms for AR-RAPG, while RAPG exhibits persistent oscillations in both direction norms and objective values for $r=0.25$.

\begin{table}[!htbp]
  \centering
  \caption{Safeguard comparison for Gaussian SPCA. The subscript $\mathrm{end}$ denotes values at termination. Success counts runs satisfying $\mathcal R_{\mathrm{end}}\leqslant2\times10^{-4}$. Other entries are arithmetic means including runs reaching the 10,000-iteration limit. Accept. counts safeguard replacements, each accompanied by a momentum reset; $L$ upd. counts all increases in $L$.}
  \label{tab:spca-safeguard}
  \begingroup
  \footnotesize
  \setlength{\tabcolsep}{1.55pt}
  \medmuskip=1.25mu\relax
  \def\itergap{\hspace{9pt}}
  \renewcommand{\arraystretch}{1.06}
  \begin{tabular}{@{}lr>{\columncolor[gray]{0.955}[1.6pt][2pt]}r@{\itergap}rr@{\hspace{7pt}}r>{\columncolor[gray]{0.955}[1.6pt][2pt]}r@{\itergap}rrrr@{}}
    \toprule
    & \multicolumn{4}{c}{RAPG} & \multicolumn{6}{c}{AR-RAPG} \\
    \cmidrule(r{3pt}){2-5}\cmidrule(l{3pt}){6-11}
    $r$ & Success & Iterations & $\|\eta_{\mathrm{end}}\|_{\mathrm F}$ & $\mathcal{R}_{\mathrm{end}}$
        & Success & Iterations & $\|\eta_{\mathrm{end}}\|_{\mathrm F}$ & $\mathcal{R}_{\mathrm{end}}$ & Accept. & $L$ upd. \\
    \midrule
    \multicolumn{11}{@{}l}{$A\in\mathbb{R}^{20\times500}$} \\
    1    & 10/10 &  1327.9 & $4.88\times10^{-7}$ & $2.00\times10^{-4}$ & 10/10 & 1074.3 & $4.87\times10^{-7}$ & $1.99\times10^{-4}$ & 1.0 & 0.0 \\
    0.5  & 10/10 &   922.5 & $9.76\times10^{-7}$ & $2.00\times10^{-4}$ & 10/10 &  758.0 & $9.74\times10^{-7}$ & $1.99\times10^{-4}$ & 1.0 & 0.0 \\
    0.25 &  0/10 & 10000.0 & $4.86\times10^{-1}$ & $4.97\times10^{1}$  & 10/10 &  655.5 & $1.33\times10^{-6}$ & $1.99\times10^{-4}$ & 5.0 & 4.0 \\
    \addlinespace[3pt]
    \multicolumn{11}{@{}l}{$A\in\mathbb{R}^{20\times5000}$} \\
    1    & 10/10 &  8526.8 & $2.16\times10^{-8}$ & $2.00\times10^{-4}$ & 10/10 & 6035.9 & $2.16\times10^{-8}$ & $2.00\times10^{-4}$ & 2.0 & 0.0 \\
    0.5  & 10/10 &  5996.6 & $4.32\times10^{-8}$ & $2.00\times10^{-4}$ & 10/10 & 4247.4 & $4.32\times10^{-8}$ & $2.00\times10^{-4}$ & 2.0 & 0.0 \\
    0.25 &  0/10 & 10000.0 & $2.55\times10^{-1}$ & $5.89\times10^{2}$  & 10/10 & 3594.9 & $5.91\times10^{-8}$ & $2.00\times10^{-4}$ & 6.1 & 4.0 \\
    \bottomrule
  \end{tabular}
  \endgroup
\end{table}

\paragraph{Comparison results.}
We next compare all five methods on the Gaussian, BRCA, and PIRvision data sets.
The fixed-step methods use $L=L_{\mathrm{base}}$, and AR-RAPG starts from $L_{\mathrm{init}}=0.5L_{\mathrm{base}}$.
Each run is capped at $30000$ iterations.

\begin{figure}[!htbp]
	\centering
	\includegraphics[width=\linewidth]{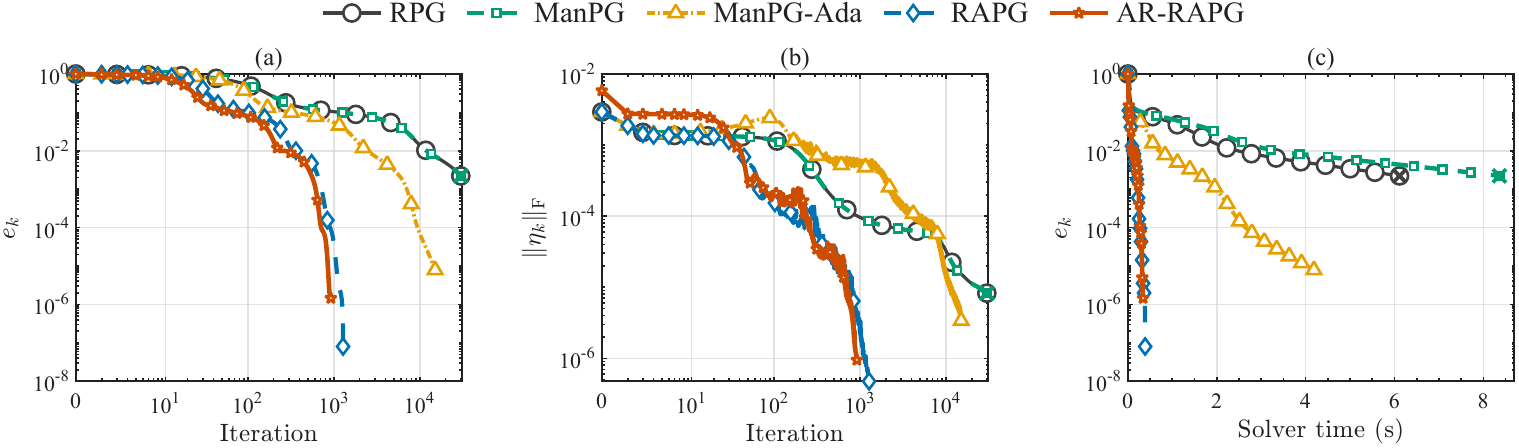}
	\caption{Convergence on the Gaussian SPCA instance with $(m,n,p,\lambda)=(20,500,4,1)$. Crosses mark methods for which at least one run reached the iteration limit.}
	\label{fig:spca-main}
\end{figure}

\Cref{tab:spca-main} summarizes the results on the three data sets, and \Cref{fig:spca-main} shows the convergence on the Gaussian instance.
ManPG-Ada, RAPG, and AR-RAPG meet the residual tolerance in all runs.
Both accelerated methods require substantially fewer iterations and less solver time than ManPG-Ada while attaining comparable final objective values.
AR-RAPG has the lowest mean iteration count on all three data sets and the lowest mean solver time on the Gaussian and BRCA data sets, whereas RAPG has a slightly lower mean solver time on PIRvision.
In contrast, RPG and ManPG reach the iteration limit in every Gaussian and BRCA run, but meet the tolerance on PIRvision.

\begin{table}[!htbp]
  \centering
  \caption{SPCA results on Gaussian and real datasets. The subscript $\mathrm{end}$ denotes values at termination. Success counts runs satisfying $\mathcal R_{\mathrm{end}}\leqslant2\times10^{-4}$. The other entries are means over ten runs, including runs that reach the iteration limit. Parentheses give sample standard deviations for iteration counts and solver times.}
  \label{tab:spca-main}
  \begingroup
  \footnotesize
  \setlength{\tabcolsep}{2.7pt}
  \renewcommand{\arraystretch}{1.07}
  \begin{tabular}{@{}l>{\columncolor[gray]{0.96}}rrrrrcr@{}}
    \toprule
    Method & Iterations & Time (s) & $F_{\mathrm{end}}$
      & $\|\eta_{\mathrm{end}}\|_{\mathrm{F}}$ & $\mathcal{R}_{\mathrm{end}}$
      & Sparsity (\%) & Success \\
    \midrule
    \multicolumn{8}{@{}l}{\textbf{Random Gaussian:} $(m,n,p,\lambda)=(20,500,4,1)$} \\
    RPG       & 30,000.0 (0.0)    & 6.114 (0.750) & 66.8483 & $8.26\times10^{-6}$ & $3.45\times10^{-3}$ & 20.18 & 0/10 \\
    ManPG     & 30,000.0 (0.0)    & 8.362 (1.432) & 66.8483 & $8.25\times10^{-6}$ & $3.45\times10^{-3}$ & 20.18 & 0/10 \\
    ManPG-Ada & 15,097.3 (361.8)  & 4.187 (0.158) & 66.8336 & $3.39\times10^{-6}$ & $1.99\times10^{-4}$ & 20.40 & 10/10 \\
    RAPG      & 1,281.7 (50.7)    & 0.389 (0.050) & 66.8335 & $4.76\times10^{-7}$ & $1.99\times10^{-4}$ & 20.30 & 10/10 \\
    AR-RAPG   & \textbf{924.2 (69.7)} & \textbf{0.343 (0.071)} & 66.8335 & $9.51\times10^{-7}$ & $1.99\times10^{-4}$ & 20.35 & 10/10 \\
    \addlinespace
    \multicolumn{8}{@{}l}{\textbf{Real CPTAC BRCA:} $(m,n,p,\lambda)=(122,7767,4,10)$} \\
    RPG       & 30,000.0 (0.0)    & 78.430 (2.666)  & 2731.8094 & $6.27\times10^{-7}$ & $2.08\times10^{-2}$ & 10.17 & 0/10 \\
    ManPG     & 30,000.0 (0.0)    & 101.056 (1.569) & 2731.8094 & $6.27\times10^{-7}$ & $2.08\times10^{-2}$ & 10.17 & 0/10 \\
    ManPG-Ada & 23,479.9 (223.1)  & 92.309 (2.834)  & 2731.7984 & $2.67\times10^{-8}$ & $1.99\times10^{-4}$ & 10.17 & 10/10 \\
    RAPG      & 6,160.3 (0.5)     & 24.948 (1.540)  & 2731.7984 & $5.94\times10^{-9}$ & $1.97\times10^{-4}$ & 10.17 & 10/10 \\
    AR-RAPG   & \textbf{856.1 (23.3)} & \textbf{3.830 (0.176)} & 2731.7984 & $1.08\times10^{-8}$ & $1.96\times10^{-4}$ & 10.17 & 10/10 \\
    \addlinespace
    \multicolumn{8}{@{}l}{\textbf{Real PIRvision:} $(m,n,p,\lambda)=(7651,56,4,0.12)$} \\
    RPG       & 12,438.3 (128.8)  & 1.601 (0.084) & 2.7985 & $1.90\times10^{-6}$ & $2.00\times10^{-4}$ & 19.20 & 10/10 \\
    ManPG     & 12,442.3 (128.8)  & 1.736 (0.068) & 2.7985 & $1.90\times10^{-6}$ & $2.00\times10^{-4}$ & 19.20 & 10/10 \\
    ManPG-Ada & 3,662.9 (90.4)    & 0.536 (0.017) & 2.7985 & $9.42\times10^{-6}$ & $1.99\times10^{-4}$ & 19.20 & 10/10 \\
    RAPG      & 430.1 (2.8)       & \textbf{0.081 (0.006)} & 2.7985 & $1.89\times10^{-6}$ & $1.99\times10^{-4}$ & 19.20 & 10/10 \\
    AR-RAPG   & \textbf{352.7 (22.2)} & 0.088 (0.012) & 2.7985 & $3.42\times10^{-6}$ & $1.98\times10^{-4}$ & 19.20 & 10/10 \\
    \bottomrule
  \end{tabular}
  \endgroup
\end{table}

\subsection{Smooth optimization}
\label{subsec:smooth-exp}

We set $h\equiv0$ and $\rho=0$. The SPD Karcher mean objective is globally geodesically convex, whereas the spherical Rayleigh quotient is geodesically convex on the region specified below.

\paragraph{Common setup.}
We compare RGD, RNAG-C, RNAG-SC, RAPG, and AR-RAPG under practical and theory-admissible parameter regimes. Both regimes uses the same $(L,\mu)$.
The practical regime uses $(\xi,\theta,A_0)=(1,1,0)$ for RAPG, $(\xi,T)=(1,0)$ for RNAG-C, and $\xi=1$ for RNAG-SC.
For theory-admissible RAPG, we select $(\xi,\theta,A_0)$ according to \Cref{rem:q-branches,rem:xi-lambda} with $\nu=1/5$. RNAG-C uses $\xi=\zeta_\Omega+3(\zeta_\Omega-\delta_\Omega)$, $T=4\xi$, and step size $1/L$ as in \cite[Corollary~5.5]{kim2022accelerated}. RNAG-SC uses the same $\xi$, following \cite[Corollary~5.7]{kim2022accelerated}. For its step size, we numerically find the upper endpoint $\widehat q_{\max}$ of the admissible interval adjacent to zero, subject to $\widehat q\leqslant\mu/L$ and the conditions in \cite[Theorem~5.6]{kim2022accelerated}, where $\widehat q=\mu s$. We then set $s=0.99\widehat q_{\max}/\mu$.
RGD and AR-RAPG retain their settings across the two regimes. All runs stop when the gradient norm is at most $10^{-6}$.

\subsubsection{Karcher mean on the SPD manifold}
\label{subsubsec:spd-exp}

\paragraph{Problem statement.}

 Let $\Spp^d$ denote the manifold of $d\times d$ SPD matrices equipped with the affine-invariant metric.
Given $A_1,\ldots,A_N\in\Spp^d$, we consider the Riemannian center of mass (Karcher mean) problem
\begin{equation}
    \min_{X\in\Spp^d} F(X)=\frac{1}{2N}\sum_{i=1}^N\mathrm d^2(X,A_i).
    \label{eq:spd-karcher}
\end{equation}
The sectional curvature of $\Spp^d$ lies in $[-1/2,0]$~\cite{criscitiello2023accelerated}.
On a geodesically convex region $\Omega$ of diameter $D_\Omega$ containing the data, the Hessian comparison bounds in \Cref{subsec:zeta-and-delta} yield $\delta_\Omega=1$ and $\zeta_\Omega=(D_\Omega/\sqrt2)\coth(D_\Omega/\sqrt2)$.
Hence $F$ is geodesically $\mu$-strongly convex and geodesically $L$-smooth on $\Omega$, with $\mu=\delta_\Omega=1$ and $L=\zeta_\Omega$.

 \paragraph{Experimental design.}
For each $(d,D_\Omega)$, we set $\Omega=B(I_d,D_\Omega/2)$, where $I_d$ is the identity matrix, and generate $N=10$ matrices $A_i=\e_{I_d}(0.4D_\Omega u_iU_i)$.
Here $u_i\sim\Unif(0,1)$, and $U_i\in\T_{I_d}\Spp^d$ is a symmetric Gaussian matrix normalized to unit norm.
Since $\mathrm d(I_d,A_i)=0.4D_\Omega u_i<D_\Omega/2$, all data points lie in $\Omega$.
For each data set, the initial point is $X_0=\e_{I_d}(0.4D_\Omega\alpha U_0)$, where $\alpha\sim\Unif(0.2,0.6)$ and $U_0$ is sampled independently in the same way as $U_i$.

Both regimes use $L=\zeta_\Omega$ and $\mu=1$.
For theory-admissible RAPG, we set $\lambda=2$ and $\xi=2\zeta_\Omega-1$ as in \Cref{rem:xi-lambda}, and take $\theta$ and $A_0$ at the lower bounds in \Cref{rem:q-branches}.
\Cref{tab:spd-theory-params} reports the resulting values.
%
AR-RAPG starts from $L_{\mathrm{init}}=0.4L$.
To compute $F_{\mathrm{ref}}$, we start from the log-Euclidean mean of the data and apply the Manopt conjugate-gradient solver~\cite{boumal2014manopt}, followed by gradient refinement with a target tolerance of $2\times10^{-8}$.

\begin{table}[!htbp]
\centering
\caption{Problem constants and theory-admissible RAPG parameters for~\eqref{eq:spd-karcher} with $\mu=1$.}
\label{tab:spd-theory-params}
\ExperimentTableSetup
\begin{tabular}{@{}rrrrr@{}}
\toprule
$D_\Omega$ & $L=\zeta_\Omega$ & $\xi$ & $\theta$ & $A_0$ \\
\midrule
3.75 & 2.678  & 4.356  & 6.373 & 258.86   \\
7.5  & 5.304  & 9.607  & 8.139 & 1567.72  \\
15   & 10.607 & 20.213 & 9.063 & 7575.08  \\
30   & 21.213 & 41.426 & 9.530 & 33090.15 \\
\bottomrule
\end{tabular}
\end{table}

\paragraph{Comparison results.}

\begin{figure}[ht]
    \centering
    \includegraphics[width=\textwidth]{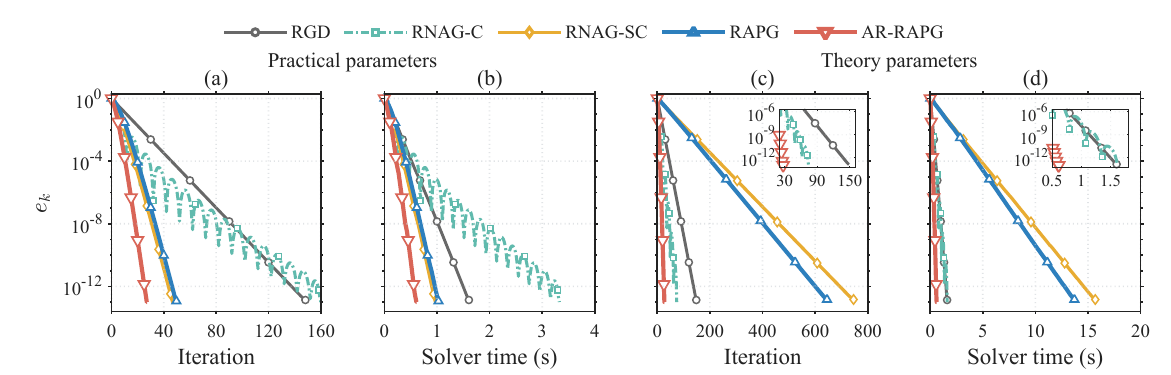}
    \caption{Convergence on the SPD Karcher mean  at $(d,D_\Omega)=(100,15)$. RNAG-C, RNAG-SC, and RAPG use practical parameters in panels (a)--(b) and theory-admissible parameters in panels (c)--(d). RGD and AR-RAPG results are shared across the two comparisons.}
    \label{fig:spd-main}
\end{figure}

With practical parameters, AR-RAPG has the lowest mean iteration count and solver time on the main instance; see \Cref{fig:spd-main}. RNAG-SC and RAPG perform similarly and improve on RGD, whereas RNAG-C oscillates and requires more solver time than RGD.
With theory-admissible parameters, RNAG-C requires fewer iterations than RGD with comparable solver time. RAPG requires fewer iterations and less time than RNAG-SC, but both require substantially more time than RGD. These comparisons involve simultaneous changes in acceleration parameters and step sizes: RAPG uses a larger $\theta$ and hence a smaller step $1/(\theta L)$, and RNAG-SC also uses a substantially smaller step.

\begin{figure}[!htbp]
\centering
\includegraphics[width=\textwidth]{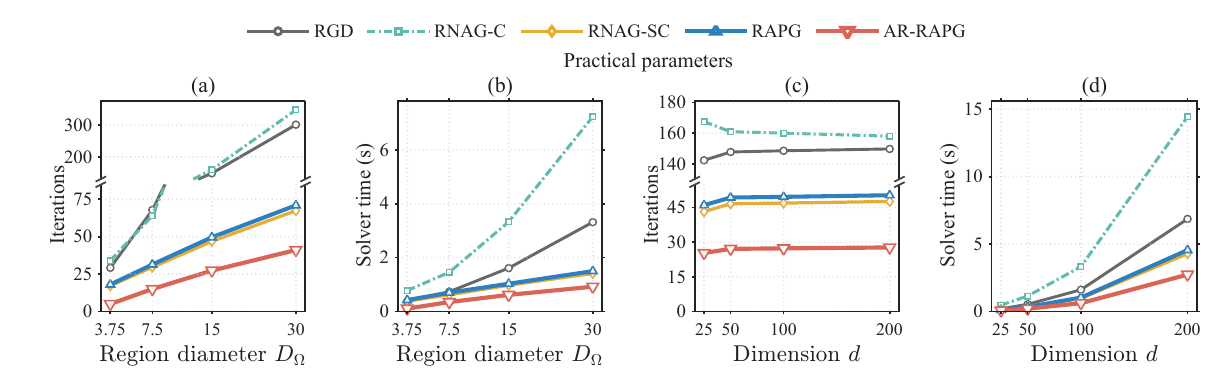}
\caption{SPD scaling with practical parameters: (a)--(b) vary $D_\Omega$ at $d=100$; (c)--(d) vary $d$ at $D_\Omega=15$.}
\label{fig:spd-scaling}
\end{figure}

AR-RAPG has the lowest mean iteration count and solver time in all practical scaling tests; see \Cref{fig:spd-scaling}. Increasing $D_\Omega$ changes both the data and $L$, with iteration counts increasing most sharply for RGD and RNAG-C. At fixed diameter, iteration counts vary little with $d$, while solver times increase as matrix operations become more costly.

The RAPG and AR-RAPG iterates $X_k,Y_k$ remained in $\Omega$ in all reported runs.

\subsubsection{Rayleigh quotient on the sphere}
\label{subsubsec:sphere-exp}

\paragraph{Problem statement.}
We consider the Rayleigh quotient minimization problem on the unit sphere
\begin{equation}
    \min_{x\in\mS^{d-1}} F(x) = -\tfrac12\, x^\top A x,
    \qquad \mS^{d-1}=\{x\in\mathbb R^d : \|x\|_2=1\},
    \label{eq:sphere-rayleigh}
\end{equation}
where $A\in\mathbb{R}^{d\times d}$ is symmetric with eigenvalues $\lambda_1>\lambda_2\geqslant\cdots\geqslant\lambda_d$. We fix $\lambda_1=1000$ and $\lambda_d=1$. The optimal value is $F^*=-\lambda_1/2$, and $F$ is geodesically $L$-smooth with $L=\lambda_1-\lambda_d=999$; see \cite[Proposition~7.1]{kim2022accelerated}.
Let $\Delta=\lambda_1-\lambda_2$ be the spectral gap and $x_*$ a unit eigenvector for $\lambda_1$. We work on the geodesic ball $\Omega=B(x_*,R_\Omega)$ with $R_\Omega=0.95\,R_{\mathrm{cert}}$, where $R_{\mathrm{cert}}=\arcsin\sqrt{\Delta/(L+\Delta)}$, and set $D_\Omega=2R_\Omega$. Then $D_\Omega<\pi/2$, so $\Omega$ is uniquely geodesically convex, and $F$ is geodesically $\mu$-strongly convex on $\Omega$ with $\mu=\Delta-(L+\Delta)\sin^2R_\Omega>0$; see \Cref{lem:sphere-rayleigh-local-convexity}. The other geometric constants are $\zeta_\Omega=1$ and $\delta_\Omega=D_\Omega\cot D_\Omega$; see \Cref{subsec:zeta-and-delta}.

\paragraph{Experimental design.}

For each $(d,\Delta)$, we set $\lambda_2=1000-\Delta$ and take $\lambda_3,\ldots,\lambda_{d-1}$ logarithmically spaced in $[1,10]$ in nonincreasing order. The eigenbasis is obtained by applying a Householder transformation to the standard basis. Each initial point is $x_0=\e_{x_*}(\alpha R_\Omega u)$, where $\alpha\sim\Unif(0.5,0.9)$ and $u\in\T_{x_*}\mS^{d-1}$ is an independently sampled unit tangent vector. In the spectral gap sweep, the Householder transformation, $\alpha$, and $u$ are reused across values of $\Delta$ so that instances and initializations are paired. We use the exact reference value $F_{\mathrm{ref}}=F^*$.
Both the practical and the theory-admissible settings use $L=999$ and $\mu=\Delta-(L+\Delta)\sin^2R_\Omega$. AR-RAPG is initialized with $L_{\mathrm{init}}=0.8L$ and uses $N_{\max}=20$. For theory-admissible RAPG, we set $\lambda=5/2$ and $\xi=(5-2\delta_\Omega)/3$ as in \Cref{rem:xi-lambda}, and take $\theta$ and $A_0$ at the lower bounds in \Cref{rem:q-branches}; the resulting values are listed in \Cref{tab:sphere-theory-params}.

\begin{table}[!htbp]
\centering
\caption{Problem constants and theory-admissible RAPG parameters for the spectral gap sweep of~\eqref{eq:sphere-rayleigh}  at $d=1000$.}
\label{tab:sphere-theory-params}
\ExperimentTableSetup
\begin{tabular}{@{}rrrrrrr@{}}
\toprule
$\Delta$ & $D_\Omega/\pi$ & $\delta_\Omega$ & $\mu$
& $\xi$ & $\theta$ & $A_0$ \\
\midrule
1  & 0.019 & 0.999 & 0.097 & 1.001 & 1 & 0.01 \\
5  & 0.043 & 0.994 & 0.487 & 1.004 & 1 & 0.05 \\
10 & 0.060 & 0.988 & 0.972 & 1.008 & 1 & 0.11 \\
20 & 0.085 & 0.976 & 1.938 & 1.016 & 1 & 0.22 \\
\bottomrule
\end{tabular}
\end{table}

\begin{figure}[!htbp]
    \centering
    \includegraphics[width=\textwidth]{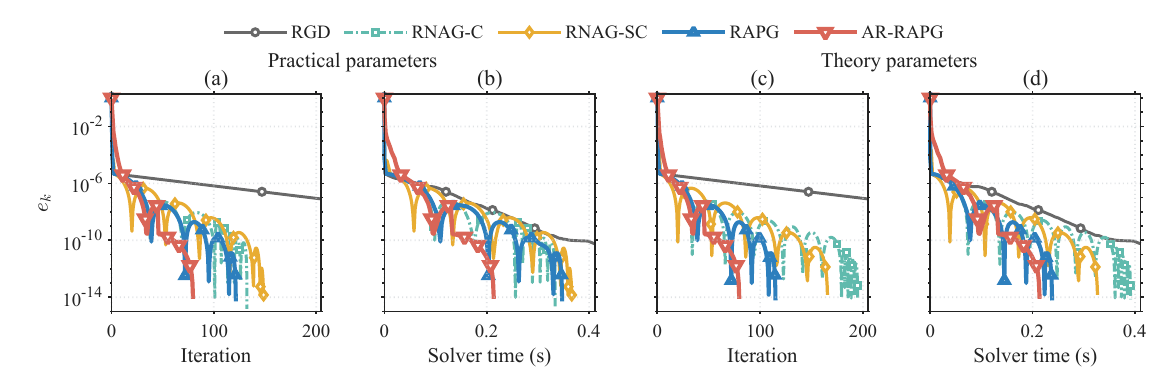}
    \caption{Rayleigh quotient  convergence at $(d,\Delta)=(1000,10)$. RNAG-C, RNAG-SC, and RAPG use practical parameters in (a)--(b) and theory-admissible parameters in (c)--(d); RGD and AR-RAPG results are shared across the comparisons.}
    \label{fig:sphere-main}
\end{figure}

\begin{figure}[!htbp]
    \centering
    \includegraphics[width=\textwidth]{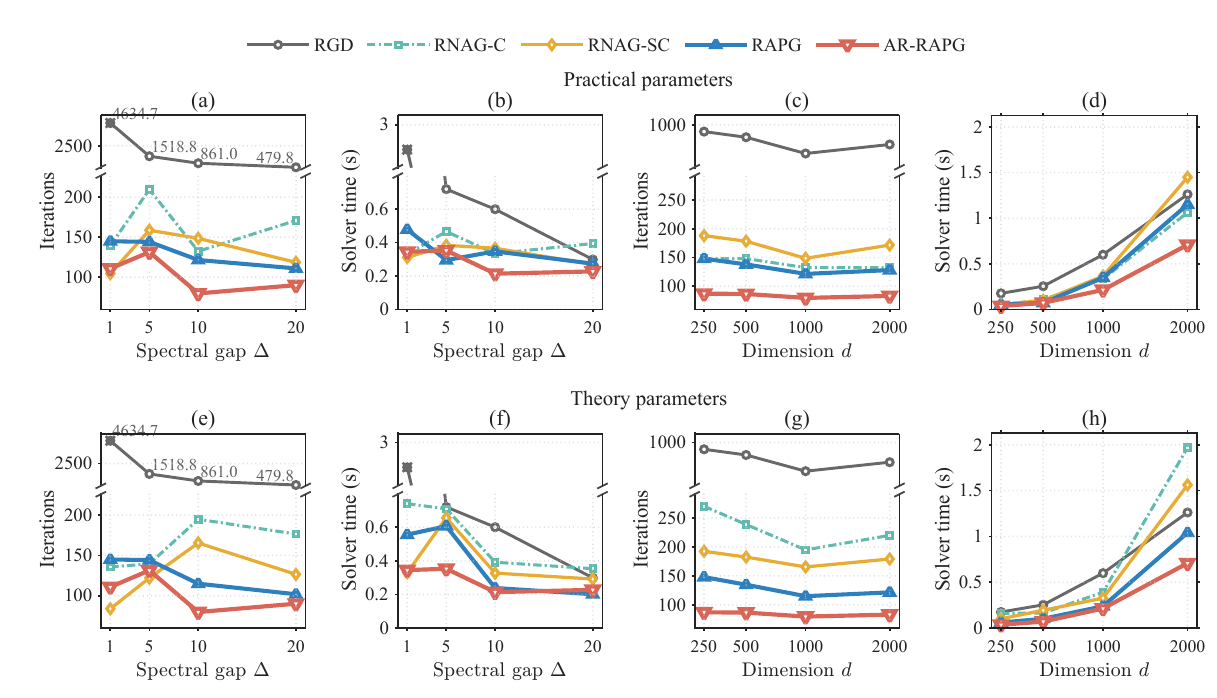}
    \caption{Rayleigh quotient scaling with practical parameters in (a)--(d) and theory-admissible parameters in (e)--(h). Columns 1--2 vary $\Delta$ at $d=1000$; columns 3--4 vary $d$ at $\Delta=10$. RGD and AR-RAPG results are shared across the comparisons. RGD reaches the iteration limit in four of ten runs at $\Delta=1$.}
    \label{fig:sphere-scaling}
\end{figure}

\paragraph{Comparison results.}
Under both parameter regimes, AR-RAPG has the lowest mean iteration count and solver time on the main instance $(d,\Delta)=(1000,10)$ and throughout the dimension sweep, followed by RAPG; see \Cref{fig:sphere-main,fig:sphere-scaling}. In the spectral gap sweep, AR-RAPG is best at larger gaps and RNAG-SC at the smallest gaps. The theory-admissible parameters give $\theta=1$ in all cases, preserving the step size $1/L$. In all runs, the iterates $x_k$ and $y_k$ of RAPG and AR-RAPG remained in $\Omega$.

\section{Conclusions}\label{sec:conclusions}

This paper develops a FISTA-type accelerated proximal gradient method for composite optimization on Riemannian manifolds. The method formulates proximal subproblems in tangent spaces through pullbacks of the nonsmooth term, and employs a parameterization of the acceleration coefficients that extends continuously to the limiting case $\mu-\rho=0$. Under suitable convexity and geometric assumptions, a single potential-function analysis then yields an $\mathcal{O}(k^{-2})$ rate for the objective gap when $\mu-\rho\geqslant0$ and an accelerated linear rate when $\mu-\rho>0$. In the smooth setting, the method reduces to a Riemannian accelerated gradient scheme whose iteration and parameter recurrence extend continuously to $\mu=0$, so that geodesically convex and geodesically strongly convex objectives are treated by the same algorithm and analysis.

The adaptive restart variant complements these results when the convexity assumptions or parameter information required by the acceleration theory are unavailable. It retains the momentum updates of the accelerated method, while its proximal gradient safeguards guarantee that, under compactness of the initial sublevel set, every accumulation point of the safeguard reference sequence is stationary and that an $\varepsilon$-stationary point is obtained within $\mathcal{O}(\varepsilon^{-2})$ iterations. 
Numerical experiments on nonsmooth composite and smooth problems  demonstrate the effectiveness of the proposed methods and illustrate the benefits of adaptive restart in improving robustness and computational performance.

\section*{Acknowledgements}
The authors thank Chengyang Yi for insightful discussions on curvature, Zhifeng Deng for valuable comments on the manuscript, and Max Pfeffer for helpful suggestions on the experiments.

\phantomsection
\addcontentsline{toc}{section}{References}

\begingroup
\makeatletter
\def\@listi{\leftmargin\leftmargini \parsep 0pt \topsep 0pt \itemsep 2pt}
\makeatother
\bibliographystyle{abbrvurl}
\bibliography{refs}
\endgroup

\begin{appendices}

\section{Pullback weak convexity of the ball penalty on hyperbolic space}
\label{sec:Appen:weak-convex-HyperSpace}

Let $\HH$ denote $n$-dimensional hyperbolic space with sectional curvature $-1$. For a smooth curve $\gamma$ in $\HH$, write $\gamma'$ and $\frac{\mathrm D\gamma'}{\dd t}$ for its velocity and covariant acceleration, respectively. We suppress the dependence on $t$ when it is clear from context.
Let $b\in\mathbb H^n$ and $r,\tau\geqslant0$. Set $C=\overline{B}(b,r)$ for $r>0$ and $C=\{b\}$ for $r=0$, and define
\begin{equation}
    h(y)=\tau\operatorname{dist}(y,C)=\tau[\mathrm d(y,b)-r]_+,
    \qquad y\in\mathbb H^n,
    \label{eq:app-ball-penalty}
\end{equation}
where $[a]_+=\max\{a,0\}$.
Although $h$ is geodesically convex, its retraction pullback need not be locally convex, even under the exponential map. 

\begin{proposition}[Weak convexity of retraction pullback]
\label{prop:hyp-ball-retraction}
Let $R$ be a smooth retraction on $\mathbb H^n$, and fix $x\in\mathbb H^n$ and $D>0$ such that $R_x$ is defined on an open neighborhood of $\overline{B}(0_x,D)$. For $\eta,v\in\mathrm T_x\mathbb H^n$, set $z(t)=\eta +tv \in \overline{B}(0_x,D)$ and $\gamma(t)=R_x(z(t))$, and let
\begin{equation}
    \mathcal A_R(x,D)=\max_{\substack{\eta\in\overline{B}(0_x,D)\\\|v\|_x=1}}\left\|\left.\frac{\mathrm D\gamma'(t)}{\mathrm dt}\right|_{t=0}\right\|_{R_x(\eta)}, \quad \text{ and } \quad 
    \rho_R(x,D)=\tau\mathcal A_R(x,D).
    \label{eq:app-retraction-acceleration}
\end{equation}
Then $\mathcal A_R(x,D)$ is finite. For the function $h$ defined in~\eqref{eq:app-ball-penalty}, the pullback $h\circ R_x$ is $\rho_R(x,D)$-weakly convex on $\overline{B}(0_x,D)$.
\end{proposition}

\begin{proof}
The maximum in~\eqref{eq:app-retraction-acceleration} is finite by smoothness and compactness. Since covariant acceleration is quadratic in the direction $v$, every affine segment $z(t)$ in $\overline{B}(0_x,D)$ satisfies
\begin{equation}
    \left\|\frac{\mathrm D\gamma'(t)}{\mathrm dt}\right\|_{\gamma(t)}\leqslant\mathcal A_R(x,D)\|v\|_x^2,
    \qquad \gamma(t)=R_x(z(t)).
    \label{eq:app-retr-acce-bound}
\end{equation}

Let $d_b(y)=\mathrm d(y,b)$ and define $d_{b,\varepsilon}(y)=\sqrt{d_b(y)^2+\varepsilon^2}$ for $\varepsilon>0$. On $\mathbb H^n$, the function $d_b$ is geodesically convex and $1$-Lipschitz, and $d_b^2$ is smooth. The scalar map $t\mapsto\sqrt{t^2+\varepsilon^2}$ is convex, nondecreasing, and $1$-Lipschitz on $[0,\infty)$. Thus $d_{b,\varepsilon}$ is smooth and geodesically convex, with
\begin{equation}
    \operatorname{Hess}d_{b,\varepsilon}\succeq0,
    \qquad \|\operatorname{grad}d_{b,\varepsilon}\|\leqslant 1.
    \label{eq:app-regular}
\end{equation}
Let $\rho=\rho_R(x,D)$ and define
\(
    \psi_\varepsilon(\eta)=\tau d_{b,\varepsilon}(R_x(\eta))+\frac{\rho}{2}\|\eta\|_x^2.
\)
Along any segment $z(t)$ in $\overline{B}(0_x,D)$, the chain rule gives
\begin{align*}
    \frac{\mathrm d^2}{\mathrm dt^2}\psi_\varepsilon(z(t))
    & \stackrel{\hphantom{\eqref{eq:app-regular}}}{=}\tau\bigl\langle\operatorname{Hess}d_{b,\varepsilon}(\gamma(t))[\gamma'],\gamma'\bigr\rangle_{\gamma(t)}
    +\tau\left\langle\operatorname{grad}d_{b,\varepsilon}(\gamma(t)),\frac{\mathrm D\gamma'}{\mathrm dt}\right\rangle_{\gamma(t)}+\rho\|v\|_x^2\\
    & \stackrel{\eqref{eq:app-regular}}{\geqslant}-\tau\left\|\frac{\mathrm D\gamma'}{\mathrm dt}\right\|_{\gamma(t)}+\rho\|v\|_x^2
    \stackrel{\eqref{eq:app-retr-acce-bound}}{\geqslant} 0.
\end{align*}
Hence $\psi_\varepsilon$ is convex on $\overline{B}(0_x,D)$ for every $\varepsilon>0$. Define
    \(
        \psi(\eta)=\tau d_b(R_x(\eta))+\frac{\rho}{2}\|\eta\|_x^2.
    \)
 The bound $0\leqslant d_{b,\varepsilon}(y)-d_b(y)\leqslant\varepsilon$ gives $0\leqslant\psi_\varepsilon(\eta)-\psi(\eta)\leqslant\tau\varepsilon$ for all $\eta\in\overline{B}(0_x,D)$. Hence $\psi_\varepsilon$ converges uniformly to $\psi$ as $\varepsilon\downarrow0$, and convexity passes to the limit.

 By the definition of $h$, we have
    \(
        h(R_x(\eta))+\frac{\rho}{2}\|\eta\|_x^2
        =\max\left\{\psi(\eta)-\tau r,\ \frac{\rho}{2}\|\eta\|_x^2\right\}.
    \)
The right-hand side is convex on $\overline{B}(0_x,D)$ as the pointwise maximum of two convex functions. Thus $h\circ R_x$ is $\rho_R(x,D)$-weakly convex on $\overline{B}(0_x,D)$.
\end{proof}

For the exponential map, the following lemma gives an explicit bound on $\mathcal A_R(x,D)$.
We use the hyperboloid representation of $\HH=\{x\in\mathbb R^{n+1}:\langle x,x\rangle_L=-1,\ x_0>0\}$~\cite[Section~7.6]{boumal2023introduction}. The Riemannian metric at $x$ is the restriction of the Lorentzian bilinear form
\(
    \langle u,v\rangle_L=-u_0v_0+\sum_{j=1}^n u_jv_j
\)
to $\T_x\HH$. For $x,y\in\HH$, parallel transport along the unique geodesic from $x$ to $y$ is given by
\(
    \Gamma_x^y(u)=u+\frac{\langle y,u\rangle_L}{1-\langle x,y\rangle_L}(x+y),
\)
for  $u\in\T_x\HH$.

\begin{lemma}
\label{lem:hyp-exp-acceleration}
Fix $x\in\mathbb H^n$ and $D>0$. Let $z(t)=\eta+tv$, where $\eta,v\in\mathrm T_x\mathbb H^n$, and suppose that $\|z(t)\|_x\leqslant D$ on an interval $I$. Then the curve $\gamma(t)=\operatorname{Exp}_x(z(t))$ satisfies, for all $t\in I$,
\begin{equation}
    \left\|\frac{\mathrm D\gamma'(t)}{\mathrm dt}\right\|_{\gamma(t)}\leqslant\mathcal A(D)\|v\|_x^2,
    \qquad \mathcal A(D)=\frac{\sinh D\cosh D-D}{D^2}.
    \label{eq:app-acceleration-bound}
\end{equation}
\end{lemma}

\begin{proof}
If $z(t_0)=0$ for some $t_0\in I$, then $\gamma(t)=\operatorname{Exp}_x((t-t_0)v)$ is a geodesic (or point $x$ if $v=0$)  and has zero covariant acceleration. Otherwise, suppressing the dependence on $t$, set
\[
    s=\|z\|_x,\qquad e=\frac{z}{s},\qquad \alpha=\langle v,e\rangle_x,\qquad w=v-\alpha e.
\]
Then $w\perp e$ and $\|v\|_x^2=\alpha^2+\|w\|_x^2$.
Differentiating with $t$ gives
\[
    s'=\alpha,\qquad e'=\frac{w}{s},\qquad \alpha'=\frac{\|w\|_x^2}{s},\qquad w'=-\frac{\|w\|_x^2}{s}e-\frac{\alpha}{s}w.
\]
The hyperbolic exponential map yields 
$\gamma= \e_x(se) = \cosh s\,x+\sinh s\,e$.

 Let $\Gamma_x^{\gamma(t)}$ denote parallel transport along the radial geodesic from $x$ to $\gamma(t)$. Set $E(t)=\Gamma_x^{\gamma(t)}(e(t))$ and $W(t)=\Gamma_x^{\gamma(t)}(w(t))$. Suppressing the argument $t$, we obtain from the transport formula
    \begin{equation*}
        E=\sinh s\,x+\cosh s\,e,
        \qquad
        W=w,
    \end{equation*}
 where $W=w$ is understood as an equality in $\R^{n+1}$. Since parallel transport preserves inner products,  we have $\langle E,W\rangle_\gamma=\langle e,w\rangle_x=0$,
 $\|E\|_\gamma=1$, and $\|W\|_\gamma=\|w\|_x$. Differentiating the preceding expressions with $t$ and using $e=\cosh s\,E-\sinh s\,\gamma$, we obtain
\[
    \gamma'=\alpha E+\frac{\sinh s}{s}W,\qquad
    \frac{\mathrm D E}{\mathrm dt}=\frac{\cosh s}{s}W,\qquad
    \frac{\mathrm D W}{\mathrm dt}=-\frac{\|w\|_x^2\cosh s}{s}E-\frac{\alpha}{s}W.
\]
Define 
\[
    \mathcal A(s)=\frac{\sinh s\cosh s-s}{s^2},\qquad
    \mathcal B(s)=\frac{s\cosh s-\sinh s}{s^2},
\]
with continuous extensions $\mathcal A(0)=\mathcal B(0)=0$. 
Since  $\frac{\mathrm d}{\mathrm dt}(\frac{\sinh s}{s})
= s'(t)\mathcal B(s) 
= \alpha\mathcal B(s)$, we obtain 
\begin{equation}
    \frac{\mathrm D\gamma'}{\mathrm dt}=-\mathcal A(s)\|w\|_x^2E+2\mathcal B(s)\alpha W.
    \label{eq:app-acceleration}
\end{equation}

To estimate the coefficients, observe that, for $s>0$,
\[
    \frac{\mathrm d}{\mathrm ds}\bigl(s^2\mathcal B(s)\bigr)=s\sinh s\geqslant0,\qquad
    \frac{\mathrm d}{\mathrm ds}\bigl(s^2[\mathcal A(s)-2\mathcal B(s)]\bigr)=2\sinh s(\sinh s-s)\geqslant 0.
\]
Both quantities being differentiated vanish at zero, so $0\leqslant2\mathcal B(s)\leqslant\mathcal A(s)$. 
Moreover, 
$ \frac{\mathrm d\mathcal A(s)}{\mathrm ds}=\frac{2\cosh s}{s}\mathcal B(s)\geqslant 0$ 
for $s>0$, so $\mathcal A$ is nondecreasing on $[0,D]$. Using the orthogonality of $E$ and $W$ in~\eqref{eq:app-acceleration}, we conclude that
\begin{align*}
    \left\|\frac{\mathrm D\gamma'}{\mathrm dt}\right\|_\gamma^2
    &=\mathcal A(s)^2\|w\|_x^4+4\mathcal B(s)^2\alpha^2\|w\|_x^2\\
    &\leqslant\mathcal A(s)^2\|w\|_x^2\bigl(\|w\|_x^2+\alpha^2\bigr)
    \leqslant\mathcal A(D)^2\|v\|_x^4.
\end{align*}
Taking square roots yields~\eqref{eq:app-acceleration-bound}.
\end{proof}

\begin{corollary}[Weak convexity of the exponential pullback]
\label{prop:hyp-ball-pullback}
Let $h$ be the ball penalty defined in~\eqref{eq:app-ball-penalty}. For every $x\in\mathbb H^n$ and $D>0$, the pullback $h\circ\operatorname{Exp}_x$ is $\rho(D)$-weakly convex on $\overline{B}(0_x,D)$, where
    $\rho(D)=\tau\frac{\sinh D\cosh D-D}{D^2}$.
\end{corollary}

\begin{proof}
Apply \Cref{prop:hyp-ball-retraction} with $R=\operatorname{Exp}$ and use \Cref{lem:hyp-exp-acceleration}.
\end{proof}

\begin{remark}[A convex pullback at the ball center]
\label{rem:hyp-geodesic-vs-pullback}
At the center $b$ of $C=\overline{B}(b,r)$, the exponential pullback is given by
$(h\circ\operatorname{Exp}_b)(\eta)=\tau[\|\eta\|_b-r]_+$.
It is therefore convex on $\T_b\mathbb H^n$. The following example shows that convexity can fail at a base point away from the center.
\end{remark}

\begin{example}[A nonconvex pullback away from the ball center]
\label{ex:hyp-nonconvex-pullback}
Let $n\geqslant2$, $\tau>0$, $r\geqslant0$, and $D>0$. Fix $x\in\mathbb H^n$, choose $s\in(0,D)$ and orthonormal vectors $e_1,e_2\in\T_x\mathbb H^n$, and set
\[
    y=\operatorname{Exp}_x(se_1),\qquad
    T=\frac{\sinh s}{s},\qquad
    u=\Gamma_x^y(e_1+Te_2).
\]
For $\ell>r$, set $b=\operatorname{Exp}_y(\ell u/\|u\|_y)$ and define $h=\tau[d_b-r]_+$, where $d_b=\mathrm d(\,\cdot\,,b)$.
The curve $\gamma(t)=\operatorname{Exp}_x(se_1+t(e_1+e_2))$ satisfies $\gamma(0)=y$ and $\gamma'(0)=u$. Since $d_b(y)=\ell>r$, we have $h=\tau(d_b-r)$ near $y$. The vector $u$ points along the geodesic from $y$ to $b$, so $\grad d_b(y)=-u/\|u\|_y$ and $\langle\Hess d_b(y)[u],u\rangle_y=0$. The chain rule and \eqref{eq:app-acceleration} give
\[
    \left.\frac{\mathrm d^2}{\mathrm dt^2}h(\gamma(t))\right|_{t=0}
    =-\frac{\tau}{\|u\|_y}
    \left\langle u,\left.\frac{\mathrm D\gamma'}{\mathrm dt}\right|_{t=0}\right\rangle_y
    =-\frac{\tau K(s)}{s^3\sqrt{1+T^2}},
\]
where $K(s)=s\sinh s\cosh s+s^2-2\sinh^2s$. Since $K(0)=K'(0)=0$ and $K''(s)=4\sinh s\,(s\cosh s-\sinh s)>0$ for $s>0$, the second derivative is strictly negative. As $se_1$ lies in the interior of $\overline{B}(0_x,D)$, the pullback $h\circ\operatorname{Exp}_x$ is nonconvex on this ball.
\end{example}

\section{Pullback weak convexity of Euclidean convex functions on embedded submanifolds}
\label{sec:Appen:weak-convex}

Let $\m$ be a smooth embedded submanifold of $\R^n$ with the induced Euclidean metric. For $x\in\m$, let $\mathrm P_x$ denote the orthogonal projection from $\R^n$ onto $\T_x\m$. Fix a retraction $R:\mathcal D_R\to\m$. For a function $h$ on $\m$, let $\hat h_x=h\circ R_x$ on $\mathcal D_x$.
Throughout, $\partial$ denotes the Clarke subdifferential. The superscript $\sharp$ denotes the adjoint with respect to the induced inner products.

The following lemma relates the subdifferentials of an ambient function, its restriction to $\m$, and the associated retraction pullbacks.

\begin{lemma}\label{lem:subgradient}
    Let $\bar h:\R^n\to\R$ be locally Lipschitz and let $h=\bar h|_{\m}$. Then, for every $z\in\m$,
\begin{equation}\label{eq:app-pullback-chain-inclusion}
\partial h(z)\subseteq\mathrm P_z[\partial\bar h(z)],
\end{equation}
with equality if $\bar h$ is convex.
Now let $x\in\m$ and $\omega\in\mathcal D_x$, and set $z=R_x(\omega)$. If $R_x$ is a local diffeomorphism at $\omega$, then the retraction pullback $\hat h_x=h\circ R_x$ satisfies
\(
    \partial\hat h_x(\omega)=\mathrm D R_x(\omega)^\sharp[\partial h(z)].
\)
\end{lemma}

\begin{proof}
Viewed as a map into $\R^n$, $R_z$ satisfies $\mathrm D R_z(0_z)[u]=u$ for every $u\in\T_z\m$. Hence, for any $u\in\T_z\m$ and $v\in\R^n$,
\(
\langle \mathrm D R_z(0_z)[u],v\rangle_{\R^n}
=\langle u,\mathrm P_zv\rangle_z,
\)
which shows that $\mathrm D R_z(0_z)^\sharp=\mathrm P_z$. Applying the Clarke chain rule~\cite[Theorem~2.3.10]{Clarke1990optimization} to $\hat h_z=\bar h\circ R_z$ at $0_z$ gives~\eqref{eq:app-pullback-chain-inclusion}. If $\bar h$ is convex, it is Clarke regular, and the inclusion becomes an equality.

By the retraction property and the inverse function theorem, $R_z$ admits a smooth local inverse $R_z^{-1}$ near $z$, with $R_z^{-1}(z)=0_z$ and $\mathrm D(R_z^{-1})(z)=\operatorname{id}_{\T_z\m}$. Since $R_x$ is a local diffeomorphism at $\omega$, so is $\Phi=R_z^{-1}\circ R_x$. We have $\Phi(\omega)=0_z$ and $\mathrm D\Phi(\omega)=\mathrm D R_x(\omega)$.
Applying the Clarke chain rule to $\hat h_x=\hat h_z\circ\Phi$ at $\omega$ and to $\hat h_z=\hat h_x\circ\Phi^{-1}$ at $0_z$ gives, respectively,
    \[
        \begin{aligned}
            \partial\hat h_x(\omega)
            \subseteq \mathrm D\Phi(\omega)^\sharp[\partial\hat h_z(0_z)],\qquad 
            \partial\hat h_z(0_z)
            \subseteq (\mathrm D\Phi(\omega)^{-1})^\sharp[\partial\hat h_x(\omega)].
        \end{aligned}
    \]
Applying $\mathrm D\Phi(\omega)^\sharp$ to the second inclusion gives the reverse of the first. Thus, using $\partial\hat h_z(0_z)=\partial h(z)$, we obtain
    \(
        \partial\hat h_x(\omega)
        =\mathrm D R_x(\omega)^\sharp[\partial h(z)].
    \)
\end{proof}

The next proposition gives a common radius and weak convexity parameter for retraction pullbacks based at points near any fixed point.

\begin{proposition}\label{lem:weakly-convex}
    Let $\bar h:\R^n\to\R$ be convex and set $h=\bar h|_{\m}$. For every $x_0\in\m$, there exist an open neighborhood $U$ of $x_0$, a radius $r>0$, and a parameter $\rho\geqslant0$ such that, for every $x\in U$:
    \begin{enumerate}[(i)]
        \item\label{LB3-1} 
         $B(0_x,r)\subseteq\mathcal D_x$ and $U\subseteq R_x(B(0_x,r))$. The restriction $R_x|_{B(0_x,r)}$ is a diffeomorphism onto its image, which is open in $\m$.
        \item\label{LB3-2} 
        The pullback $\hat h_x=h\circ R_x$ is $\rho$-weakly convex on $B(0_x,r)$.
    \end{enumerate}
\end{proposition}

\begin{proof}
    To prove (i), define $\Psi:\mathcal D_R\to\m\times\m$ by $\Psi(x,\eta)=(x,R_x(\eta))$. In local bundle coordinates, the retraction identities give the invertible differential $\mathrm D\Psi(x_0,0_{x_0})[u,v]=(u,u+v)$. By the inverse function theorem, $\Psi$ maps some open neighborhood $\mathcal W$ of $(x_0,0_{x_0})$ diffeomorphically onto an open neighborhood of $(x_0,x_0)$.
    Choose an open neighborhood $V$ of $x_0$ with compact closure and a radius $r>0$ such that
    \[
        \mathcal K=\{(x,\eta):x\in\overline V,\ \eta\in\overline{B}(0_x,r)\}
        \subseteq\mathcal W.
    \]
     The image under $\Psi$ of $\{(x,\eta):x\in V,\ \eta\in B(0_x,r)\}$ is open and contains $(x_0,x_0)$. We may therefore choose an open neighborhood $U\subseteq V$ of $x_0$ such that $U\times U$ lies in this image. For each $x\in U$, we have $B(0_x,r)\subseteq\mathcal D_x$ by construction. 
     Since $\Psi$ preserves the first component, restricting $\Psi$ to $\{x\}\times B(0_x,r)$ shows that $R_x$ maps $B(0_x,r)$ diffeomorphically onto the open set $R_x(B(0_x,r))$. This set contains $U$ by the choice of $U$, proving (i).

    To prove (ii), regard $R_x$ as a map into $\R^n$. By smoothness and compactness of $\mathcal K$, there exists $C\geqslant0$ such that
    \[
        \|\mathrm D^2R_x(\eta)[v,v]\|
        \leqslant C\|v\|^2,
        \qquad \forall\, (x,\eta)\in\mathcal K,\ v\in\T_x\m.
    \]
    Since $\bar h$ is finite and convex on $\R^n$, its subgradients are uniformly bounded in norm by some $G\geqslant0$ on the compact set $R(\mathcal K)$. 

    Fix $x\in U$, $\eta,\omega\in B(0_x,r)$, and $\hat g\in\partial\hat h_x(\omega)$. Set $z=R_x(\omega)$. By (i) and \Cref{lem:subgradient}, there exists $g\in\partial\bar h(z)$ such that $\hat g=\mathrm D R_x(\omega)^\sharp(\mathrm P_zg)$. Since $\mathrm D R_x(\omega)[\eta-\omega]\in\T_z\m$, we have
    \[
        \langle g,\mathrm D R_x(\omega)[\eta-\omega]\rangle
        =\langle\hat g, \eta-\omega\rangle.
    \]
    Taylor's theorem along the segment from $\omega$ to $\eta$ gives
    \[
        \|R_x(\eta)-R_x(\omega)-\mathrm D R_x(\omega)[\eta-\omega]\|
        \leqslant\frac{C}{2}\|\eta-\omega\|^2.
    \]
    Combining this estimate with the convex subgradient inequality for $\bar h$ yields
    \[
        \begin{aligned}
            \hat h_x(\eta)
            &\geqslant \hat h_x(\omega)
                +\langle g,R_x(\eta)-R_x(\omega)\rangle\\
            &\geqslant \hat h_x(\omega)
                +\langle\hat g,\eta-\omega \rangle -\frac{GC}{2}\|\eta-\omega\|^2.
        \end{aligned}
    \]
    Set $\rho=GC$. Adding $\rho\|\eta\|^2/2$ to both sides of the above inequality and expanding $\|\eta\|^2=\|\omega+(\eta-\omega)\|^2$ gives
    \[
        \hat h_x(\eta)+\frac{\rho}{2}\|\eta\|^2
        \geqslant \hat h_x(\omega)+\frac{\rho}{2}\|\omega\|^2
        +\langle\hat g+\rho\omega,\eta-\omega\rangle,
        \qquad \forall\, \eta, \omega\in B(0_x,r).
    \]
    Hence $\hat h_x+\rho\|\cdot\|^2/2$ is convex on $B(0_x,r)$. Since $G$ and $C$ are independent of $x$, the same $\rho=GC$ works for every $x\in U$.
\end{proof}

The weak convexity conclusion cannot in general be strengthened to convexity. The restriction of the Euclidean $\ell_1$ norm to the sphere illustrates this distinction: its retraction pullbacks are locally weakly convex, whereas neither geodesic convexity nor convexity of the exponential pullbacks holds locally.

\begin{lemma}\label{lem:norm_1}
Let $n\geqslant2$ and equip $\mathbb S^{n-1}$ with the induced Euclidean metric. The function $h(x)=\|x\|_1$ satisfies the following properties.
\begin{enumerate}[(i)]
\item\label{norm1-1} $h$ is strictly geodesically concave on each open orthant of the sphere.
\item\label{norm1-2} For every $x\in\mathbb S^{n-1}$, $h$ is not geodesically convex on any neighborhood of $x$, and hence $h\circ\e_x$ is not convex on any neighborhood of $0_x$.
\item\label{norm1-3} For every retraction $R$, $h\circ R_x$ is weakly convex on a neighborhood of $0_x$.
\end{enumerate}
\end{lemma}

\begin{proof}
    Let $\gamma$ be a unit-speed geodesic in the same open orthant on $\mathbb S^{n-1}$. On any interval where the sign of each coordinate of $\gamma(t)$ is constant, we have $h(\gamma(t))=\langle s,\gamma(t)\rangle$ for a fixed $s\in\{-1,1\}^n$. Since $\gamma''(t)=-\gamma(t)$, it follows that
    \(
        (h\circ\gamma)''(t)=-h(\gamma(t))<0.
    \)
    This proves (i).

    For (ii), it remains to consider points  $x\in\mathbb S^{n-1}$ with at least one zero coordinate.
    Suppose $x_i=0$ for some $i$. Then the $i$th standard basis vector $e_i$ lies in $\T_x\mathbb S^{n-1}$. Along the geodesic $\gamma(t)=\e_x(te_i)=\cos(t)x+\sin(t)e_i$, we have
        $h(\gamma(t))=\cos(t)h(x)+\sin(t)$
     for    $0<t<\frac{\pi}{2}$.
    Hence $(h\circ\gamma)''(t)=-h(\gamma(t))<0$ on this interval. Every neighborhood of $x$ contains a short minimizing subarc of $\gamma$ with positive parameter, so $h$ fails to be geodesically convex on that neighborhood. Likewise, $(h\circ\e_x)(te_i)=h(\gamma(t))$ rules out convexity of $h\circ\e_x$ on any neighborhood of $0_x$.

    Finally, (iii) follows from \Cref{lem:weakly-convex}, because $\ell_1$ norm is convex in $\R^{n}$.
\end{proof}

\section{Sectional curvature of the oblique manifold}
\label{sec:A-curvature}

For integers $n\geqslant2$ and $p\geqslant1$, equip the oblique manifold \(\operatorname{OB}(p,n)=\{X\in\R^{n\times p}:\operatorname{diag}(X^{\mathsf T}X)=\mathbf{1}_p\}\) with the Riemannian metric induced by the Frobenius inner product. With this metric, $\operatorname{OB}(p,n)$ is isometric to the Riemannian product $(\mathbb{S}^{n-1})^p$ and hence has dimension $p(n-1)$. If $p(n-1)<2$, sectional curvature is undefined because the tangent spaces contain no two-dimensional planes.
For $p(n-1)\geqslant2$, let $\operatorname{sec}_X(\Pi)$ denote the sectional curvature of a two-dimensional plane $\Pi\subseteq\T_X\operatorname{OB}(p,n)$. 

\begin{lemma}\label{lem:oblq}
Suppose $p(n-1)\geqslant2$. The sectional curvatures of $\operatorname{OB}(p,n)$ lie in $[0,1]$. The sectional curvature is identically $0$ for $n=2$ and identically $1$ for $p=1$. When $n\geqslant3$ and $p\geqslant2$, both bounds are attained at every point.
\end{lemma}

\begin{proof}
Let $\Pi$ be a two-dimensional tangent plane at $X=(x_1,\ldots,x_p)\in \operatorname{OB}(p,n)$, and choose an orthonormal basis $\{u,v\}$ of $\Pi$. Denote $u=(u_1,\ldots,u_p)$ and $v=(v_1,\ldots,v_p)$, where $u_i,v_i\in\T_{x_i}\mathbb{S}^{n-1}$. The product curvature formula~\cite[Exercise~1.68]{chow2006hamilton} and the curvature tensor of the unit sphere give
\[
    \operatorname{sec}_X(\Pi)=\sum_{i=1}^p\left(\|u_i\|^2\|v_i\|^2-\langle u_i,v_i\rangle^2\right).
\]
By the Cauchy--Schwarz inequality and $\|v_i\|\leqslant\|v\|=1$,
\[
    0\leqslant\operatorname{sec}_X(\Pi)
    \leqslant\sum_{i=1}^p\|u_i\|^2\|v_i\|^2
    \leqslant\sum_{i=1}^p\|u_i\|^2=1.
\]

If $n=2$, the vectors $u_i$ and $v_i$ lie in a one-dimensional space, so every summand vanishes. If $p=1$, orthonormality gives $\operatorname{sec}_X(\Pi)=1$. Finally, when $n\geqslant3$ and $p\geqslant2$, at any $X$ an orthonormal pair supported in a single factor attains the upper bound $1$, and a pair of unit vectors supported in distinct factors attains the lower bound $0$.
\end{proof}

\section{Local strong convexity of the spherical Rayleigh quotient}
\label{sec:app-sphere-rayleigh}

\begin{lemma}
\label{lem:sphere-rayleigh-local-convexity}
Let $A\in\mathbb R^{d\times d}$ with $d\geqslant2$ be symmetric with eigenvalues $\lambda_1>\lambda_2\geqslant\cdots\geqslant\lambda_d$, and let $x_*$ be a unit eigenvector associated with $\lambda_1$. Set $\Delta=\lambda_1-\lambda_2$ and $L=\lambda_1-\lambda_d$. If $0<R_\Omega<\arcsin\sqrt{\Delta/(L+\Delta)}$, then the geodesic ball $\Omega=B(x_*,R_\Omega)\subset\mS^{d-1}$ is uniquely geodesically convex, and the function $f(x)=-\tfrac12\,x^\top Ax$ is geodesically $\mu$-strongly convex on $\Omega$ with $\mu=\Delta-(L+\Delta)\sin^2R_\Omega>0$.
\end{lemma}

\begin{proof}
Since $0<\Delta\leqslant L$, we have $R_\Omega<\pi/4$, so $\Omega$ is uniquely geodesically convex. Fix $x\in\Omega$ and a unit vector $\eta\in\T_x\mS^{d-1}$, and write $\vartheta=\mathrm d(x,x_*)$. We have $\eta^\top x=0$, $x^\top x_*=\cos\vartheta$, and
\(
    |\eta^\top x_*|
    =|\eta^\top(x_*-\cos\vartheta\,x)|
    \leqslant\|x_*-\cos\vartheta\,x\|_2
    =\sin\vartheta.
\)
The spectral bounds therefore give
\[
    \begin{aligned}
        x^\top Ax
        &\geqslant\lambda_1\cos^2\vartheta+\lambda_d\sin^2\vartheta
        =\lambda_1-L\sin^2\vartheta,\\
        \eta^\top A\eta
        &\leqslant \lambda_1 \sin^2\vartheta + \lambda_2 \cos^2\vartheta 
        =
        \lambda_2+\Delta\sin^2\vartheta.
    \end{aligned}
\]
Under the metric induced by the ambient Euclidean inner product, the Riemannian Hessian satisfies
\[
    \begin{aligned}
        \langle\Hess f(x)[\eta],\eta\rangle
        =x^\top Ax-\eta^\top A\eta
        \geqslant\Delta-(L+\Delta)\sin^2\vartheta
        \geqslant\Delta-(L+\Delta)\sin^2R_\Omega
        =\mu.
    \end{aligned}
\]
The radius condition ensures $\mu>0$, which completes the proof.
\end{proof}

\end{appendices}

\end{document}